\documentclass[reqno]{amsart}
\usepackage{amsfonts,amsmath,amssymb,hyperref,mathtools,enumitem}
\usepackage{datetime}
\usepackage[dvipsnames]{xcolor}
\usepackage{tikz}

\usepackage[margin=1.1in]{geometry}
\numberwithin{equation}{section}

\usepackage{hyperref}
\hypersetup{colorlinks=true, pdfstartview=FitV, linkcolor=BrickRed,citecolor=BrickRed, urlcolor=BrickRed, linktoc=page}
\mathtoolsset{showonlyrefs=true}

\usepackage{graphicx}
\graphicspath{ {./} }

\def\eps{\varepsilon}
\def\k{\kappa}

\def\La{\Lambda}

\def\ka{\kappa}

\def\varep{\varepsilon}
\def\vphi{\varphi}

\def\dt{\partial_t}
\def\divergence{\text{div}}
\def\dr{\partial_r}
\def\dz{\partial_z}
\def\done{\partial_1}
\def\dtwo{\partial_2}
\def\drho{\partial_\rho}
\def\dphi{\partial_\phi}
\def\bvphi{\bar{\vphi}}
\def\tvphi{\tilde{\vphi}}

\providecommand{\ip}[1]{\langle#1\rangle}
\providecommand{\abs}[1]{\left\lvert#1\right\rvert}

\providecommand{\norm}[1]{\left\|#1\right\|}

\newcommand{\blue}[1]{#1}

\newcommand{\Z}{\mathbb{Z}}

\newcommand{\R}{\mathbb{R}}
\newcommand{\PP}{\mathbb{P}}

\renewcommand{\div}{{\mathrm{div }}}
\newcommand{\curl}{{\mathrm{curl }}}

\providecommand{\abs}[1]{\left\lvert#1\right\rvert}
\providecommand{\norm}[1]{\left\|#1\right\|}

\newtheorem{theorem}{Theorem}[section]
\newtheorem{lemma}[theorem]{Lemma}
\newtheorem{corollary}[theorem]{Corollary}
\newtheorem{proposition}[theorem]{Proposition}

\newtheorem{remark}[theorem]{Remark}

\begin{document}

\title{Increased lifespan for 3d compressible Euler flows with rotation}

\author{Haram Ko}
\address{Brown University, Providence, RI, USA}
\email{haram\_ko@brown.edu}

\author{Benoit Pausader}
\address{Brown University, Providence, RI, USA}
\email{benoit\_pausader@brown.edu}

\author{Ryo Takada}
\address{The University of Tokyo, Tokyo, Japan}
\email{r-takada@g.ecc.u-tokyo.ac.jp}

\author{Klaus Widmayer}
\address{University of Vienna, Vienna, Austria \& University of Zurich, Zurich, Switzerland}
\email{klaus.widmayer@univie.ac.at \& klaus.widmayer@math.uzh.ch}

\begin{abstract}
We consider the compressible Euler equation with a Coriolis term and prove a lower bound on the time of existence of solutions in terms of the speed of rotation, sound speed and size of the initial data. Along the way, we obtain precise dispersive decay estimates for the linearized equation. In the incompressible limit, this improves current bounds for the incompressible Euler-Coriolis system as well.
\end{abstract}

\maketitle

\setcounter{tocdepth}{1}

\tableofcontents

\section{Introduction}

\subsection{Long time existence for the compressible Euler system with rotation}

A 3-dimensional compressible inviscid fluid can behave quite violently, exhibiting finite time blowup \cite{Sideris85}, and its dynamics has been extensively studied and still remains a field of fruitful works, see \cite{BucShkVic23-b, BucShkVic23-a,  Christo07,LukSpe24, ShkVicol24} and \cite{BucCaoGom-radial, CaoGomSta-nonradial,MerRapRodSze-1, MerRapRodSze-2}. For geophysical fluids, however, it turns out that the rotation of the background plays an important role that can often not be neglected, and the more relevant (isentropic) model is given by:
\begin{equation}\label{eq:CER}
\begin{split}
	\partial_tn+\hbox{div}\left[nu\right]&=0,\\
	\partial_tu+u\cdot\nabla u+\vec{k}\times u+\frac{1}{n}\nabla p&=0.
\end{split}
\end{equation}
Here, $n(x,t)\ge0$ denotes the density, $u(x,t)\in\mathbb{R}^3$ the local velocity, $\vec{k}\in\R^3$ the speed and axis of rotation of the frame and $p=p(n)\ge0$ the pressure which we assume to satisfy a polytropic pressure law
\begin{equation}
p(n)=\frac{1}{\gamma}n^\gamma,\qquad \gamma>1.
\end{equation}

In this article, we study how the stability of the equilibrium $n\equiv n_0,u\equiv 0$ in \eqref{eq:CER} is affected by two important parameters: the Mach number, related to $n_0$, and the Rossby number, related to the speed of rotation $\vert \vec{k}\vert$. We establish the complete linear analysis of the system including the optimal and finer linear decay and the properties of the transformation to dispersive variables. This extends and refines previous analysis in \cite{Mu,NgoScr18}.

To this effect, we fix coordinates such that $\vec{k}=\eps^{-1}\vec{e}_3$, $\eps>0$, and (instead of $n$) work with the ``sound speed'' unknown
\begin{equation}\label{SoundSpeedUnknown}
\sigma=\sigma(n):=\sqrt{p'(n)}=n^\alpha,\qquad \alpha:=\frac{\gamma-1}{2}
\end{equation}
which satisfies
\begin{equation}
\partial_t\sigma+u\cdot\nabla\sigma+\alpha\sigma \div(u)=0,\qquad \frac{1}{n}\nabla p=\frac{1}{\alpha}\sigma\cdot\nabla\sigma.
\end{equation}
To study the stability of the equilibrium $n\equiv n_0,u\equiv 0$ in \eqref{eq:CER}, we write $\sigma=c+\alpha\rho$ with $c=\sqrt{p'(n_0)}$ the ``effective sound speed''. We obtain that \eqref{eq:CER} is then equivalent to
\begin{equation}\label{eq:CER2}
\begin{aligned}
	\partial_t\rho+u\cdot\nabla\rho+(c+\alpha\rho) \div(u)&=0, \\
	\partial_tu+u\cdot\nabla u+\eps^{-1} \vec{e}_3\times u+(c+\alpha\rho)\nabla\rho&=0.
\end{aligned}
\end{equation}

We note that, for $\varepsilon=\infty$, we recover the classical compressible Euler equation. When $c=\infty$, the first equation gives $\hbox{div}(u)=0$, and considering the curl of the second equation, we obtain the Incompressible Euler-Coriolis system \eqref{EC}. We consider the intermediate regime and our main result provides a lower bound on the time of existence of solutions of \eqref{eq:CER2} in terms of the parameters $c$, $\varepsilon$ and the norm of the initial data.

\begin{theorem}\label{main-theorem}
Let $q > 2$. For $m > \frac{7}{2}$ there exists $M=M(m,q)$ such that the solutions to the Cauchy problem to \eqref{eq:CER2} with initial data $(\rho_0, u_0) \in H^m$ exist at least up to time 
\begin{equation}\label{eq:main-assumption-high-reg}
	T\geq M  \eps^{-\frac{1}{q-1}}\min\{1,(c\eps)^{^{\frac{3}{q-1}}}\}\norm{(\rho_0,u_0)}_{H^m}^{-\frac{q}{q-1}}.
\end{equation}
In particular,
\begin{itemize}
	\item if $c\eps \geq 1$, then $T$ satisfies $T\geq M\eps^{-\frac{1}{q-1}}\norm{(\rho_0,u_0)}_{H^m}^{-\frac{q}{q-1}}$,
	\item if $c\eps \leq 1$, then $T$ satisfies $T\geq M \eps^{\frac{2}{q-1}} c^{\frac{3}{q-1}}\Vert (\rho_0,  u_0)\Vert_{H^m}^{-\frac{q}{q-1}}$.
\end{itemize}
When $m \in (\frac{5}{2}, \frac{7}{2})$, the solutions exist at least up to time
\begin{equation}\label{eq:main-assumption-low-reg}
	T \geq M \eps^{-\frac{ m-\frac{5}{2} }{ q - (m-\frac{5}{2}) }} \min\{1,(c\eps)^{\frac{3(m-\frac{5}{2})}{ q - (m-\frac{5}{2})} }\} \norm{(\rho_0,u_0)}_{H^m}^{-\frac{q}{ q - (m-\frac{5}{2}) }}.
\end{equation}
In particular,
\begin{itemize}
	\item if $c\eps \geq 1$, then $T$ satisfies $T \geq M \eps^{-\frac{ m-\frac{5}{2} }{ q - (m-\frac{5}{2}) }} \norm{(\rho_0,u_0)}_{H^m}^{-\frac{q}{ q - (m-\frac{5}{2}) }}$,
	\item if $c\eps \leq 1$, then $T$ satisfies $T \geq M \eps^{\frac{ 2(m-\frac{5}{2}) }{ q - (m-\frac{5}{2}) }} c^{\frac{3(m-\frac{5}{2})}{ q - (m-\frac{5}{2})} } \norm{(\rho_0,u_0)}_{H^m}^{-\frac{q}{ q - (m-\frac{5}{2}) }}$.
\end{itemize}
\end{theorem}

We make a few remarks on the theorem.
\begin{remark}
\begin{enumerate}
\item Theorem \ref{main-theorem} states that any time can be reached for any data, by choosing suitable parameters. Since small-data blowup is prevalent for many hyperbolic systems \cite{Speck-SmallDataShockBook}, this highlights an important effect of rotation.
\item The lower-bounds above all scale as a high power of the norm. For large norm, they can be combined with the bound from local well-posedness, $T\gtrsim \norm{(\rho_0,u_0)}_{H^{5/2+}}^{-1}$.
	
	Adapting an argument of Sideris \cite{Sideris85}, one sees that the latter scaling is optimal for large data, when the nonlinear terms dominate and can lead to finite-time blow-up (see \cite[Section 1.2]{Ko-Thesis} for more details). For smaller initial data, the linear terms dominate, and Theorem \ref{main-theorem} provides extended time of existence which are, however, probably not optimal. In the model case of the purely compressible Euler system, one still has small data/finite time blow-up, but the time of existence is exponential \cite{Alinhac99-1, Alinhac99-2}, while for the other model of incompressible Euler-Coriolis \eqref{EC}, one has global existence \cite{GPW23, RenTian24}.
\item The stabilizing effect of rotation in \eqref{eq:CER} is yet another example of a modification to the (pure) compressible Euler equations that delays the formation of singularity. This is in line with the works on the small data/global existence for the Euler-Maxwell system \cite{GuoPau11,GuoIonPau16} which shows that adding an electrostatic force suppresses the tendency to shock formation for small perturbations in compressible fluids. In light of this, and of the results in the incompressible case \cite{GPW23, RenTian24}, it would be interesting to investigate whether, in the case $\varepsilon=c=1$, one could have a small data/global existence result. \blue{This would require a substantially more involved analysis since the dispersion relation in \eqref{def-Lambdas} is highly degenerate. The improved linear decay in Section \ref{DecayAndStrichartzSec}, providing (non-uniform) decay at the optimal rate $t^{-3/2}$ provides an encouraging first step in this direction.}
	\item Introducing the sound speed unknown in \eqref{SoundSpeedUnknown} leads to simple energy estimates, which remain uniform as we transition from small data to large data. 
	For non-power-like barotropic pressure laws, one needs modified energies which have different scalings for small and large data.
\item The dispersive stabilization effect persists in the presence of viscosity. We refer to \cite {FujWat25} for such a result for the compressible Navier-Stokes-Coriolis system when $c\varepsilon\ge1$. This idea has already been known and used for the {\it incompressible} Navier-Stokes-Coriolis system, see \cite{KLT14a} for a similar long time existence result and \cite{Ko24} for small data-global existence independent of viscosity.
	\item \blue{ The dependence of $M=M(m,q)$ on $m$ and $q$ is technical, to close the bootstrap. It includes constants used in Littlewood-Paley theorem, Sobolev embedding, energy estimates, and the Strichartz estimates. The restriction $q>2$ is due to the lack of endpoint Strichartz estimates. The dependence on $m$ is limited, and we expect that the regularity of the initial data will not affect $M$ (beyond $m>7/2$). }
\end{enumerate}
\end{remark}

For fixed $\varepsilon$, our analysis holds uniformly in $1\le c<\infty$ and in the limit provides a long-time existence result for the {\it incompressible} Euler-Coriolis system
\begin{equation}\label{EC}
\partial_tu+u\cdot\nabla u+\varepsilon^{-1}\vec{e_3}\times u+\nabla p=0,\qquad\hbox{div}(u)=0,
\end{equation}
which improves on previous results \cite{KLT14b,Ryo16}, more precisely, we get:
\begin{theorem}\label{MainThmEC}

Let $q > 2$ and $m > 5/2, m \neq \frac{7}{2}$. There exists $M=M(m,q)$ such that the solutions to \eqref{EC} with initial data $u_0\in H^m$ exist at least up to time 
\begin{equation}\label{eq:main-assumption-EC}
	T\geq M(m,q)\cdot
	\begin{cases}
	\eps^{-\frac{1}{q-1}}\norm{u_0}_{H^m}^{-\frac{q}{q-1}}&\hbox{ if }m>7/2,\\
	\eps^{-\frac{ m-\frac{5}{2} }{ q - (m-\frac{5}{2}) }} \norm{u_0}_{H^m}^{-\frac{q}{ q - (m-\frac{5}{2}) }}&\hbox{ if }5/2<m<7/2.
	\end{cases}
\end{equation}

\end{theorem}

Alternatively, one can fix $T_\ast$ and the initial data $u_0$, and the result guarantees a solution on $[0,T_\ast]$ so long as
\begin{equation*}
\begin{split}
\varepsilon^{-1}\gtrsim 1+T_\ast^{q-1}\Vert u_0\Vert_{H^m}^q,
\end{split}
\end{equation*}
when $m>7/2$, and similarly for $5/2<m<7/2$.

\subsection{Main strategy}

Our method combines energy estimates and dispersion analysis. The energy estimates give a precise blow-up criterion in terms of the norm $B:=\Vert (\nabla\rho,\nabla u)\Vert_{L^1_tL^\infty_x}$, as well as estimate on the growth of ``high-regularity'' norms $H^m$, see \eqref{energy-estimate}. This allows to obtain local well-posedness and a time of existence  $T\sim \norm{(\rho_0,u_0)}_{H^m}^{-1}$ independent of the linear parameters.

The dispersive features of \eqref{eq:CER2} allow to improve these arguments through better bounds on the above quantity $B$. These are obtained via Strichartz estimates in the Duhamel formulation of strong solutions. Since the equations are quasilinear, this entails a loss of derivatives that is compensated by using a higher energy norm coming from the energy estimates. Compared to prior works, two key refinements allow us to obtain significantly improved bounds (see also the discussion below): on the one hand, the use of inhomogeneous Strichartz estimates allows for a longer lifespan of solutions. On the other hand, a bootstrap that combines both the energy norm and $B$ simultaneously (a framework similar to \cite{JurjaKo25}) allows to be more efficient in terms of derivative loss, improving also the regularity requirements.
At high regularity $m>7/2$, when the loss of derivative can be absorbed, this leads to $T\sim \varepsilon^{-1+}\norm{(\rho_0,u_0)}_{H^m}^{-2+}$, while for lower regularity $5/2<m<7/2$, we ``interpolate'' between the Strichartz norm for small frequencies and the energy estimates for high-frequencies, leading to the intermediate scaling \eqref{eq:main-assumption-low-reg}.

Underlying these arguments is a careful study of the linearized dynamics, which leads to the aforementioned Strichartz estimates. While these are partly known since \cite{FujWat25}, we go an extra step and carry out a precise analysis of the linear decay. Here one finds two significantly different, cylindrically symmetric, dispersion relations $\Omega$ and $\Sigma$ (see Section \ref{SecLin} for the definition), which give rise to different dispersive behavior. Taking into account different behaviors in various frequency regimes, which moreover differs between the two, we obtain the full $t^{-\frac{3}{2}}$ decay for both when properly localized. This averages to an overall decay rate $t^{-1}$ for both dispersive modes (see also \cite{FujWat25}), which we show to be sharp (in Proposition \ref{OptimalDecay} in Appendix \ref{appdx-optimal-decay-rate}).

\subsubsection{Comparison to previous works}

In terms of the method, our work extends works on the incompressible Euler equation \cite{Dutrif05, KLT14b, Ryo16} (see also \cite{CDGG,EW15,JW20,Tak21,Cha23}) to the compressible model, in the sense that we analyze the decay of the linearized model and apply Strichartz estimates to close the bootstrap. However, not only are there two distinct dispersion relations instead of one as in the incompressible model, but they are also no longer homogeneous, which requires additional efforts to fully understand the linearized dynamics.

Several prior works have studied the lifespan of solutions to \eqref{eq:CER2} when the parameters $c,\eps$ satisfy the particular relation $c\eps\gtrsim 1$ \cite{NgoScr18, Mu}, or have investigated certain limiting systems \cite{NecTan20, Chaudh22}. Compared to these works, Theorem \ref{main-theorem} improves on the regularity, dependence on parameters and initial data. We simply require the initial data to be in $H^m, m > 5/2$, without restriction on $c$ and $\eps$.
Even in the restricted scenarios previously considered ($c\eps \gtrsim 1, m \geq 4$ in \cite{Mu}), the exponents for Rossby number (resp. for the norm of the initital data) has been improved from $-\frac{1}{7}$ (resp. $-\frac{4}{7}$), to $-1^+$ (resp. $-2^+$). \blue{ On the other hand, \cite{Chaudh22, NecTan20} proved that in certain singular limit of $\eps \to 0$ and for certain class of initial data, the solutions converges to two-dimensional flows. Although the domain in these papers is different from ours, it is interesting to compare as the results show different aspect of fast rotation. }

We note that, while analyzing the compressible \textit{Navier-Stokes} system, \cite{FujWat25} demonstrated the optimal $t^{-1}$ decay rate for one of the two dispersion relations ($\Omega$, see Section \ref{SecLin} for the definition), using perturbative methods. \blue{ The same authors showed \cite{FujWat24} that they can obtain global solutions for {\it unit} viscosity 
in the fast rotation and small Mach number limit $0 \ll \eps^{-1} \ll c$. Obtaining global bounds independent of the viscosity remains an interesting open problem, but we also refer to \cite{Ko24} for a similar analysis of the (incompressible) Navier-Stokes Coriolis system which holds even in the inviscid limit.}


\subsection{Structure of the paper}
In Section \ref{SecEE}, we state and prove the main energy estimate needed for local wellposedness. In Section \ref{SecLin} we analyze the linearized system of \eqref{eq:CER2}, namely, its dispersive structure. Corresponding Strichartz estimates will follow from the dispersive decay of two dispersion relations, whose proof is given in Section \ref{SecDisp}. As a preparation for the proofs, we give a version of stationary phase theorem in Section \ref{SecQuantStatPhi} which takes a general form and is versatile enough to confront two dispersion relations here which behave very differently. We conclude by proving Theorem \ref{main-theorem} and \ref{MainThmEC} in Section \ref{SecMainThm}.

\subsection{Notation}
For some fixed, sufficiently large constant $K > 1$, we will say $x \sim \lambda$ to indicate $\left\vert \frac{x}{\lambda} \right\vert \in (K^{-1}, K)$, and $x \nsim \lambda$ to mean $\left\vert \frac{x}{\lambda} \right\vert \notin (K^{-1}, K)$. Likewise, $x \lesssim \lambda$ and $x \ll \lambda$ mean $|x| \leq K\lambda$ and $|x| \leq K^{-1} \lambda$, respectively, and $x \gtrsim \lambda$ and $x \gg \lambda$ mean $|x| \geq K^{-1} \lambda$ and $|x| \geq K\lambda$, respectively. In these cases, $\lambda$ will always be positive.

\section{Energy estimates}\label{SecEE}
The natural energy associated to solutions to \eqref{eq:CER} is given by
\begin{equation}
\begin{split}
E(t)& := \int_{\mathbb{R}^3}[\frac{n\vert u\vert^2}{2}+\mathfrak{e}(n)]dx,\qquad\mathfrak{e}^{\prime\prime}(x):=\frac{p^\prime(x)}{x},\quad \mathfrak{e}(n_0)=\mathfrak{e}^\prime(n_0)=0,\\
&=\int_{\R^3} \frac{1}{2} [c+\alpha\rho]^{\frac{1}{\alpha}}|u|^2 + \frac{1}{2\alpha(2\alpha+1)} [(c+\alpha\rho)^{2+\frac{1}{\alpha}}-c^{2+\frac{1}{\alpha}}-(2\alpha+1)c^{1+\frac{1}{\alpha}}\rho] dx.
\end{split}
\end{equation}
This is (formally) a conserved quantity, and yields a priori uniform $L^2$ bounds for the unknowns $(\rho,u)$ satisfying \eqref{eq:CER2}, since
\begin{equation}
\frac{1}{2} c^{1/\alpha}\int |u|^2 + \rho^2 dx\le E(t)\equiv E(0)\lesssim\Vert (\rho(0),u(0))\Vert_{L^2_x}^2,
\end{equation}
at least provided that $\vert \alpha\rho\vert\le \frac{c}{2}$. 
However, more relevant for our purposes are the following energy bounds adapted to \eqref{eq:CER2}:
\begin{proposition}\label{Energy-estimate}
Let $s\geq0$. There exists $K = K(s)\geq 0$ such that if  $(\rho,u)\in C_t([0,T],H^s_x(\mathbb{R}^3))$ solve \eqref{eq:CER2}, then for $t\in[0,T]$, there holds that 
\begin{equation}\label{energy-estimate}
	\norm{(\rho,u)(t)}_{H^s} \leq \norm{(\rho,u)(0)}_{H^s} \exp\left(K(s)\int_0^t \norm{ (\nabla\rho,\nabla u)(\tau)}_{L^\infty} d\tau\right).
\end{equation}
\end{proposition}
The energy estimates in Proposition \ref{Energy-estimate} follow using a standard argument involving the following commutator lemma from \cite{KatoPon88,KLT14b}.
\begin{lemma}\label{Hs-commutator-estimate}
For any $s \geq 0$ and $f \in \dot{W}^{1,\infty}(\R^3) \cap H^s(\R^3)$, $g \in L^{\infty}(\R^3) \cap H^{s-1}(\R^3)$,
\begin{equation*}
	|| (1 - \Delta)^{s/2} (fg) - f (1-\Delta)^{s/2}g ||_{L^2} \lesssim ||\nabla f||_{\infty} ||g||_{H^{s-1}} + ||f||_{H^s} ||g||_{\infty}.
\end{equation*}
\end{lemma}

\begin{proof}[Proof of Proposition \ref{Energy-estimate}]
The proof follows the standard energy estimates. Define $W = (\rho,u)^\intercal$, and take the $H^s$ dot product with \eqref{eq:CER2}. The linear terms except the time derivatives cancel out since
\begin{equation*}
	\langle \divergence(cu), \rho \rangle_{H^s}+ \langle c\nabla\rho, u \rangle_{H^s} + \langle \varepsilon^{-1} u \times e_3, u \rangle_{H^s}  = 0,
\end{equation*}
from integration by parts and the property of cross product. Hence, we get
\begin{equation*}
	\frac{1}{2}\frac{d}{dt} \Vert W(t) \Vert_{H^s}^2 = - \langle u\cdot\nabla\rho, \rho \rangle_{H^s} - \langle \alpha\rho\div(u), \rho \rangle_{H^s} - \langle u\cdot\nabla u, u \rangle_{H^s} - \langle \alpha\rho\nabla\rho, u \rangle_{H^s}.
\end{equation*}
We look at the third term first. Using Lemma \ref{Hs-commutator-estimate},
\begin{align*}
	\langle u\cdot\nabla u, u \rangle_{H^s}
	& = \int (u\cdot\nabla) (1 - \Delta)^{\frac{s}{2}} u \cdot (1-\Delta)^{\frac{s}{2}}u dx \\
	& \hspace{2mm} + \int [ (1-\Delta)^{\frac{s}{2}} (u\cdot\nabla)u - u\cdot (1-\Delta)^{\frac{s}{2}}\nabla u] \cdot (1-\Delta)^{\frac{s}{2}}u dx \\
	& = - \frac{1}{2}\int \text{div}(u) \vert (1-\Delta)^{\frac{s}{2}}u \vert^2 dx + \int [ (1-\Delta)^{\frac{s}{2}} (u\cdot\nabla)u - u\cdot (1-\Delta)^{\frac{s}{2}}\nabla u] \cdot (1-\Delta)^{\frac{s}{2}}u dx \\
	& \lesssim \Vert \text{div}(u) \Vert_\infty \Vert u \Vert_{H^s}^2
	+ (\Vert \nabla u \Vert_\infty \Vert \nabla u \Vert_{H^{s-1}} + \Vert u \Vert_{H^s} \Vert \nabla u \Vert_\infty) \Vert u \Vert_{H^s} \\
	& \lesssim \Vert \nabla u \Vert_\infty \Vert u \Vert_{H^s}^2.
\end{align*}
The first term can be dealt similarly:
\begin{align*}
	\langle u\cdot\nabla\rho, \rho \rangle_{H^s} & = \int [(1-\Delta)^{\frac{s}{2}} u\cdot\nabla\rho - u\cdot\nabla(1-\Delta)^{\frac{s}{2}}\rho] (1-\Delta)^{\frac{s}{2}} \rho dx + \int u\cdot \frac{1}{2}\nabla |(1-\Delta)^{\frac{s}{2}}\rho|^2 dx \\
	& \lesssim (||\nabla u||_\infty ||\nabla\rho||_{H^{s-1}} + ||u||_{H^s} ||\nabla\rho||_\infty) ||\rho||_{H^s} + \left\vert \int \text{div}(u) |(1-\Delta)^{\frac{s}{2}}\rho|^2 dx \right\vert \\
	& \lesssim ||\nabla u||_\infty ||\rho||_{H^s}^2 + ||\nabla\rho||_\infty ||\rho||_{H^s} ||u||_{H^s}.
\end{align*}
Lastly, we deal with the remaining two terms together after pulling out $\alpha$:
\begin{align*}
	\langle \rho \text{ div}(u), \rho \rangle_{H^s} + \langle \rho\nabla\rho, u \rangle_{H^s}
	& = \int \left[ (1-\Delta)^{\frac{s}{2}}(\rho\text{div}(u)) - \rho (1-\Delta)^{\frac{s}{2}} \text{div}(u) \right] (1-\Delta)^{\frac{s}{2}}\rho dx \\
	& + \int [(1-\Delta)^{\frac{s}{2}}\text{div}(u)] \rho (1-\Delta)^{\frac{s}{2}}\rho dx + \int \rho \nabla (1-\Delta)^{\frac{s}{2}}\rho \cdot (1-\Delta)^{\frac{s}{2}} u dx \\
	& + \int [(1-\Delta)^{\frac{s}{2}}(\rho\nabla\rho) - \rho (1-\Delta)^{\frac{s}{2}}\nabla\rho] (1-\Delta)^{\frac{s}{2}} u dx.
\end{align*}
Integration by parts reduces the second line to $-\int [(1-\Delta)^{\frac{s}{2}} \rho] \nabla\rho \cdot (1-\Delta)^{\frac{s}{2}} u dx$. Therefore,
\begin{align*}
	|\langle \rho \text{ div}(u), \rho \rangle_{H^s} + \langle \rho\nabla\rho, u \rangle_{H^s}| \lesssim & (||\nabla\rho||_\infty ||\text{div}(u)||_{H^{s-1}} + ||\rho||_{H^s} ||\text{div}(u)||_\infty) ||\rho||_{H^s} \\
	& + ||\nabla\rho||_\infty ||\rho||_{H^s} ||u||_{H^s} + (||\nabla\rho||_\infty ||\nabla\rho||_{H^{s-1}} + ||\rho||_{H^s} ||\nabla\rho||_\infty) ||u||_{H^s} \\
	\lesssim & ||\nabla u||_\infty ||\rho||_{H^s}^2 + ||\nabla\rho||_\infty ||\rho||_{H^s} ||u||_{H^s}.
\end{align*}
Therefore, we have
\begin{equation*}
	\frac{d}{dt} \Vert W(t) \Vert_{H^s}^2 \lesssim ||\nabla u||_\infty (||\rho||_{H^s}^2 + ||u||_{H^s}^2) + ||\nabla\rho||_\infty ||\rho||_{H^s} ||u||_{H^s}
	\leq (||\nabla\rho||_\infty + ||\nabla u||_\infty) (||\rho||_{H^s}^2 + ||u||_{H^s}^2),
\end{equation*}
and hence,
\begin{equation}
	\frac{d}{dt} ||W(t)||_{H^s} \lesssim (||\nabla\rho(t)||_\infty + ||\nabla u(t)||_\infty) W(t).
\end{equation}
From Grownall's inequality, we obtain \eqref{energy-estimate}.
\end{proof}

\section{Linear analysis}\label{SecLin}
In this section we consider solutions to the linearization of \eqref{eq:CER2}, i.e.\
\begin{equation}\label{eq:CER2_lin}
\begin{aligned}
	\partial_t\rho+c\div(u)&=0, \\
	\partial_tu+\eps^{-1} \vec{e}_3\times u+c\nabla\rho&=0.
\end{aligned}
\end{equation}
By rescaling 
\begin{equation}\label{eq:rescale_lin}
(\rho,u)(t,x)\mapsto (\rho,u)(c^{-1}t,x),
\end{equation}
it suffices to consider (for $\kappa=c\eps>0$)
\begin{equation}\label{eq:CER2_lin2}
\begin{aligned}
	\partial_t\rho+\div(u)&=0, \\
	\partial_tu+\kappa^{-1} \vec{e}_3\times u+\nabla\rho&=0.
\end{aligned}
\end{equation} 
\begin{lemma}\label{lem:disp-rels}
The system \eqref{eq:CER2_lin2} is dispersive with dispersion relations $\pm\Sigma(\xi),\pm\Omega(\xi)$, $\xi\in\R^3$, given by
\begin{equation}\label{def-Lambdas}
	\begin{aligned}    
		\Sigma := \frac{d_2 + d_1}{2},& \qquad \Omega:= \frac{d_2 - d_1}{2},\\
		d_1(\xi):=\sqrt{\abs{\xi_h}^2 + (\xi_3-\kappa^{-1})^2},&\qquad d_2(\xi):= \sqrt{\abs{\xi_h}^2 + (\xi_3+\kappa^{-1})^2}.
	\end{aligned}   
\end{equation}
More precisely, the constant coefficient equation \eqref{eq:CER2_lin2} is equivalent to a system 
\begin{equation}\label{eq:CER2_lin3}
	\partial_t(\widehat{\rho},\widehat{u})(\xi,t)=L(\xi)(\widehat{\rho},\widehat{u})(\xi,t),
\end{equation}
where the matrix $L(\xi)\in\R^{4\times 4}$ is diagonalized with purely imaginary eigenvalues $\pm i\La(\xi)$, $\La\in\{\Sigma,\Omega\}$.
\end{lemma}
\begin{proof}
It is convenient (and useful below) to work with the unknown $V := (\rho,|\nabla|^{-1}\div(u),\abs{\nabla}^{-1}\vec{e}_3\cdot\curl(u),u_3)^\intercal\in\R^4$, for which \eqref{eq:CER2_lin2} can equivalently be expressed as 
\begin{equation}\label{eq:CER2_lin4}
	\partial_t\widehat{V}(\xi,t)=M(\xi)\widehat{V}(\xi,t),
\end{equation}
where
\begin{equation*}
	M(\xi):=\begin{pmatrix} 0 & -|\xi| & 0 & 0 \\
		\abs{\xi} & 0 & \kappa^{-1} & 0 \\
		0 & -\kappa^{-1} & 0 & \kappa^{-1}\frac{i\xi_3}{|\xi|} \\
		-i \xi_3 & 0 & 0 & 0 \end{pmatrix}.
\end{equation*}
Its characteristic polynomial is given by
\begin{equation*}
	\tau^4 + (\kappa^{-2} + |\xi|^2)\tau^2 + \kappa^{-2}\xi_3^2,
\end{equation*}
and its zeros thus satisfy
\begin{equation}\label{evalue-squares}
	\tau^2=\frac12\left(-(\kappa^{-2} + |\xi|^2)\pm\sqrt{(\kappa^{-2} + |\xi|^2)^2-4\kappa^{-2}\xi_3^2}\right)=-\frac14\left(d_1^2+d_2^2\mp 2d_1d_2\right).
\end{equation}
The corresponding eigenvectors $E_{\pm\Lambda}$, $\La\in\{\Sigma,\Omega\}$ follow by direction computation:
\begin{equation}\label{eq:def-eigenvectors}
	E_{\tau} = \begin{pmatrix}
		1 \\ -i\tau|\xi|^{-1} \\ \k(\tau^2 - |\xi|^2)|\xi|^{-1} \\ -\tau^{-1}\xi_3
	\end{pmatrix}.
\end{equation}
\end{proof}

\begin{remark}
\begin{enumerate}
	\item For simplicity of notation we have suppressed the dependence of $\Sigma,\Omega$ on $\kappa$ unless explicitly relevant, where we write $\La=\La_\kappa$ for $\La\in\{\Sigma,\Omega\}$. 
	\item As is clear by the above rescaling, \eqref{eq:CER2_lin} is dispersive with dispersion relations $\pm c\Sigma_{c\eps}(\xi),\pm c\Omega_{c\eps}(\xi)$.
	\item Since the functions $d_j(\xi)$, $j=1,2$, are simply the Euclidean distances from $\xi$ to the points $\pm\kappa^{-1}\vec{e}_3$, the level sets of $\Sigma$ resp.\ $\Omega$ have geometric interpretations as ellipsoids, resp.\ hyperboloids. 
	
	\item We highlight some immediate properties of the dispersion relations, which are relevant to our further analysis: $d_j$, $j=1,2$, and thus also $\Sigma,\Omega$ are smooth except at $\pm\kappa^{-1}\vec{e}_3$. Moreover,
	\begin{equation}
		\Sigma(\xi) \geq \kappa^{-1} \geq \abs{\Omega(\xi)}=\textnormal{sign}(\xi_3)\Omega(\xi),\qquad \blue{ \Sigma^2(\xi)=\Omega^2(\xi)\Leftrightarrow \mathcal{D}:=d_1^2d_2^2=0, }
	\end{equation}
	and
	\begin{equation}
		\Sigma(\xi)\Omega(\xi)= \kappa^{-1}\xi_3.
	\end{equation}
	\item The eigenvalue $c\Omega_{c\eps}(\xi)$ is related to the inertial waves due to rotation, whereas $c\Sigma_{c\eps}(\xi)$ is related to the sound waves. Indeed, it follows from \eqref{def-Lambdas} that, when $\kappa\vert\xi\vert\ge1$,
	\begin{equation}\label{ExpansionsSigmaOmega}
	\begin{split}
	\left\vert c\Sigma/(c\vert\xi\vert)-1\right\vert=O\left( (\kappa\vert \xi\vert)^{-2}\right),\qquad \left\vert c\Omega/(\varepsilon^{-1}\xi_3/\vert\xi\vert)-1\right\vert=O\left( (\kappa\vert \xi\vert)^{-2}\right).
	\end{split}
	\end{equation}
\end{enumerate}
\end{remark}

\subsection{Dispersive unknowns and their boundedness}
The dispersive unknowns $U_{\pm\Lambda}$, $\Lambda\in\{\Omega,\Sigma\}$, of \eqref{eq:CER2_lin2} will naturally be defined via amplitudes of the corresponding eigenvectors, suitably normalized. As in the proof of Lemma \ref{lem:disp-rels}, we consider the change of unknowns
\begin{equation}
\mathcal{T}_1: (\rho,u)\mapsto V = (\rho,\alpha,\beta,\gamma) :=(\rho,|\nabla|^{-1}\div(u),\abs{\nabla}^{-1}\vec{e}_3\cdot\curl(u),u_3)^\intercal,
\end{equation}
composed with the change of variables $\mathcal{T}_2$ that diagonalizes the associated equation \eqref{eq:CER2_lin4}, as represented by the matrix of eigenvectors  \eqref{eq:def-eigenvectors} (and a slight abuse of notation)
\begin{gather}
\mathcal{T}_2:V\mapsto\mathcal{F}^{-1}\left(\mathcal{T}_2(\xi)\widehat{V}\right), \label{def-dispersive-variables} \\
\mathcal{T}_2(\xi):=\begin{pmatrix}
	a_\Sigma E_{\Omega} & a_\Sigma E_{-\Omega} & a_\Omega E_{\Sigma} & a_\Omega E_{-\Sigma}
\end{pmatrix}, \quad
a_{\Lambda} = \frac{\k|\xi_h|}{[ (\k^2 \Lambda^2 - 1)^2 + (\k|\xi_h|)^2 ]^{1/2}}, \quad \Lambda \in \{\Sigma, \Omega\}.
\end{gather}
We denote by $\mathcal{T}:=\mathcal{T}_2^{-1}\circ\mathcal{T}_1$ the resulting linear change of unknowns and let
\begin{equation}
(U_\Omega,U_{-\Omega},U_\Sigma,U_{-\Sigma})^\intercal:=\mathcal{F}^{-1}\left(\mathcal{T}(\xi)(\widehat{\rho},\widehat{u})^\intercal\right).
\end{equation}
We define the associated projectors
\begin{equation}\label{eq:projections}
\PP_{\pm\La}(\rho,u):=U_{\pm\La},\qquad \La\in\{\Sigma,\Omega\}.
\end{equation}
We record that if $(\rho,u)$ satisfy \eqref{eq:CER2_lin2}, then by construction $U_{\pm\La}$ diagonalize the linear evolution, i.e.\
\begin{equation}
\partial_tU_{\mu\La}=i\mu \La U_{\mu\La},\qquad \mu\in\{-,+\}.
\end{equation}
\begin{proposition}\label{prop:Lp-equivalence}
The linear change of unknowns $\mathcal{T}$ is well-defined and invertible. Moreover, $\mathcal{T}$ is $L^p$-bounded for any $ p\in(1, \infty)$, i.e.\ there exist $c_p,C_p>0$ such that 
\begin{equation}
	\norm{\mathbb{P}_{\pm\La}(\rho,u)}_{L^p}=\norm{U_{\pm\La}}_{L^p}\leq c_p\norm{(\rho,u)}_{L^p}\leq C_p\norm{(U_\Omega, U_{-\Omega}, U_\Sigma, U_{-\Sigma})}_{L^p},\qquad \La\in\{\Sigma,\Omega\}.
\end{equation}
In particular, the constants involved in the norm equivalence are independent of $\kappa$.
\end{proposition}
We remark that since $\mathcal{T}$ is linear, its boundedness also extends to $L^p$-based Sobolev spaces $1<p<+\infty$. Moreover, it can be shown that it is an isometry on $L^2$, namely that
\begin{equation}
\norm{(\rho,u)}_{L^2}=2\norm{(U_\Omega, U_{-\Omega}, U_\Sigma, U_{-\Sigma})}_{L^2},
\end{equation}
but this will not be of particular relevance here.\\

To prove Proposition \ref{prop:Lp-equivalence} it will be convenient to use the following concise form.

\begin{lemma}\label{transformation-in-angles}
There holds that
\begin{gather}
	\begin{pmatrix}
		\rho \\ \frac{|\nabla|}{|\nabla_h|}\beta
	\end{pmatrix} = 2
	\begin{pmatrix}
		\cos\theta_1 & -\sin\theta_1 \\
		\sin\theta_1 & \cos\theta_1
	\end{pmatrix}
	\begin{pmatrix}
		\frac{U_{\Omega} + U_{-\Omega}}{2} \\ \Re(U_{\Sigma})
	\end{pmatrix}, \qquad \theta_1 \in [-\frac{\pi}{2},0], \label{From-U12-to-rho-beta} \\
	\begin{pmatrix}
		\frac{|\nabla|\alpha - \partial_3 \gamma}{|\nabla_h|} \\ \gamma
	\end{pmatrix} = 2
	\begin{pmatrix}
		i\sin\theta_2 & \cos\theta_2 \\
		-\cos\theta_2 & -i\sin\theta_2
	\end{pmatrix}
	\begin{pmatrix}
		\frac{U_{\Omega} - U_{-\Omega}}{2} \\ \Im(U_{\Sigma})
	\end{pmatrix}, \quad \theta_2 \in [-\frac{\pi}{2}, \frac{\pi}{2}].\label{From-U12-to-alpha-gamma}
\end{gather}
$\cos\theta_j = \cos(\theta_j(\xi)), \sin\theta_j = \sin(\theta_j(\xi)), j=1,2$ are smooth except at the singular points $(0, \pm\k^{-1})$.
\end{lemma}
\begin{proof}
Writing the transformation \eqref{def-dispersive-variables} in detail, it follows that
\begin{equation}\label{U12-to-rho_a_b_c}
	\begin{split}
		\rho & = \frac{\kappa|\nabla_h|}{[ (\kappa^2 \Sigma^2 - 1)^2 + (\kappa|\nabla_h|)^2 ]^{1/2}} (U_{\Omega} + U_{-\Omega}) + \frac{\kappa|\nabla_h|}{[ (\kappa^2 \Omega^2 - 1)^2 + (\kappa|\nabla_h|)^2 ]^{1/2}} \cdot 2\Re(U_{\Sigma}), \\
		\alpha & = \frac{\kappa|\nabla_h|/|\nabla| }{[ (\kappa^2 \Sigma^2 - 1)^2 + (\kappa|\nabla_h|)^2 ]^{1/2}} \cdot -i\Omega (U_{\Omega} - U_{-\Omega}) + \frac{\kappa|\nabla_h|/|\nabla| }{[ (\kappa^2 \Omega^2 - 1)^2 + (\kappa|\nabla_h|)^2 ]^{1/2}} \Sigma \cdot 2\Im(U_{\Sigma}), \\
		\frac{|\nabla|}{|\nabla_h|}\beta & = \frac{-(\kappa^2 \Sigma^2 - 1)}{[ (\kappa^2 \Sigma^2 - 1)^2 + (\kappa|\nabla_h|)^2 ]^{1/2}} (U_{\Omega} + U_{-\Omega}) + \frac{-(\kappa^2 \Omega^2 - 1)}{[ (\kappa^2 \Omega^2 - 1)^2 + (\kappa|\nabla_h|)^2 ]^{1/2}} \cdot 2\Re(U_{\Sigma}), \\
		\gamma & = \frac{-\kappa|\nabla_h|}{[ (\kappa^2 \Sigma^2 - 1)^2 + (\kappa|\nabla_h|)^2 ]^{1/2}} \kappa\Sigma (U_{\Omega} - U_{-\Omega}) + \frac{-\kappa|\nabla_h|}{[ (\kappa^2 \Omega^2 - 1)^2 + (\kappa|\nabla_h|)^2 ]^{1/2}} i\kappa\Omega \cdot 2\Im(U_{\Sigma}).
	\end{split}
\end{equation}
We write the transformation concisely as
\begin{equation}\label{UtoV-transformation-matrix}
	\begin{pmatrix}
		\rho \\ \frac{|\nabla|\alpha - \partial_3 \gamma}{|\nabla_h|} \\ \frac{|\nabla|}{|\nabla_h|}\beta \\ \gamma
	\end{pmatrix} = 2
	\begin{pmatrix}
		b_{\rho,1} & & b_{\rho,2} & \\
		& b_{\alpha,1} & & b_{\alpha,2} \\
		b_{\beta,1} & & b_{\beta,2} & \\
		& b_{\gamma,1} & & b_{\gamma,2}
	\end{pmatrix}
	\begin{pmatrix}
		\frac{U_{\Omega}+U_{-\Omega}}{2} \\ \frac{U_{\Omega}-U_{-\Omega}}{2} \\ \Re(U_{\Sigma}) \\ \Im(U_{\Sigma})
	\end{pmatrix}.
\end{equation}
where $b_{i,j}$'s are Fourier multipliers. By abusing the notations where we use $b_{i,j}$ for $b_{i,j}(\xi)$, it is immediate that $b_{\rho,j}^2 + b_{\beta,j}^2 \equiv 1$, $j=1,2$. Putting $r=|\xi_h|=(\xi_1^2+\xi_2^2)^{1/2}, z=\xi_3$, from
\begin{equation*}
	(\kappa^2 \Sigma^2 - 1)^2 + \k^2 r^2 = \frac{1}{4}(\k^2|\xi|^2 - 1 + \kappa^2\sqrt{\mathcal{D}})^2 + \k^2 r^2 = \frac{1}{2}(\kappa^4 \mathcal{D} + (\k^2|\xi|^2 - 1)\kappa^2 \sqrt{\mathcal{D}}) = \frac{\kappa^2 \sqrt{\mathcal{D}}}{2}(\kappa^2 \sqrt{\mathcal{D}} + \k^2|\xi|^2 - 1),
\end{equation*}
and similarly $(\kappa^2 \Omega^2 - 1)^2 + \k^2 r^2 = \frac{\kappa^2 \sqrt{\mathcal{D}}}{2}(\kappa^2 \sqrt{\mathcal{D}} - (\k^2|\xi|^2 - 1))$, so we have
\begin{gather*}
	(\kappa^2 \Sigma^2 - 1)^2 + \k^2 r^2 + (\kappa^2 \Omega^2 - 1)^2 + \k^2 r^2 = \kappa^4 \mathcal{D}, \\
	[(\kappa^2 \Sigma^2 - 1)^2 + \k^2 r^2] \cdot [(\kappa^2 \Omega^2 - 1)^2 + \k^2 r^2] = \frac{\kappa^4 \mathcal{D}}{4}(\kappa^4 \mathcal{D} - (\k^2|\xi|^2 - 1)^2) = \kappa^4 \mathcal{D} \cdot \k^2 r^2.
\end{gather*}
Therefore,
\begin{equation*}
	b_{\rho,1}^2 + b_{\rho,2}^2
	= \frac{\k^2 r^2}{(\kappa^2 \Sigma^2 - 1)^2 + \k^2 r^2} + \frac{\k^2r^2}{(\kappa^2 \Omega^2 - 1)^2 + \k^2 r^2}
	= \frac{\k^2 r^2 \cdot \kappa^4 \mathcal{D}}{\kappa^4 \mathcal{D} \cdot \k^2 r^2} = 1.
\end{equation*}
$b_{\beta,1}^2 + b_{\beta,2}^2 \equiv 1$ follows immediately from the above equalities.

We now show similar equalities for $b_{\alpha,j}$ and $b_{\gamma,j}$'s. First,
\begin{align*}
	|\xi_h|^2\left( \left\vert \frac{-i \Omega +i\k\xi_3\Sigma}{|\xi_h|} \right\vert^2 + \left[-\k\Sigma\right]^2 \right) & = [\Omega - \k\xi_3\Sigma]^2 + \k^2 |\xi_h|^2 \Sigma^2 \\
	& = \Omega^2 - 2\k\xi_3\Omega\Sigma + \k^2 |\xi|^2 \Sigma^2 \\
	& = \Omega^2 - 2\k\xi_3\Omega \Sigma + [\kappa^2 \Omega^2 + \kappa^2 \Sigma^2 - 1]\Sigma^2 \\
	& = \kappa^2 \Sigma^4 + \kappa^2 \Omega^2 \Sigma^2 - \Sigma^2 + \Omega^2 - 2\k\xi_3 \Omega \Sigma \hspace{5mm} \text{(using } \Omega \Sigma = \kappa^{-1}\xi_3,) \\
	& = \kappa^2 \Sigma^4 - \Sigma^2 + \Omega^2 - \xi_3^2 \hspace{10mm} \text{(using } \Omega^2 = \kappa^{-2} + |\xi|^2 - \Sigma^2,) \\
	& = \kappa^2 \Sigma^4 - 2\Sigma^2 + \kappa^{-2} + |\xi|^2 - \xi_3^2 \\
	& = \kappa^{-2} ((\kappa^2\Sigma^2 - 1)^2 + \k^2|\xi_h|^2) \\
	& = |\xi_h|^2 \frac{(\kappa^2\Sigma^2 - 1)^2 + (\k|\xi_h|)^2}{(\k|\xi_h|)^2},
\end{align*}
and the other one can be computed symmetrically by exchanging the roles of $\Omega$ and $\Sigma$. These prove that $|b_{\alpha,j}|^2 + |b_{\gamma,j}|^2 \equiv 1$. Moreover, since $\kappa^2 \Sigma^2 = \frac{1}{2}(\k^2 |\xi|^2 + 1 + \kappa^2 \sqrt{\mathcal{D}})$, $\kappa^2 \Omega^2 = \frac{1}{2}(\k^2 |\xi|^2 + 1 - \kappa^2 \sqrt{\mathcal{D}})$,
\begin{align*}
	|b_{\gamma,1}|^2 + |b_{\gamma,2}|^2
	& = \frac{\k^2 r^2}{\kappa^4 \mathcal{D} \cdot \k^2 r^2} ([(\kappa^2 \Omega^2 - 1)^2 + \k^2 r^2] \kappa^2 \Sigma^2 + [(\kappa^2 \Sigma^2 - 1)^2 + \k^2 r^2] \kappa^2 \Omega^2) \\
	& = \frac{1}{\kappa^4 \mathcal{D}} [\frac{\kappa^2 \sqrt{\mathcal{D}}}{2}(\kappa^2 \sqrt{\mathcal{D}} - (\k^2|\xi|^2 - 1)) \cdot \frac{1}{2}(\k^2 |\xi|^2 + 1 + \kappa^2 \sqrt{\mathcal{D}}) \\
	& \hspace{10mm} + \frac{\kappa^2 \sqrt{\mathcal{D}}}{2} (\kappa^2 \sqrt{\mathcal{D}} + \k^2|\xi|^2 - 1) \cdot \frac{1}{2}(\k^2 |\xi|^2 + 1 - \kappa^2 \sqrt{\mathcal{D}})] \\
	& = \frac{1}{4\kappa^2 \sqrt{\mathcal{D}}} [(\kappa^2\sqrt{\mathcal{D}} + 1)^2 - (\k^2|\xi|^2)^2 + (\k^2|\xi|^2)^2 - (\kappa^2\sqrt{\mathcal{D}} - 1)^2] = 1.
\end{align*}
$|b_{\alpha,1}|^2 + |b_{\alpha,2}|^2 \equiv 1$ follows directly. Finally, note that $b_{\gamma,1} = -\kappa\Sigma \cdot b_{\rho,1}$ and $b_{\gamma,2} = -i\kappa\Omega \cdot b_{\rho,2}$.(see \eqref{U12-to-rho_a_b_c}) Therefore, to justify \eqref{From-U12-to-rho-beta} and \eqref{From-U12-to-alpha-gamma}, we only need to check the signs. $b_{\beta,1} \leq 0 \leq b_{\beta,2}$ follows from the fact that $\kappa\Sigma \geq 1 \geq \kappa\Omega \geq -1$, and $b_{\rho,1}, b_{\rho,2}$ are obviously nonnegative. Hence, \eqref{From-U12-to-rho-beta} follows.

\eqref{From-U12-to-alpha-gamma} is proved similarly. From
\begin{gather*}
	(\kappa\Omega - \k^2 z\Sigma) \cdot \kappa\Sigma = \k z(1-\kappa^2\Sigma^2), \\
	\Sigma^2 \geq z^2 = \k z \cdot \Omega \cdot \Sigma,
\end{gather*}
$\Sigma - \k z\Omega$ is always positive, and $\Omega - \k z\Sigma$ has the opposite sign to the sign of $z$. Hence, $b_{\alpha,2}$ is always positive, but $b_{\alpha,1}$ is purely imaginary with the imaginary part having the same sign as the sign of $z$. Also, $b_{\gamma,1}$ is clearly always negative, and $b_{\gamma,2}$ is purely imaginary with the imaginary part having the opposite sign to the sign of $z$. This proves \eqref{From-U12-to-alpha-gamma}.
\end{proof}

\begin{proof}[Proof of Proposition \ref{prop:Lp-equivalence}]
The invertibility and boundedness of change of variable from $(\rho,u)$ to the left-hand side of \eqref{UtoV-transformation-matrix}, $(\rho, \frac{|\nabla|\alpha - \partial_3 \gamma}{|\nabla_h|}, \frac{|\nabla|}{|\nabla_h|}\beta, \gamma)$, follows from direct computations and the boundedness of Riesz transforms, since it is given by
$(\rho,u)\mapsto \mathcal{F}^{-1}\left(\mathcal{T}_1(\xi)(\widehat{\rho},\widehat{u})\right)$, where
\begin{equation}
	\mathcal{T}_1(\xi)=\begin{pmatrix} 1 &0  & 0 & 0 \\
		0 &i\xi_1\abs{\xi_h}^{-1} & i\xi_2\abs{\xi_h}^{-1} & 0 \\ 
		0 & -i\xi_2\abs{\xi_h}^{-1}  & i\xi_1\abs{\xi_h}^{-1}  & 0 \\
		0& 0 & 0 & 1 \end{pmatrix}, \quad \xi_h = (\xi_1,\xi_2).
\end{equation}
To prove Proposition \ref{prop:Lp-equivalence} it thus suffices to study matrices in \eqref{From-U12-to-rho-beta} and \eqref{From-U12-to-alpha-gamma}. We can further reduce the problem into showing that that $b_{\rho,1} = \cos\theta_1, b_{\rho,2} = -\sin\theta_1, b_{\gamma,1} = -\cos\theta_2$, and $b_{\gamma,2} = -i\sin\theta_2$ are $L^p$ Fourier multipliers, because the inverse transformations are again rotation matrices using the same angles, so that this will prove the norm equivalence. See Appendix \ref{appdx-Lp-equiv} for the rest of the proof.
\end{proof}

\subsection{Decay and Strichartz estimates}\label{DecayAndStrichartzSec}

The dispersion relations \eqref{def-Lambdas} are degenerate, and as such do not yield amplitude decay at the full $3d$ rate $t^{-\frac{3}{2}}$. 
The degeneracy is quantified by the following lemma.
\begin{lemma}\label{Hess}
There holds that
\begin{equation}\label{DetHess}
	\det \nabla^2\Lambda = \kappa^{-2} \frac{r^2\Lambda}{(d_1d_2)^4},\quad \Lambda \in \{\Sigma, \Omega\}.
\end{equation}

\end{lemma}

\begin{proof}[Proof of Lemma \ref{Hess}]

We start from the formulas:
\begin{equation*}
\begin{split}
	d_1\nabla d_1=re_r+(z-\kappa^{-1})e_z,\quad d_2\nabla d_2=re_r+(z+\kappa^{-1})e_z,\qquad\nabla^2d=\frac{1}{d}\left[Id-\nabla d\otimes\nabla d\right],
\end{split}
\end{equation*}
and we can evaluate on the orthonormal basis $(e_r,e_\theta,e_z)$:
\begin{equation*}
\begin{split}
	d_1\nabla^2d_1(e_r,e_r)&=1-r^2/d_1^2=\frac{(z-\kappa^{-1})^2}{d_1^2},\qquad d_1\nabla^2d_1(e_r,e_z)=-\frac{r(z-\kappa^{-1})}{d_1^2}\\
	d_1\nabla^2d_1(e_r,e_\theta)&=0=d_1\nabla^2d_1(e_z,e_\theta),\qquad d_1\nabla^2d_1(e_\theta,e_\theta)=1\\
	d_1\nabla^2d_1(e_z,e_z)&=1-(z-\kappa^{-1})^2/d_1^2=\frac{r^2}{d_1^2}.
\end{split}
\end{equation*}
As a consequence, we can write the matrix
\begin{equation*}
\begin{split}
	\nabla^2\Omega&=\frac{1}{2}\begin{pmatrix}\frac{(z+\kappa^{-1})^2}{d_2^3}-\frac{(z-\kappa^{-1})^2}{d_1^3}&0&-r\left(\frac{z+\kappa^{-1}}{d_2^3}-\frac{z-\kappa^{-1}}{d_1^3}\right)\\
		0&\frac{1}{d_2}-\frac{1}{d_1}&0\\
		-r\left(\frac{z+\kappa^{-1}}{d_2^3}-\frac{z-\kappa^{-1}}{d_1^3}\right)&0&r^2\left(\frac{1}{d_2^3}-\frac{1}{d_1^3}\right)
	\end{pmatrix}.
\end{split}
\end{equation*}
and similarly for $\Sigma$. \eqref{DetHess} is now immediate.
\end{proof}
These show that $\nabla^2\Sigma$ vanishes on the vertical axis, and $\nabla^2\Omega$ vanishes additionally on the horizontal plane. In order to track these degeneracies, we introduce the following localizations in the spirit of the standard Littlewood-Paley decomposition. For a standard bump function $\psi \in C^{\infty}(\mathbb{R};[0,1])$ supported on $[-2,2]$, put $\varphi(x) = \psi(\frac{x}{2}) - \psi(x)$. We define the localizations $P_k, P_{k,p}$, and $P_{k,p,q}$ as
\begin{equation*}
\mathcal{F}(P_k f)(\xi) = \varphi_k(\xi)\hat{f}(\xi), \quad \mathcal{F}(P_{k,p} f)(\xi) = \varphi_{k,p}(\xi)\hat{f}(\xi), \quad \mathcal{F}(P_{k,p,q} f)(\xi) = \varphi_{k,p,q}(\xi)\hat{f}(\xi),
\end{equation*}
for $k \in \mathbb{Z}$, $p,q \in \mathbb{Z}_- = \{n\in\mathbb{Z}, n\leq0\}$ where
\begin{equation*}
\varphi_k(\xi) = \varphi(2^{-k} |\xi|), \quad
\varphi_{k,p}(\xi) = \varphi_k(\xi) \varphi(2^{-k-p} |\xi_h|), \quad
\varphi_{k,p,q}(\xi) = \varphi_{k,p}(\xi) \varphi(2^{-k-q} \xi_3).
\end{equation*}
$P_k$ is the usual Littlewood-Paley projection onto functions with frequency of size roughly $2^k$, whereas the localization parameters $p, q$ measure the relative size of horizontal and vertical compoments of the frequencies. To give a complete picture in the medium frequency, we need an additional parameter $l$ to define $P_{k,p;l}$ and $P_{k,p,q;l}$ as
\begin{equation*}
	\mathcal{F}(P_{k,p;l} f)(\xi) = \varphi_{k,p;l}(\xi)\hat{f}(\xi), \quad
	\mathcal{F}(P_{k,p,q;l} f)(\xi) = \varphi_{k,p,q;l}(\xi)\hat{f}(\xi),
\end{equation*}
for $k,p,q$ as above and $p \leq l \leq 0$ where
\begin{equation*}
	\varphi_{k,p;l}(\xi) = \varphi_{k,p}(\xi) \varphi(2^{-k-l} (\xi_3-\k^{-1})), \quad
	\varphi_{k,p,q;l}(\xi) = \varphi_{k,p,q}(\xi) \varphi(2^{-k-l} (\xi_3-\k^{-1})).
\end{equation*}
Lastly, we define $P_{k, p; \leq p}$ and $P_{k,p,q; \leq p}$ by
\begin{equation*}
	\mathcal{F}(P_{k,p; \leq p} f) = \varphi_{k,p; \leq p}(\xi) \hat{f}(\xi), \quad
	\mathcal{F}(P_{k,p,q; \leq p} f)(\xi) = \varphi_{k,p,q; \leq p}(\xi) \hat{f}(\xi),
\end{equation*}
for $k,p,q$ as above where
\begin{equation*}
	\varphi_{k,p; \leq p}(\xi) = \varphi_{k,p}(\xi) \psi(2^{-k-p} (\xi_3-\k^{-1})), \quad
	\varphi_{k,p,q; \leq p}(\xi) = \varphi_{k,p,q}(\xi) \psi(2^{-k-p} (\xi_3-\k^{-1})).
\end{equation*}

A key step in our analysis are the following fine linear decay estimates:

\begin{proposition}\label{prop:full_lin_decay}
There exists $C>0$ such that
\begin{equation}\label{eq:Sigma_decay_full}
	\Vert e^{it\Sigma_\kappa} P_{k,p} f \Vert_{L^\infty} \leq C t^{-\frac{3}{2}} \Vert P_{k,p} f \Vert_{L^1} \begin{cases}
		\k 2^{\frac{5}{2}k-p} \qquad & \k 2^k \gg 1,\\
		2^{\frac{3}{2}k - p} & \k 2^k \sim 1, \\
		\k^{-\frac{5}{2}} 2^{-k-p} \qquad & \k 2^k \ll 1,
	\end{cases}
\end{equation}
\begin{equation}\label{eq:Omega_decay_full}
	\Vert e^{it\Omega_\kappa} P_{k,p,q} f \Vert_{L^\infty} \leq C t^{-\frac{3}{2}} \Vert P_{k,p,q} f \Vert_{L^1} \begin{cases}
		\k^{\frac{3}{2}} 2^{3k-p-\frac{q}{2}} \qquad & \k 2^k \gg 1,\\
		2^{\frac{3}{2}k - p - \frac{q}{2}} & \k 2^k \sim 1, \\
		\k^{-3} 2^{-\frac{3}{2}k-p-\frac{q}{2}} \qquad & \k 2^k \ll 1.
	\end{cases}
\end{equation}
For the medium frequencies, we have finer (in terms of coefficients) decay,
\begin{equation}
	\begin{split}
		\Vert e^{it\Sigma_\kappa} P_{k,p;l} f \Vert_{L^\infty} & \lesssim 2^{\frac{3}{2}k - p +2l} t^{-\frac{3}{2}} \Vert P_{k,p;l} f \Vert_{L^1}, \qquad \k 2^k \sim 1, p \leq l \leq 0, \\
		\Vert e^{it\Sigma_\kappa} P_{k, p; \leq p} f \Vert_{L^\infty} & \lesssim 2^{\frac{3}{2}k + p} t^{-\frac{3}{2}} \Vert P_{k,p; \leq p} f \Vert_{L^1},
	\end{split}
\end{equation}
\begin{equation}
	\begin{split}
		\Vert e^{it\Omega_\kappa} P_{k,p,q;l} f \Vert_{L^\infty} & \lesssim 2^{\frac{3}{2}k - p + 2l - \frac{q}{2}} t^{-\frac{3}{2}} \Vert P_{k,p,q;l} f \Vert_{L^1}, \quad \k 2^k \sim 1, p \leq l \leq 0, \\
		\Vert e^{it\Omega_\kappa} P_{k,p,q; \leq p} f \Vert_{L^\infty} & \lesssim 2^{\frac{3}{2}k + p - \frac{q}{2}} t^{-\frac{3}{2}} \Vert P_{k,p,q; \leq p} f \Vert_{L^1}.
	\end{split}
\end{equation}
\end{proposition}

The proof is postponed to Section 5. Although additional localization in $p,q$ and $l$ is needed for the above stronger decay, in the proof of Theorem \ref{main-theorem} the following decay for functions with $P_k$ localizations only will be sufficient.

\begin{proposition}\label{prop:lin_decay}
There exists $C>0$ such that 
\begin{align}
	\norm{e^{it\Omega_\kappa}P_kf}_{L^\infty}&\leq C \kappa^{-2} t^{-1}\ip{2^k\kappa}^3\norm{P_kf}_{L^1},\label{eq:Omega_decay}\\
	\norm{e^{it\Sigma_\kappa}P_kf}_{L^\infty}&\leq C  \kappa^{-\frac{5}{3}} t^{-1}\ip{2^k\kappa}^{\frac{7}{3}} 2^{\frac{k}{3}}\norm{P_kf}_{L^1}.\label{eq:Sigma_decay}
\end{align}   
\end{proposition}
\begin{proof}
The proof is direct by balancing the decay in Proposition \ref{prop:full_lin_decay} with set size estimates. For example,
\begin{align*}
	\norm{e^{it\Omega_\kappa}P_kf}_{L^\infty}
	& \lesssim \sum_{p,q \leq 0, \max\{p,q\}=0} \norm{e^{it\Omega_\kappa}P_{k,p,q}f}_{L^\infty} \\
	& \lesssim \sum_{p,q \leq 0, \max\{p,q\}=0} \norm{P_{k,p,q}f}_{L^1} \min\{ t^{-\frac{3}{2}} \k^{-3} \ip{\k 2^k}^{\frac{9}{2}} 2^{-\frac{3}{2}k - p - \frac{q}{2}}, 2^{3k+2p+q} \} \\
	& \lesssim \norm{P_k f}_{L^1} t^{-1} \k^{-2} \ip{\k 2^k}^3.
\end{align*}
The case of $\Sigma_\kappa$ follows similarly.
\end{proof}
With only $P_k$ localization, this decay at rate $t^{-1}$ is optimal; see Appendix \ref{appdx-optimal-decay-rate}.

\begin{remark}
Although \eqref{eq:Sigma_decay} is the result one can obtain directly from \eqref{eq:Sigma_decay_full}, modifying the proof to use only 2 directions instead of all 3 during the stationary phase argument, one can obtain slightly better $t^{-1}$ decay for the  $\Sigma$ evolution,
\begin{equation*}
	\Vert e^{it\Sigma_\kappa} P_{k} f \Vert_{L^\infty} \leq C \k^{-1} \ip{\k 2^k} 2^k t^{-1} \Vert P_{k} f \Vert_{L^1}.
\end{equation*}
In particular, in the proof of \eqref{eq:Sigma_decay_full}, we are sacrificing in terms of coefficients to get $t^{-3/2}$ decay by exploiting the very faint decay in the last direction. The difference will not be so relevant for the proof of Theorem \ref{main-theorem}.
\end{remark}

With the rescaling \eqref{eq:rescale_lin} and $\kappa=c\eps$, this implies the following decay rates for the dispersion relations of \eqref{eq:CER2_lin}:
\begin{corollary}\label{cor:lin_decay}
With $C>0$ as in Proposition \ref{prop:lin_decay}, there holds that
\begin{equation}
	\begin{aligned}    
		\norm{e^{itc\Omega_{c\eps}}P_kf}_{L^\infty}&\leq C c^{-3}\eps^{-2} t^{-1}\ip{2^kc\eps}^3\norm{P_kf}_{L^1},\\
		\norm{e^{itc\Sigma_{c\eps}}P_kf}_{L^\infty}&\leq C  c^{-\frac{8}{3}}\eps^{-\frac{5}{3}} t^{-1}\ip{2^kc\eps}^{\frac{7}{3}} 2^\frac{k}{3}\norm{P_kf}_{L^1}.
	\end{aligned} 
\end{equation}  
\end{corollary}

From the amplitude decay of Proposition \ref{prop:lin_decay} resp.\ Corollary \ref{cor:lin_decay}, we can also obtain the following Strichartz estimates, which will be instrumental in the proof of Theorem \ref{main-theorem}. We state them here directly for the dispersion relations of \eqref{eq:CER2_lin}.
\begin{theorem}\label{Strichartz}
Let $(a,b), (q,r)\in [2,\infty]^2$ be admissible in the sense that
\begin{equation}
	\frac{2}{a} + \frac{2}{b} = 1=\frac{2}{q} + \frac{2}{r} ,\qquad (a,b),(q,r) \neq (2,\infty).
\end{equation}
Then there holds for $f\in\mathcal{S}(\R^3)$ that
\begin{equation}\label{eq:L2Strichartz}
	\begin{aligned}    
		\Vert e^{itc\Omega_{c\eps}} P_k f \Vert_{L_t^q L_x^r} &\lesssim \blue{ (\langle 2^k c\eps\rangle^3 \eps^{-2}c^{-3} )^{\frac{1}{q}} \Vert P_k f \Vert_{L_x^2(\mathbb{R}^3)}, }\\
		\Vert e^{itc\Sigma_{c\eps}} P_k f \Vert_{L_t^q L_x^r} &\lesssim (\langle 2^k c\eps \rangle^{\frac{7}{3}} 2^{\frac{1}{3}k} \varepsilon^{-\frac{5}{3}} c^{-\frac{8}{3}} )^{\frac{1}{q}} \Vert P_k f \Vert_{L_x^2(\mathbb{R}^3)},
	\end{aligned}    
\end{equation}
and for $F\in C_t\mathcal{S}_x(([0,T],\R^3)$ that
\begin{equation}\label{eq:dualStrichartz}
	\begin{aligned}    
		\left\Vert \int_0^t e^{i(t-s)c\Omega_{c\eps}} P_k F(s,x) ds \right\Vert_{L_t^q L_x^r} &\lesssim \blue{ (\langle 2^k c\eps\rangle^3 \eps^{-2}c^{-3} )^{\frac{1}{q} + \frac{1}{a}} \Vert P_k F \Vert_{L_t^{a'} L_x^{b'}}, }\\
		\left\Vert \int_0^t e^{i(t-s)c\Sigma_{c\eps}} P_k F(s,x) ds \right\Vert_{L_t^q L_x^r} &\lesssim (\langle 2^k c\eps \rangle^{\frac{7}{3}} 2^{\frac{1}{3}k} \varepsilon^{-\frac{5}{3}} c^{-\frac{8}{3}} )^{\frac{1}{q} + \frac{1}{a}} \Vert P_k F \Vert_{L_t^{a'} L_x^{b'}},
	\end{aligned}    
\end{equation}
where $a', b'$ denote the H\"older conjugates of $a, b$, respectively.
\end{theorem}
\begin{proof}
This follows by combining Corollary \ref{cor:lin_decay} with a standard $TT'$ argument -- see e.g.\ \cite{KeelTao98}.
\end{proof}

\begin{remark}
The bounds in Theorem \ref{Strichartz} show that we can gain smallness in the linear evolution, albeit with derivative loss, provided that $\eps$ is small and that $c \gtrsim \varepsilon^{-\frac{2}{3}-}$. This, in particular, includes the case $c\varepsilon=1, \varepsilon \to 0$ investigated by Ngo and Scrobogna in \cite{NgoScr18}, and the case $c\varepsilon \gtrsim 1, \varepsilon \to 0$ studied by Mu in \cite{Mu}.
\end{remark}

The proof of Theorem \ref{main-theorem} using these estimates is reserved to the last section. We first prove Proposition \ref{prop:full_lin_decay} in Section \ref{SecDisp}, using a quantified stationary phase theorem in Section \ref{SecQuantStatPhi}.

\subsubsection*{The incompressible case} We remark that in the limit $c\to\infty$ we also recover the standard Strichartz estimates for the linearized dynamics in the incompressible case \cite[Theorem 1.1]{KLT14b}: From \eqref{eq:Omega_decay} we recall that
\begin{equation}
     \Vert e^{itc\Omega_{c\eps}} P_k f \Vert_{L_t^q L_x^r} \leq C(\langle 2^k c\eps\rangle^3 \epsilon^{-2}c^{-3} )^{\frac{1}{q}} \Vert P_k f \Vert_{L_x^2(\mathbb{R}^3)}. \label{CER-Strichartz}
\end{equation}
Since (see \eqref{ExpansionsSigmaOmega}),
\begin{equation}
 c\Omega_{c\eps} \to \omega_\eps := \eps^{-1}\frac{\xi_3}{|\xi|}, \qquad c\to\infty,
\end{equation}
where in $\omega_\eps$ we recognize the dispersion relation of the incompressible dynamics \eqref{EC}, it follows from the weak convergence $e^{itc\Omega_{c\eps}} P_k f \rightharpoonup e^{it\omega_\eps} P_k f$, $c\to\infty$, that
\begin{equation}\label{EC-Strichartz}
    \Vert e^{it\omega_{\eps}} P_k f \Vert_{L_t^q L_x^r} \lesssim (\eps 2^{3k})^{\frac{1}{q}} \norm{P_k f}_{L^2}.
\end{equation}
The inhomogeneous Strichartz estimate follows.

\section{Quantified stationary phase theorem}\label{SecQuantStatPhi}
Here we state concrete conditions for a stationary phase theorem in two dimensions, which will be used throughout the next section.\\

Recall that $\vphi$ is used to denote a bump function around $1$ (likewise for $\vphi_1, \vphi_2, ...$). $\bvphi \in C^{\infty}(-2,2)$ is the function satisfying $\bar{\vphi}(x) + \underset{n\geq0}{\sum}\vphi(2^{-n}x) \equiv 1$, and $\tvphi(x) = \vphi(x)/x$ (likewise for $\tvphi_1, \tvphi_2$). As is well known, qualitative properties of $\vphi$ and $\tvphi$ (and $\tilde{\tvphi}$) are all the same.

\begin{theorem}\label{stationary-phase-thm}
(2-dimensional Stationary phase theorem) Let $L\in C^2(\R^2)$, and assume that there exist $M>0$ and $f_j\in C^0(\R^2)$, $ L_j\in C^0(\R)$, $j\in\{1,2\}$, such that
\begin{gather*}
	|L| \leq M, \quad \iint |\partial_1 \partial_2 L | dx \leq M, \\
	|\partial_j L|, \frac{|L|}{|f_j|}(x_1,x_2) \leq L_j(x_j), \quad \int L_j(x_j) dx_j \leq M.
\end{gather*}
Given $\Phi\in C^3(\textnormal{supp}(L))$, assume moreover that there exist $a,A>0$, such that on the support of each $L(x_1,x_2) \cdot g_1(\done\Phi) \cdot g_2(\dtwo\Phi)$, where $g_1(\lambda) \in \{\bvphi(a^{-1}\lambda)\} \cup \{\vphi(a^{-1}2^{-n}\lambda) : n \geq 0\}$ and $g_2(\lambda) \in \{\bvphi(A^{-1}\lambda)\} \cup \{\vphi(A^{-1}2^{-m}\lambda) : m \geq 0\}$, one of the following two conditions holds for $\{j,k\} = \{1,2\}$:
\begin{enumerate}
\item (Diagonal-dominant case) 
\begin{equation}\label{diag-dom}
	|\det\nabla^2\Phi| \sim |\partial_1^2 \Phi| |\partial_2^2 \Phi|,\quad
	|\partial_k\partial_j^2\Phi| \lesssim \frac{1}{|f_k|}|\partial_j^2 \Phi|, \quad
	|\partial_1^2\Phi|^{\frac{1}{2}} \sim a, \quad
	|\partial_2^2\Phi|^{\frac{1}{2}} \sim A,
\end{equation}
Note that $|\done\dtwo\Phi| \lesssim aA$ follows directly.
\item (Off-diagonal-dominant case) 
\begin{equation}\label{nondiag-dom}
	\qquad |\det\nabla^2\Phi| \sim |\done \dtwo \Phi|^2, \hspace{3mm}
	|\partial_j^2\Phi| \lesssim |\done\dtwo\Phi|, \hspace{3mm} |\partial_k\partial_j^2\Phi| \lesssim \frac{1}{|f_j|}|\done\dtwo\Phi|, \hspace{3mm}
	|\done\dtwo\Phi|^{\frac{1}{2}} \sim a = A.
\end{equation}
In addition,
\begin{equation}\label{small-area-cond}
	|\partial_i^2 L| \lesssim \frac{|L|}{|f_i|^2}, \qquad \frac{1}{|f_i|} \lesssim a, \quad i \in \{1,2\}
\end{equation}
hold in the support of $L(x) \cdot \bvphi(a^{-1}\partial_{3-i}\Phi) \cdot \vphi(a^{-1}2^{-n}\partial_i\Phi)$ for each $n \geq 0$.
\end{enumerate}
Then, if $x \mapsto (\nabla\Phi)(x)$ is injective on the support of $\Phi \cdot L$,
\begin{equation}
	\left\vert \iint e^{i\Phi} L(x_1,x_2) dx_1dx_2 \right\vert \lesssim (\det\nabla^2\Phi)^{-\frac{1}{2}} M.
\end{equation}
\end{theorem}

\begin{remark}
\begin{enumerate}
	\item We stress that only fairly natural conditions are imposed on the amplitude and the phase. This is more prominent in the diagonal-dominant case, \eqref{diag-dom}, where they only need to follow type of homogeneity conditions.
	\item The latter of the conditions \eqref{small-area-cond} is based on the idea that when performing the stationary phase estimate in the region localized by $\bvphi(a^{-1}\partial_j\Phi)$, the area has to be small in some sense. In applications of the theorem in the following sections, the first condition will naturally hold since $L$ will make $|f_j| \sim |x_j|$, and the latter condition will follow after trimming out the large derivative regions and small time cases.
	\item The conditions here are a bit more general than what we need in the following sections. Namely, $f_j$'s will be constants in our applications, but are allowed to be functions as this generalization does not interrupt the proof. The injectivity requirement, however, is stricter than necessary, and is here only to prevent the possibility of infinite preimage. This can be relaxed in the usual way, for example, by requiring an upper bound on $|(\nabla\Phi)^{-1}(y)|$ for a.e. $y$.
	\item As one can check, the proof works by decomposing the plane into the central region around stationary point ($I$ in the proof), two strips from the center (one vertical and one horizontal, $I_n$ and $I^m$'s), and diagonal regions outside ($I_{n,m}$'s). This will be the recurring pattern we also use to position ourselves into a situation where we can apply the theorem in the following sections. The assumptions \eqref{diag-dom} and \eqref{nondiag-dom} claims the strips are (roughly) vertical and horizontal, which can always be enforced. The next proposition deals with $I_{n,m}$'s whose proof follows the spirit of the usual (non-)stationary phase argument.
\end{enumerate}
\end{remark}

\begin{proposition}\label{stationary-phase-lemma}
    Let $L\in C^2(\R^2)$, and assume that there exist $M>0$ and $f_j\in C^0(\R^2)$, $ L_j\in C^0(\R)$, $j\in\{1,2\}$, such that
    \begin{gather*}
        |L| \leq M, \quad \iint |\partial_1 \partial_2 L | dx \leq M, \\
        |\partial_j L|, \frac{|L|}{|f_j|}(x_1,x_2) \leq L_j(x_j), \quad \int L_j(x_j) dx_j \leq M.
    \end{gather*}
    Given $\Phi\in C^3(\textnormal{supp}(L))$ and some $n,m \geq 0$, assume moreover that there exist $a,A>0$, such that on the support of $L(x_1,x_2) \cdot \vphi_1(a^{-1}2^{-n} \done\Phi) \cdot \vphi_2(A^{-1}2^{-m} \dtwo\Phi)$ one of the following two conditions holds:
    \begin{equation}
    \begin{aligned}
        \text{(1)} \quad  |\det\nabla^2\Phi| &\sim |\partial_1^2 \Phi| |\partial_2^2 \Phi|,\quad
        |\partial_k\partial_j^2\Phi| \lesssim \frac{1}{|f_k|}|\partial_j^2 \Phi|,\quad |\partial_1^2\Phi|^{\frac{1}{2}} \sim a,\quad |\partial_2^2\Phi|^{\frac{1}{2}} \sim A,\\
        \text{(2)}\quad  |\det\nabla^2\Phi| &\sim |\done \dtwo \Phi|^2, \quad
        |\done^2\Phi|, |\dtwo^2\Phi| \lesssim |\done\dtwo\Phi|, \quad |\partial_k\partial_j^2\Phi| \lesssim \frac{1}{|f_j|}|\partial_k\partial_j\Phi|,\quad |\done\dtwo\Phi|^{\frac{1}{2}} \sim a = A,
    \end{aligned}
    \end{equation}
    where $\{j,k\}=\{1,2\}$.
    Then, if $x \mapsto (\nabla\Phi)(x)$ is injective on the support of $\Phi \cdot L$,
    \begin{equation}\label{stationary-phase-estimate}
        \left\vert \iint e^{i\Phi} \vphi_1(a^{-1}2^{-n} \partial_1\Phi) \vphi_2(A^{-1}2^{-m} \partial_2\Phi) L(x_1,x_2) dx_1dx_2 \right\vert \lesssim (aA)^{-1} M 2^{-n-m}.
    \end{equation}
\end{proposition}

\begin{proof}[Proof of Theorem \ref{stationary-phase-thm}]
Decompose the integral into $I + \underset{n \geq 0}{\sum} I_n + \underset{m \geq 0}{\sum} I^m + \underset{n,m \geq 0}{\sum} I_{n,m}$ where
\begin{gather*}
	I = \iint e^{i\Phi} \bvphi(a^{-1} \done\Phi) \bvphi(A^{-1} \dtwo\Phi) L dx, \\
	I_n = \iint e^{i\Phi} \vphi(a^{-1}2^{-n} \done\Phi) \bvphi(A^{-1} \dtwo\Phi) L dx, \\
	I^m = \iint e^{i\Phi} \bvphi(a^{-1} \done\Phi) \vphi(A^{-1}2^{-m} \dtwo\Phi) L dx, \\
	I_{n,m} = \iint e^{i\Phi} \vphi(a^{-1}2^{-n} \done\Phi) \vphi(A^{-1}2^{-m} \dtwo\Phi) L dx.
\end{gather*}
For notational convenience, from hereon we will write single integral sign instead of two. Then,
\begin{equation*}
	|I| = \left\vert \int e^{i\Phi} \bar{\vphi}(a^{-1}\done\Phi) \bar{\vphi}(A^{-1}\dtwo\Phi) L(x_1,x_2) dx_1dx_2 \right\vert \lesssim (\det\nabla^2\Phi)^{-\frac{1}{2}} M,
\end{equation*}
by changing the variables. $I_{n,m}$'s are dealt with by Proposition \ref{stationary-phase-lemma}. For $I^m$, we showcase the harder case \eqref{nondiag-dom} which involves using condition \eqref{small-area-cond} (the case \eqref{diag-dom} does not need this extra assumption). Via integration by parts,
\begin{equation*}
	|I^m| = a^{-1}2^{-m} \left\vert \int e^{i\Phi} \dtwo\left[ \bar{\vphi}(a^{-1}\done\Phi) \tvphi(a^{-1}2^{-m}\dtwo\Phi) L(x_1,x_2) \right] dx \right\vert.
\end{equation*}
Note that $a=A$ in case \eqref{nondiag-dom}. The easiest term is when the derivative falls on $\tvphi(a^{-1}2^{-m}\dtwo\Phi)$:
\begin{align*}
	& a^{-1}2^{-m} \left\vert \int e^{i\Phi} \bar{\vphi}(a^{-1}\done\Phi) \cdot a^{-1}2^{-m}\dtwo^2\Phi \cdot \tvphi'(a^{-1}2^{-m}\dtwo\Phi) L(x_1,x_2) dx \right\vert \\
	\lesssim & \hspace{1mm} a^{-1}2^{-m} \int \bvphi(a^{-1}\done\Phi) \cdot a^{-1}2^{-m} a^2 \cdot |\tvphi'|(a^{-1}2^{-m}\dtwo\Phi) \cdot dx \cdot M \lesssim a^{-2}2^{-m}M,
\end{align*}
by changing the variable from $(x_1,x_2)$ to $(\done\Phi,\dtwo\Phi)$. When it falls on $L$, we do another integration by parts and get
\begin{align*}
	& a^{-1}2^{-m} \left\vert \int e^{i\Phi} \bar{\vphi}(a^{-1}\done\Phi) \tvphi(a^{-1}2^{-m}\dtwo\Phi) \dtwo L(x_1,x_2) dx \right\vert \\
	= & \hspace{1mm} a^{-2}2^{-2m} \left\vert \int e^{i\Phi} \dtwo[ \bar{\vphi}(a^{-1}\done\Phi) \tilde{\tvphi}(a^{-1}2^{-m}\dtwo\Phi) \dtwo L(x_1,x_2) ] dx \right\vert \\
	\lesssim & \hspace{1mm} a^{-2}2^{-2m} \left[ \int a^{-1}|\done\dtwo\Phi| |\bvphi'|(a^{-1}\done\Phi) \tilde{\tvphi}(a^{-1}2^{-m}\dtwo\Phi) L_2(x_2) dx \right. \\
	& \hspace{15mm} + \int a^{-1}2^{-m}|\dtwo^2\Phi| |\tilde{\tvphi}'|(a^{-1}2^{-m}\dtwo\Phi) L_2(x_2) dx \\
	& \hspace{15mm} \left. + \int \bvphi(a^{-1}\done\Phi) \tilde{\tvphi}(a^{-1}2^{-m}\dtwo\Phi) \frac{M}{|f_2|^2} dx \right] \lesssim a^{-2}2^{-m}M.
\end{align*}
Finally, if the derivative falls on $\bvphi(a^{-1}\done\Phi)$,
\begin{align*}
	& a^{-1}2^{-m} \left\vert \int e^{i\Phi} a^{-1}\done\dtwo\Phi\bar{\vphi}'(a^{-1}\done\Phi) \tvphi(a^{-1}2^{-m}\dtwo\Phi) L(x_1,x_2) dx \right\vert \\
	= & \hspace{1mm} a^{-2}2^{-2m} \left\vert \int e^{i\Phi} \dtwo[ a^{-1}\done\dtwo\Phi\bar{\vphi}'(a^{-1}\done\Phi) \tilde{\tvphi}(a^{-1}2^{-m}\dtwo\Phi) L(x_1,x_2) ] dx \right\vert \\
	\lesssim & \hspace{1mm} a^{-2}2^{-2m} \left[ \int a^{-1} \frac{1}{|f_2|}|\done\dtwo\Phi| |\bvphi'|(a^{-1}\done\Phi) \tilde{\tvphi}(a^{-1}2^{-m}\dtwo\Phi) \cdot M dx \right. \\
	& \hspace{15mm} + \int a^{-2} |\done\dtwo\Phi|^2 |\bvphi''|(a^{-1}\done\Phi) \tilde{\tvphi}(a^{-1}2^{-m}\dtwo\Phi) \cdot M dx \\
	& \hspace{15mm} + \int a^{-2}2^{-m}|\done\dtwo\Phi| |\dtwo^2\Phi| |\bvphi'|(a^{-1}\done\Phi) |\tilde{\tvphi}'|(a^{-1}2^{-m}\dtwo\Phi) \cdot M dx \\
	& \hspace{15mm} \left. + \int a^{-1} |\done\dtwo\Phi| \tilde{\tvphi}(a^{-1}2^{-m}\dtwo\Phi) L_2(x_2) dx \right] \lesssim a^{-2}2^{-m} M,
\end{align*}
which are all summable in $m$. $I_n$'s can be treated similarly.
\end{proof}

Now we prove Proposition \ref{stationary-phase-lemma}
\begin{proof}[Proof of Proposition \ref{stationary-phase-lemma}]
Noting that the sizes of $\done\Phi, \dtwo\Phi$ are specified, we integrate by parts in both directions. Then, the integral in \eqref{stationary-phase-estimate} becomes
\begin{equation*}
	- a^{-1}2^{-n}\cdot A^{-1}2^{-m} \cdot \int e^{i\Phi} \done\dtwo \left[ \tvphi_1(a^{-1}2^{-n} \partial_1\Phi) \tvphi_2(A^{-1}2^{-m} \partial_2\Phi) L(x_1,x_2) \right] dx.
\end{equation*}
Hence, from hereon, we only need to show that the integral above is bounded by $M$. First,
\begin{gather*}
	\left\vert \int e^{i\Phi} \partial_j \left[ \tvphi_1(a^{-1}2^{-n} \done\Phi) \tvphi_2(A^{-1}2^{-m} \dtwo\Phi) \right] \partial_k[L(x_1,x_2)] dx \right\vert \lesssim M, \\
	\left\vert \int e^{i\Phi} \partial_j \left[ \tvphi_1(a^{-1}2^{-n} \done\Phi) \right] \partial_k[\tvphi_2(A^{-1}2^{-m} \dtwo\Phi)] L(x_1,x_2) dx \right\vert \lesssim M, \\
	\left\vert \int e^{i\Phi} \tvphi_1(a^{-1}2^{-n} \done\Phi) \tvphi_2(A^{-1}2^{-m} \dtwo\Phi) \hspace{1mm} \partial_j\partial_k[L(x_1,x_2)] dx \right\vert \lesssim M,
\end{gather*}
for $\{j,k\} = \{1,2\}$, by changing the variables from $(x_1,x_2)$ to $(\done\Phi,x_k), (\dtwo\Phi,x_k)$, or $(\done\Phi,\dtwo\Phi)$ appropriately, and applying the assumptions on $L$ and $\Phi$. When the derivative from second integration by parts falls on the second order derivative of the phase from the first integration by parts,
\begin{equation*}
	\left\vert \int e^{i\Phi} \dtwo(a^{-1}2^{-n} \done^2\Phi) \tvphi_1'(a^{-1}2^{-n} \done\Phi) \tvphi_2(A^{-1}2^{-m} \dtwo\Phi) L(x_1,x_2) dx \right\vert \lesssim M
\end{equation*}
follows by changing the variables from $(x_1,x_2)$ to $(\done\Phi, x_2)$ in the case \eqref{diag-dom}, or from $(x_1,x_2)$ to $(x_1,\done\Phi)$ in the case \eqref{nondiag-dom}. One can check that the analogous claim with the roles of $\done$ and $\dtwo$ swapped holds as well.

For the only remaining case when both derivatives fall on the same phase localization term, we do another integration by parts. We will show
\begin{equation*}
	\left\vert \int e^{i\Phi} a^{-1}2^{-n} \done^2\Phi \cdot a^{-1}2^{-n} \dtwo\done\Phi \cdot \tvphi_1''(a^{-1}2^{-n} \done\Phi) \tvphi_2(A^{-1}2^{-m} \dtwo\Phi) L(x_1,x_2) dx \right\vert \lesssim M(2^{-n-m} + 2^{-2n}),
\end{equation*}
and one can check that the similar result for swapping the role of $\done$ and $\dtwo$ holds as well. For notational convenience, we write $\tvphi_1$ instead of $\tvphi_1''$ as there is no qualitative difference. We integrate by parts once more in $\dtwo$ and deal with
\begin{equation*}
	A^{-1}2^{-m} \int \left\vert \dtwo \left[ a^{-2}2^{-2n} \done^2\Phi \cdot \dtwo\done\Phi \cdot \tvphi_1(a^{-1}2^{-n} \done\Phi) \tilde{\tvphi}_2(A^{-1}2^{-m}\dtwo\Phi) L(x_1,x_2) \right] \right\vert dx.
\end{equation*}
1) When $\dtwo$ falls on $L$,
\begin{align*}
	& A^{-1}2^{-m} \int \left\vert a^{-2}2^{-2n} \done^2\Phi \cdot \dtwo\done\Phi \cdot \tvphi_1(a^{-1}2^{-n} \done\Phi) \tilde{\tvphi}_2(A^{-1}2^{-m}\dtwo\Phi) \dtwo L(x_1,x_2) \right\vert dx \\
	& \lesssim A^{-1}2^{-m} \int a^{-1}A 2^{-2n} |\done^2\Phi| |\tvphi_1(a^{-1}2^{-n} \done\Phi)| L_2(x_2) dx \\
	& \lesssim A^{-1} 2^{-m} \cdot a^{-1}A2^{-2n} \cdot a2^n M = M 2^{-n-m}.
\end{align*}
2) When $\dtwo$ falls on $\tilde{\tvphi}_2$,
\begin{align*}
	& A^{-1}2^{-m} \int \left\vert a^{-2}2^{-2n} \done^2\Phi \cdot \dtwo\done\Phi \cdot \tvphi_1(a^{-1}2^{-n} \done\Phi) \dtwo[\tilde{\tvphi}_2(A^{-1}2^{-m}\dtwo\Phi)] L(x_1,x_2) \right\vert dx \\
	& \lesssim A^{-1}2^{-m} \int aA 2^{-2n} |\tvphi_1(a^{-1}2^{-n} \done\Phi)| \cdot A2^{-m} |\tilde{\tvphi}_2'(A^{-1}2^{-m}\dtwo\Phi)| dx \cdot M \\
	& \lesssim aA 2^{-2m-2n} \cdot a 2^n A 2^{m} \cdot (aA)^{-2} \cdot A^{-1} 2^{-m} \cdot A2^m \cdot M = M 2^{-n-m}.
\end{align*}
3) When $\dtwo$ falls on $\tvphi_1$,
\begin{align*}
	& A^{-1}2^{-m} \int \left\vert a^{-2}2^{-2n} \done^2\Phi \cdot \dtwo\done\Phi \cdot \dtwo[\tvphi_1(a^{-1}2^{-n} \done\Phi)] \tilde{\tvphi}_2(A^{-1}2^{-m}\dtwo\Phi) L(x_1,x_2) \right\vert dx \\
	& \lesssim A^{-1}2^{-m} \int aA 2^{-2n} \cdot A2^{-n} |\tvphi_1'(a^{-1}2^{-n} \done\Phi)| \cdot |\tilde{\tvphi}_2(A^{-1}2^{-m}\dtwo\Phi)|dx \cdot M \\
	& \lesssim aA 2^{-m-3n}\cdot a2^n \cdot A2^m \cdot (aA)^{-2} \cdot M = M 2^{-2n}.
\end{align*}
4) When $\dtwo$ falls on $\dtwo\done\Phi$, in the case \eqref{diag-dom}
\begin{align*}
	& A^{-1}2^{-m} \int \left\vert a^{-2}2^{-2n} \done^2\Phi \cdot \dtwo^2\done\Phi \cdot \tvphi_1(a^{-1}2^{-n} \done\Phi) \tilde{\tvphi}_2(A^{-1}2^{-m}\dtwo\Phi) L(x_1,x_2) \right\vert dx \\
	& \lesssim A^{-1}2^{-m} \int 2^{-2n} \cdot \frac{1}{|f_1|}|\dtwo^2\Phi| \cdot |\tilde{\tvphi}_2(A^{-1}2^{-m}\dtwo\Phi)| |L(x_1,x_2)| dx \\
	& \lesssim A^{-1}2^{-m} \cdot 2^{-2n} \cdot A2^m \cdot M = M 2^{-2n}.
\end{align*}
In the case \eqref{nondiag-dom},
\begin{align*}
	& A^{-1}2^{-m} \int \left\vert a^{-2}2^{-2n} \done^2\Phi \cdot \dtwo^2\done\Phi \cdot \tvphi_1(a^{-1}2^{-n} \done\Phi) \tilde{\tvphi}_2(A^{-1}2^{-m}\dtwo\Phi) L(x_1,x_2) \right\vert dx \\
	& \lesssim A^{-1}2^{-m} \int 2^{-2n} \cdot \frac{1}{|f_2|}|\done\dtwo\Phi| \cdot |\tilde{\tvphi}_2(A^{-1}2^{-m}\dtwo\Phi)| |L(x_1,x_2)| dx \\
	& \lesssim A^{-1}2^{-m} \cdot 2^{-2n} \cdot A2^m \cdot M = M 2^{-2n}.
\end{align*}
5) Finally, when $\dtwo$ falls on $\done^2\Phi$,
\begin{align*}
	& A^{-1}2^{-m} \int \left\vert a^{-2}2^{-2n} \dtwo\done^2\Phi \cdot \dtwo\done\Phi \cdot \tvphi_1(a^{-1}2^{-n} \done\Phi) \tilde{\tvphi}_2(A^{-1}2^{-m}\dtwo\Phi) L(x_1,x_2) \right\vert dx \\
	& \lesssim A^{-1}2^{-m} \int a^{-1}A 2^{-2n} \cdot \frac{1}{|f_2|}|\done^2\Phi| |\tvphi_1(a^{-1}2^{-n} \done\Phi)| |L(x_1,x_2)| dx \\
	& \lesssim A^{-1}2^{-m} \cdot a^{-1}A2^{-2n} \cdot a2^n \cdot M = M 2^{-n-m},
\end{align*}
in the case \eqref{diag-dom}. With \eqref{nondiag-dom},
\begin{align*}
	& A^{-1}2^{-m} \int \left\vert a^{-2}2^{-2n} \dtwo\done^2\Phi \cdot \dtwo\done\Phi \cdot \tvphi_1(a^{-1}2^{-n} \done\Phi) \tilde{\tvphi}_2(A^{-1}2^{-m}\dtwo\Phi) L(x_1,x_2) \right\vert dx \\
	& \lesssim A^{-1}2^{-m} \int 2^{-2n} \cdot \frac{1}{|f_1|}|\done\dtwo\Phi| |\tvphi_1(a^{-1}2^{-n} \done\Phi)| |L(x_1,x_2)| dx \\
	& \lesssim A^{-1}2^{-m-2n} \cdot a2^n \cdot M = M 2^{-n-m}.
\end{align*}
\end{proof}

\begin{remark}
The condition $|\done^2\Phi|, |\dtwo^2\Phi| \lesssim |\done\dtwo\Phi|$ in the case \eqref{nondiag-dom} can be removed in practice by scaling the variables as one can see in the applications below.
\end{remark}

\section{Dispersive decay estimate}\label{SecDisp}

This section is dedicated to the proof of Proposition \ref{prop:full_lin_decay}. Using Young's inequality, it is enough to prove the decay of
\begin{gather*}
	\int e^{it\Sigma + i\mathbf{x}\cdot\xi} \vphi(2^{-k}|\xi|) \vphi(2^{-k-p}|(\xi_1,\xi_2)|) d\xi, \\
	\int e^{it\Omega + i\mathbf{x}\cdot\xi} \vphi(2^{-k}|\xi|) \vphi(2^{-k-p}|(\xi_1,\xi_2)|) \vphi(2^{-k-q}\xi_3) d\xi.
\end{gather*}
Both $\Sigma$ and $\Omega$ are invariant under rotation around the $\xi_3$-axis, and so are the localization terms with $\vphi$'s. By using this axisymmetry we can assume that $\mathbf{x} = (x_1,0,-x_3)$, and it is natural to use cylindrical coordinates $\xi \to (r,\theta,z)$. Integrating in $\theta$ and using standard results on Bessel functions as in \cite[p.194]{GPW23}, the problem boils down to estimating
\begin{equation}
\begin{aligned}    
	I_{\Sigma, \pm} & :=\int_{\R} \int_0^\infty e^{it\Sigma - ix_3 z \pm ix_1 r} \vphi(2^{-k}|(r,z)|) \vphi(2^{-k-p}r) H_\pm(x_1 r) r dr dz,\\
	I_{\Omega, \pm} & :=\int_{\R} \int_0^\infty e^{it\Omega - ix_3 z \pm ix_1 r} \vphi(2^{-k}|(r,z)|) \vphi(2^{-k-p}r) \vphi(2^{-k-q}z) H_\pm(x_1 r) r dr dz,
\end{aligned}    
\end{equation}
where
\begin{equation}\label{decay-of-H}
	\left\vert \left(\frac{d}{dx}\right)^n H_\pm(x) \right\vert \lesssim \langle x \rangle^{-\frac{1}{2}-n}.
\end{equation}
We will localize further when we deal with the medium frequencies.

Throughout the section, it will be helpful to distinguish the following sizes in three main subcases:
\begin{align*}
	\textnormal{low frequencies }\k 2^k \ll 1 & \quad \Rightarrow \quad \Sigma \sim \k^{-1}, \quad \Omega \sim \xi_3, \quad \sqrt{\mathcal{D}} \sim \k^{-2}, \quad d_1,d_2 \sim \k^{-1} \\
	\textnormal{medium frequencies }\k 2^k \sim 1 & \quad \Rightarrow \quad \Sigma \sim \k^{-1}, \quad \Omega \sim \xi_3, \quad \sqrt{\mathcal{D}} \lesssim \k^{-2}, \quad \max\{d_1,d_2\} \sim \k^{-1} \\
	\textnormal{high frequencies }\k 2^k \gg 1 & \quad \Rightarrow \quad \Sigma \sim |\xi|, \quad \Omega \sim \frac{\xi_3}{\k|\xi|}, \quad \sqrt{\mathcal{D}} \sim |\xi|^2, \quad d_1,d_2 \sim 2^k.
\end{align*}

\subsection{Decay of $\Sigma$}\label{SecSigma}

We prove $t^{-\frac{3}{2}}$ decay of the $L^{\infty}$ norm of the following integral:
\begin{equation}
	I_\Sigma := \iint e^{i\Phi} \vphi(2^{-k}|r,z|) \vphi(2^{-k-p}r) H(x_1 r) r drdz.
\end{equation}
Here, $\Phi = t\Sigma - (r,z)\cdot(x_1,x_3)$. In general, we need to consider $t\Sigma -x_3 z \pm x_1 r$, but we only deal with the negative sign case as there's no qualitative difference. Also, whenever needed, we will assume $z \geq 0$ using the symmetry without loss of generality. We write here some derivatives of $\Phi$ that will appear often:
\begin{equation}
\begin{gathered}\label{size-of-Phi-Sigma}
	\dr\Phi = t\frac{r\Sigma}{d_1 d_2} - x_1, \\
	\dz\Phi = \frac{t}{2}\left( \frac{z-\k^{-1}}{d_1} + \frac{z+\k^{-1}}{d_2} \right) - x_3, \\
	\dr^2 \Phi = \frac{t}{2} \left(\frac{(z-\k^{-1})^2}{d_1^3} + \frac{(z+\k^{-1})^2}{d_2^3}\right), \\
	\dz^2 \Phi = \frac{t}{2} \left(\frac{r^2}{d_1^3} + \frac{r^2}{d_2^3}\right), \\
	\dr^3\Phi = -\frac{3tr}{2} \left( \frac{(z-\k^{-1})^2}{d_1^5} + \frac{(z+\k^{-1})^2}{d_2^5} \right), \\
	\dr\dz^2\Phi = \frac{tr}{2}\left( \frac{2(z-\k^{-1})^2 - r^2}{d_1^5} + \frac{2(z+\k^{-1})^2 - r^2}{d_2^5} \right), \\
	\dz\dr^2\Phi = t\frac{z-\k^{-1}}{2}\frac{2r^2 - (z-\k^{-1})^2}{d_1^5} + t\frac{z+\k^{-1}}{2}\frac{2r^2 - (z+\k^{-1})^2}{d_2^5},
\end{gathered}
\end{equation}
whose sizes will vary according to the subcases. We also note here
\begin{equation}\label{basic-dr-loss}
	\begin{gathered}
		|\partial_r^n (H_\pm(x_1r))| \lesssim \frac{x_1^n}{\langle x_1 r\rangle^{n+\frac{1}{2}}} \leq r^{-n}\langle x_1 r\rangle^{-\frac{1}{2}} \sim (2^{-k-p})^n \langle x_1 r\rangle^{-\frac{1}{2}}, \\
		| \partial_r^n ( \varphi(2^{-k-p}r) ) | \lesssim (2^{-k-p})^n, \\
		| \partial_r^n ( \varphi(2^{-k}|(r,z)|) ) | \lesssim (2^{-k})^n,
	\end{gathered}
\end{equation}
where the first line follows from \eqref{decay-of-H}.

\subsubsection{Low frequency, $\k 2^k \ll 1$.} If $x_1 \nsim t\k2^{k+p}$, we integrate by parts in $\dr$ first. When the derivative falls on terms not related to the phase, \eqref{basic-dr-loss} shows that there is at most $2^{-k-p}$ loss from each $\dr$ derivative. Also, using \eqref{size-of-Phi-Sigma}, we get $|\dr\Phi| \gtrsim t\k 2^{k+p}$ from $\dr\Sigma \sim 2^{k+p} \frac{\k^{-1}}{\k^{-2}} = \k 2^{k+p}$, $\dr^2 \Phi \sim t\k$, and $\dr^3 \Phi \sim t2^{k+p}\k^3 \ll t\k 2^{-k+p}$, so that, in an averaged sense, $\dr$ derivative always yields $2^{-k-p}$ loss, i.e.
\begin{gather*}
	\left\vert \frac{\partial_r^2 \Phi}{(\partial_r \Phi)^2} \right\vert \lesssim t^{-1} \k^{-1} 2^{-2k-2p}, \\
	\left\vert \frac{\partial_r^3 \Phi}{(\partial_r \Phi)^3} \right\vert \lesssim t^{-2} \k^{-2} 2^{-4k-2p} \leq (t^{-1}\k^{-1} 2^{-2k-2p})^2.
\end{gather*}
Thus,
\begin{equation*}
	|I_\Sigma| \lesssim 2^{3k+2p} \cdot \min\{ 1, (t\k2^{2k+2p})^{-2} \} \leq \k^{-\frac{3}{2}} 2^{-p} t^{-\frac{3}{2}} \ll \k^{-\frac{5}{2}} 2^{-k-p} t^{-\frac{3}{2}}.
\end{equation*}
Hence, now we assume $x_1 \sim t\k2^{k+p}$. Now, from \eqref{size-of-Phi-Sigma},
\begin{gather*}
	\dr^2 \Phi = \frac{t}{2} \left(\frac{(z-\k^{-1})^2}{d_1^3} + \frac{(z+\k^{-1})^2}{d_2^3}\right) \sim t\k, \\
	\dz^2 \Phi = \frac{t}{2} \left(\frac{r^2}{d_1^3} + \frac{r^2}{d_2^3}\right) \sim t\k^3 2^{2k+2p}, \\
	\dr\dz^2\Phi = t\frac{r}{2}\left( \frac{2(z-\k^{-1})^2 - r^2}{d_1^5} + \frac{2(z+\k^{-1})^2 - r^2}{d_2^5} \right) \sim t\k^3 2^{k+p}, \\
	\dz\dr^2\Phi = t\frac{z-\k^{-1}}{2}\frac{2r^2 - (z-\k^{-1})^2}{d_1^5} + t\frac{z+\k^{-1}}{2}\frac{2r^2 - (z+\k^{-1})^2}{d_2^5} \lesssim t\k^2, \\
	\det \nabla^2_{r,z}\Phi = t^2\k^{-2} \frac{r^2}{d_1^3 d_2^3} \sim t^2 \k^4 2^{2k+2p},
\end{gather*}
where the last line follows from \eqref{DetHess}, but by ignoring the contribution from $\partial_\theta^2\Sigma$. Together with $L(r,z) = \vphi(2^{-k}|r,z|) \vphi(2^{-k-p}r) H(x_1 r) r$, $L_r(r) = \vphi(2^{-k-p}r) H(x_1 r)$, $L_z(z) = 2^{-k}M \bvphi(2^{-k}z)$, and $M = 2^{k+p} \langle t\k 2^{2k+2p} \rangle^{-1/2}$, the conditions of Proposition \ref{stationary-phase-lemma} are all satisfied by using $f_r = 2^{k+p}$, $f_z=2^k$, $a = (t\k)^{1/2}$, and $A = (t\k^3 2^{2k+2p})^{1/2}$, and we have
\begin{equation}
	|I_\Sigma| \lesssim (t^2 \k^4 2^{2k+2p})^{-\frac{1}{2}} 2^{k+p} \langle t\k 2^{2k+2p} \rangle^{-1/2} \leq \k^{-\frac{5}{2}} 2^{-k-p} t^{-\frac{3}{2}}.
\end{equation}

\vspace{5mm}

\subsubsection{High frequency, $\k 2^k \gg 1$.}

If $x_1 \nsim t2^{p}$, since $|\dr\Sigma| \sim 2^p$, we have $|\dr\Phi| \gtrsim t 2^p$, and thus can integrate by parts along $\dr$. Along with the common losses \eqref{basic-dr-loss}, $|\dr^2 \Phi| \lesssim t\cdot 2^{-k}$, and $\dr^3 \Phi \lesssim t \cdot 2^{-2k+p}$ by \eqref{size-of-Phi-Sigma} which shows that the derivatives yield $2^{-k-p}$ loss at worst. Thus, two integration by parts give us
\begin{equation}
	|I_\Sigma| \lesssim 2^{3k+2p} \cdot \min\{ 1, (t^{-1}2^{-k-2p})^2 \} \leq 2^{\frac{3}{2}k - p} t^{-\frac{3}{2}} \ll \k 2^{\frac{5}{2}k - p} t^{-\frac{3}{2}}.
\end{equation}
Hence, from hereon, we assume $x_1 \sim t 2^{p}$. We now change the variable to the spherical variables $(\rho, \phi)$(but since we are in dimension 2 with $r$ and $z$, one can also think of this as the polar coordinates), and use the derivatives $\drho = \frac{r}{\rho}\dr + \frac{z}{\rho}\dz$, $\dphi = z\dr - r\dz$.
From integration by parts in $\rho$, using
\begin{equation}
	I_{\Sigma,n} := \int e^{i\Phi} \vphi(t^{-1}\k 2^{k-n} \drho\Phi)
	\vphi(2^{-k}\rho) \vphi(2^{-p} \sin\phi) H(x_1 \rho\sin\phi) \rho^2 \sin\phi \hspace{1mm} d\rho d\phi,
\end{equation}
we get
\begin{align*}
	|I_{\Sigma,n}| = & \hspace{1mm} t^{-1}\k 2^{k-n} \left\vert \int e^{i\Phi} \drho \left\{ \tvphi(t^{-1}\k 2^{k-n} \drho\Phi) \vphi(2^{-k}\rho) \vphi(2^{-p} \sin\phi) H(x_1 \rho\sin\phi) \rho^2 \sin\phi \right\} d\rho d\phi \right\vert \\
	\leq & \hspace{1mm} (t^{-1}\k 2^{k-n})^2 \int |\drho^2\Phi| |\tvphi'(t^{-1}\k 2^{k-n} \drho\Phi)| |\vphi(2^{-k}\rho) \vphi(2^{-p} \sin\phi) H(x_1 \rho\sin\phi) \rho^2 \sin\phi| d\rho d\phi \\
	& + t^{-1}\k 2^{k-n} \int \big\vert \drho \left\{ \vphi(2^{-k}\rho) \vphi(2^{-p} \sin\phi) H(x_1 \rho\sin\phi) \rho^2 \sin\phi \right\} \big\vert d\rho d\phi \\
	\lesssim & \hspace{1mm} t^{-1}\k 2^{k-n} 2^{2k+2p} (t 2^{k+2p})^{-\frac{1}{2}} \leq 2^{-n}\k 2^{\frac{5}{2}k - p} t^{-\frac{3}{2}}.
\end{align*}
Summing over $n$ for $2^n \gtrsim 1$, we obtain the desired result. Hence, it remains to deal with the region where $|\drho\Phi| \ll t\k^{-1}2^{-k}$, which in particular makes
\begin{equation}
	|(r,z)\cdot(x_1,x_3)| = \left| \rho\drho\Phi - t\Sigma - \frac{t}{2\k}\left( \frac{z-\k^{-1}}{d_1} - \frac{z+\k^{-1}}{d_2} \right) \right| \sim t 2^k.
\end{equation}
We additionally localize according to the size of $\dphi\Phi$ to remove the region where its size is large, using
\begin{equation*}
	I_\Sigma^m = \int e^{i\Phi} \vphi(t^{-1}\k 2^{-p-m} \dphi\Phi) \bvphi(t^{-1}\k 2^k \drho\Phi) \vphi(2^{-k}\rho) \vphi(2^{-p} \sin\phi) H(x_1 \rho\sin\phi) \rho^2 \sin\phi \hspace{1mm} d\rho d\phi.
\end{equation*}
The integration by parts give
\begin{align*}
	|I_\Sigma^m| = & \hspace{1mm} t^{-1}\k 2^{-p-m} \\
	& \quad \cdot \left\vert \int e^{i\Phi} \dphi\left\{ \tvphi(t^{-1}\k 2^{-p-m} \dphi\Phi) \bvphi(t^{-1}\k 2^k \drho\Phi) \vphi(2^{-k}\rho) \vphi(2^{-p} \sin\phi) H(x_1 \rho\sin\phi) \rho^2 \sin\phi \right\} d\rho d\phi \right\vert \\
	\leq & \hspace{1mm} (t^{-1}\k 2^{-p-m})^2 \int |\dphi^2\Phi| |\tvphi'(t^{-1}\k 2^{-p-m} \dphi\Phi)| |\vphi(2^{-k}\rho) \vphi(2^{-p} \sin\phi) H(x_1 \rho\sin\phi) \rho^2 \sin\phi| d\rho d\phi \\
	& + t^{-1}\k 2^{-p-m} \int t^{-1}\k 2^k |\dphi\drho\Phi| |\bvphi'(t^{-1}\k 2^k \drho\Phi)| |\vphi(2^{-k}\rho) \vphi(2^{-p} \sin\phi) H(x_1 \rho\sin\phi) \rho^2 \sin\phi| d\rho d\phi \\
	& + t^{-1}\k 2^{-p-m} \int |\dphi \left\{ \vphi(2^{-k}\rho) \vphi(2^{-p} \sin\phi) H(x_1 \rho\sin\phi) \rho^2 \sin\phi \right\}| d\rho d\phi \\
	\lesssim & \hspace{1mm} t^{-1}\k 2^{-p-m} 2^{3k+p} (t 2^{k+2p})^{-\frac{1}{2}} = 2^{-m}\k 2^{\frac{5}{2}k - p} t^{-\frac{3}{2}}.
\end{align*}
Summing over $m$ for $2^m \gtrsim 1$, we obtain another acceptable bound, and can then assume that $|\dphi\Phi| \ll t\k^{-1}2^p$. From
\begin{equation}\label{Sigma-high-drhodphi-dphi-relation}
	|\drho\dphi\Phi - \frac{1}{\rho}\dphi\Phi| = \left\vert t\frac{r}{2\k\rho} \left( \frac{r^2}{d_2^3} + \frac{z(z+\k^{-1})}{d_2^3} - \frac{r^2}{d_1^3} - \frac{z(z-\k^{-1})}{d_1^3} \right)\right\vert \lesssim t\k^{-2}2^{-2k+p} \ll t\k^{-1} 2^{-k+p},
\end{equation}
we also have $|\drho \dphi \Phi| \ll t\k^{-1} 2^{-k+p}$. Moreover,
\begin{gather}
	\drho^2\Phi = \frac{t}{(\k|\xi|)^2} \frac{1}{2} \left( \frac{r^2}{d_2^3} + \frac{r^2}{d_1^3} \right) \sim t\k^{-2} 2^{-3k+2p}, \\
	\dphi^2\Phi = t\k^{-1}z\frac{\Omega}{d_1d_2} - (\k^{-1}r)^2 \frac{t}{2}\left( \frac{1}{d_2^3} + \frac{1}{d_1^3} \right) + (r,z) \cdot (x_1, x_3) \sim t2^k.
\end{gather}
We are now left with redefined $I_\Sigma$,
\begin{equation}
	\tilde{I}_\Sigma = \int e^{i\Phi} \bvphi(t^{-1}\k 2^k \drho\Phi) \bvphi(t^{-1}\k 2^{-p} \dphi\Phi) \vphi(2^{-k}\rho) \vphi(2^{-p} \sin\phi) H(x_1 \rho\sin\phi) \rho^2 \sin\phi \hspace{1mm} d\rho d\phi.
\end{equation}
In practice, thanks to the precise localizations that determine the sizes of the second and third order derivatives, we can start applying Theorem \ref{stationary-phase-thm}. We elaborate why we can proceed as such in full detail for this case (high frequency of $\Sigma$), and let the reader check how they are repeated for other cases, as these are essentially repetitions of the proof of Theorem \ref{stationary-phase-thm} as detailed below. Put $a = (t\k^{-2}2^{-3k+2p})^{1/2}, A = (t2^k)^{1/2}$, $b = t\k^{-1}2^{-k}, B = t\k^{-1}2^p$ and consider four different cases. \\

\noindent
\textcircled{1} $a \geq b, A \geq B$ : This means the above steps already yield small enough region. Simply
\begin{equation}
	|\tilde{I}_\Sigma| \lesssim \int \bvphi(b^{-1}\drho\Phi) \bvphi(B^{-1}\dphi\Phi) |L(\rho,\phi)| d\rho d\phi \lesssim bB \cdot (\det\nabla^2_{\rho,\phi}\Phi)^{-1} M \leq aA \cdot (aA)^{-2} M \leq \k 2^{\frac{5}{2}k-p} t^{-\frac{3}{2}}.
\end{equation}
Here, $L(\rho,\phi) = \vphi(2^{-k}\rho) \vphi(2^{-p} \sin\phi) H(x_1 \rho\sin\phi) \rho^2 \sin\phi$ and $M = 2^{2k+p} \langle t 2^{k+2p} \rangle^{-1/2}$. $L, M, a$, and $A$ are intended to resemble the notations in Proposition \ref{stationary-phase-lemma} and Theorem \ref{stationary-phase-thm} as these will be used as such in step \textcircled{4}.\\

\noindent
\textcircled{2} $a \geq b, A \leq B$ : Solving this in terms of $t$ gives
\begin{equation}
	\k 2^{\frac{1}{2}k - p} \leq t^{\frac{1}{2}} \leq 2^{-\frac{1}{2}k + p} \quad \Rightarrow \quad \k 2^{k-2p} \leq 1,
\end{equation}
which is against the assumption.\\

\noindent
\textcircled{3} $a \leq b, A \geq B$ : Since $a$ is smaller, we decompose $\tilde{I}_\Sigma = \tilde{I}_{\Sigma,-1} + \underset{n \geq 0}{\sum}\tilde{I}_{\Sigma,n}$ where
\begin{equation}
	\tilde{I}_{\Sigma,n} = \int e^{i\Phi} \vphi(a^{-1}2^{-n}\drho\Phi) \bvphi(B^{-1}\dphi\Phi) L(\rho,\phi) d\rho d\phi,
\end{equation}
and integrate by parts in $\drho$. If $\drho$ falls on $L$,
\begin{equation}
	a^{-1}2^{-n} \int e^{i\Phi} \tvphi(a^{-1}2^{-n}\drho\Phi) \bvphi(B^{-1}\dphi\Phi) \drho L \hspace{1mm} d\rho d\phi \lesssim a^{-1}2^{-n} \cdot B \cdot A^{-2} \cdot M \leq (aA)^{-1}M 2^{-n}.
\end{equation}
If $\drho$ falls on $\bvphi(B^{-1}\dphi\Phi)$,
\begin{equation}
\begin{split}
	& \left\vert a^{-1} 2^{-n} \int e^{i\Phi} \tvphi(a^{-1}2^{-n}\drho\Phi) B^{-1} \drho\dphi\Phi \cdot \bvphi'(B^{-1}\dphi\Phi) L(\rho,\phi) d\rho d\phi \right \vert \\
	&\quad \lesssim a^{-1} 2^{-n} \int \tvphi(a^{-1}2^{-n}\drho\Phi) 2^{-k} |\bvphi'|(B^{-1}\dphi\Phi) |L(\rho,\phi)| d\rho d\phi \\
	&\quad \lesssim a^{-1}2^{-n} \int |\bvphi'|(B^{-1}\dphi\Phi) L_\rho (\rho) d\rho d\phi \lesssim a^{-1} 2^{-n} \cdot B \cdot A^{-2} \cdot M \leq (aA)^{-1} M 2^{-n}.
\end{split}
\end{equation}
Here, \eqref{Sigma-high-drhodphi-dphi-relation} was used in the first inequality and $L_\rho = 2^{-k} \cdot M \vphi(2^{-k}\rho)$ is a corresponding choice of $L_1$ in Proposition \ref{stationary-phase-lemma} where one can check that it indeed satisfies the required condition. Lastly, if $\drho$ falls on $\vphi(a^{-1}2^{-n}\drho\Phi)$,
\begin{equation}
\begin{split}
	& \left\vert a^{-1}2^{-n} \int e^{i\Phi} a^{-1}2^{-n}\drho^2\Phi \tvphi'(a^{-1}2^{-n}\drho\Phi) \bvphi(B^{-1}\dphi\Phi) L d\rho d\phi \right\vert \\
	&\quad  \lesssim a^{-1}2^{-n} \int a^{-1}2^{-n} |\tvphi'|(a^{-1}2^{-n}\drho\Phi) \bvphi(B^{-1}\dphi\Phi) \hspace{1mm} a^2 d\rho d\phi \cdot M \\
	&\quad \lesssim a^{-1}2^{-n} \cdot B \cdot A^{-2} \cdot M \leq (aA)^{-1}M \cdot 2^{-n}.
\end{split}
\end{equation}
All of these are summable in $n$, resulting in
\begin{equation}
	\sum_{n \geq 0} |\tilde{I}_{\Sigma,n}| \lesssim (t^2 \k^{-2} 2^{-2k+2p})^{-1/2} 2^{2k+p} \langle t 2^{k+2p} \rangle^{-1/2} \leq \k 2^{\frac{5}{2}k-p} t^{-\frac{3}{2}}.
\end{equation}
$\tilde{I}_{\Sigma,-1}$ can be treated the same way as in \textcircled{1}.

\noindent
\textcircled{4} $a \leq b, A \leq B$ : This is the case where we need stationary phase estimates in both directions, which will be handled by Theorem \ref{stationary-phase-thm}. Direct calculations show that
\begin{equation}
\begin{split}
	& \dphi\drho^2\Phi = \frac{tr}{2(\k|\xi|)^2} \left( \frac{2z}{d_2^3} + \frac{2z}{d_1^3} + \frac{3\k^{-1}r^2}{d_2^5} - \frac{3\k^{-1}r^2}{d_1^5} \right) \lesssim t\k^{-2}2^{-3k+p}, \\
	& \drho\dphi^2\Phi = -\frac{tz}{2\k^2|\xi|}\left[\frac{z+\k^{-1}}{d_2^3} + \frac{z-\k^{-1}}{d_1^3}\right] + \frac{r^2}{2\k^2|\xi|} \left( \frac{1}{d_2^3} + \frac{1}{d_1^3} \right) \\ & \hspace{15mm} - \frac{3tr^2}{2\k^2|\xi|} \left( \frac{\k^{-1}(z+\k^{-1})}{d_2^5} - \frac{\k^{-1}(z-\k^{-1})}{d_1^5} \right) + \frac{1}{|\xi|}(r,z)\cdot(x_1,x_3) \sim t.
\end{split}
\end{equation}
In particular, $|\drho\dphi\Phi|^2 \ll |\drho^2\Phi||\dphi^2\Phi|$, $|\dphi\drho^2\Phi| \lesssim 2^{-p} |\drho^2\Phi|$, and $|\drho\dphi^2\Phi| \lesssim 2^{-k} |\dphi^2\Phi|$. Applying Theorem \ref{stationary-phase-thm} with \eqref{diag-dom} for $L(\rho,\phi) = \vphi(2^{-k}\rho) \vphi(2^{-p}\sin\phi) H(x_1 \rho\sin\phi) \rho^2\sin\phi$ by observing that it satisfies all conditions with $M = 2^{2k+p} \langle t 2^{k+2p} \rangle^{-1/2}$, we have
\begin{equation}
	|\tilde{I}_\Sigma| \lesssim (t^2 \k^{-2} 2^{-2k+2p})^{-1/2} 2^{2k+p} \langle t 2^{k+2p} \rangle^{-1/2} \leq \k 2^{\frac{5}{2}k-p} t^{-\frac{3}{2}},
\end{equation}
which is consistent with other estimates, and thus the final bound. \vspace{5mm}

\subsubsection{Medium frequency, $\k 2^k \sim 1$.} In this regime we are faced with the lack of smoothness of the dispersion relations at $r=0$, $z=\pm\kappa^{-1}$, see also \eqref{DetHess}. In order to treat these singularities one at a time, without loss of generality we consider the stationary phase integral only on the upper half-plane and localize in the size of $z-\k^{-1}$ by $2^{k+l}$, i.e.\ we consider
\begin{equation}
\begin{aligned}
    I_{\Sigma,l}&:= \int_{z=0}^\infty\int_{r=0}^\infty e^{i\Phi} \varphi(2^{-k}|r,z|) \varphi(2^{-k-p}r) \varphi(2^{-k-l}(z-\k^{-1})) H(x_1 r) r drdz,\quad p\leq l\leq 0,\\
     I_{\Sigma,\leq p} & := \int_{z=0}^\infty\int_{r=0}^\infty e^{i\Phi} \varphi(2^{-k}|r,z|) \varphi(2^{-k-p}r) \bvphi(2^{-k-p}(z-\k^{-1})) H(x_1 r) r drdz, \quad l \leq p.
\end{aligned}    
\end{equation}
Here, we will consider the cases $p \leq l \leq 0$ and $l \leq p$ separately, as that is the turning point of relative size between $r$ and $z-\k^{-1}$. In this section, we will alternatively use $(\rho,\phi)$ as the spherical coordinates centered at $(0,\k^{-1})$. Note that $d_1$ and $\rho$ indicate exactly same quantities and will be used interchangeably.  In terms of the above cylindrical coordinates, the spherical derivatives take the form
\begin{equation*}
	\drho = \frac{r}{d_1}\dr + \frac{z-\k^{-1}}{d_1}\dz, \quad \dphi = (z-\k^{-1})\dr - r\dz.
\end{equation*}
It is convenient to record the following derivatives of the phase before separating the cases:
\begin{equation}\label{derivatives-La2-medfreq}
\begin{split}
	\drho\Phi & = \frac{t}{2}\left( 1 + \frac{r^2+z^2-\k^{-2}}{d_1 d_2} \right) - \frac{1}{d_1} (r,z-\k^{-1}) \cdot (x_1,x_3), \quad \drho^2\Phi  = 2t \frac{\k^{-2}r^2}{d_1^2 d_2^3}, \\
	\dphi\Phi & = -t \frac{\k^{-1}r}{d_2} - (z-\k^{-1}, -r) \cdot (x_1,x_3), \\
	\dphi^2\Phi + d_1 \drho\Phi & = t \frac{(r^2+z^2-\k^{-2})^2 + d_1 d_2^3}{d_2^3}, \quad
	d_1\drho\dphi\Phi - \dphi\Phi  = t \frac{\k^{-1}r}{d_2^2} \frac{r^2+z^2-\k^{-2}}{d_2}.
\end{split}
\end{equation}

\noindent\textbf{3-1.} $\boldsymbol{p \leq l \leq 0.}$ This is when $d_1 \sim z-\k^{-1} \sim 2^{k+l}$. When $x_1 \nsim t 2^{p-l}$, from $\dr\Sigma \sim 2^{p-l}$ we obtain that $|\dr\Phi| \gtrsim t 2^{p-l}$. Since (see \eqref{size-of-Phi-Sigma})
\begin{gather*}
	\frac{\dr^2\Phi}{|\dr\Phi|} \lesssim \frac{t 2^{-k-l}}{t 2^{p-l}} = 2^{-k-p}, \\
	\left| \frac{\dr^3\Phi}{\dr\Phi} \right| \lesssim t^{-1} 2^{l-p} \cdot \frac{t}{2} \left(\frac{r(z-\k^{-1})^2}{d_1^5} + \frac{r(z+\k^{-1})^2}{d_2^5}\right) \lesssim 2^{-2k-2l} \leq (2^{-k-p})^2,
\end{gather*}
the gain from one integration by parts in $\dr$ is $(t2^{p-l})^{-1} \cdot 2^{-k-p} = t^{-1} 2^{-k-2p+l}$, up to twice. Hence, integration by parts in $\dr$ gives
\begin{equation*}
	|I_{\Sigma,l}| \lesssim 2^{3k+2p+l} \cdot \min\{1, (t2^{k+2p-l})^{-2}\} \leq 2^{\frac{3}{2}k - p + \frac{5}{2}l} t^{-\frac{3}{2}},
\end{equation*}
so we now assume $x_1 \sim t 2^{p-l}$ and restrict the size of first order derivatives. Let
\begin{align*}
	I_{\Sigma,l,n} :=& \int e^{i\Phi} \varphi((t 2^n)^{-1} \drho\Phi) \varphi(2^{-k}|r,z|) \varphi(2^{-k-p}r) \varphi(2^{-k-l}(z-\k^{-1})) H(x_1 r) r drdz \\
	=& i(t 2^n)^{-1} \int e^{i\Phi} \drho[ \tvphi((t 2^n)^{-1} \drho\Phi) \varphi(2^{-k}|r,z|) \varphi(2^{-k-p}r) \varphi(2^{-k-l}(z-\k^{-1})) H(x_1 r) r ] drdz.
\end{align*}
If the derivative falls on anything else than $\tvphi((t 2^n)^{-1} \drho\Phi)$, we can simply integrate to gain the set size, and this gives $(t 2^n)^{-1} 2^{-k-l} 2^{3k+2p+l} \langle t 2^{k+2p-l} \rangle^{-\frac{1}{2}} \leq 2^{\frac{3}{2}k + p + \frac{1}{2}l} 2^{-n} t^{-\frac{3}{2}}$. If the derivative falls on $\tvphi$, we change variables to $(\rho,\phi)$ first and then to $(\drho\Phi, \phi)$ to get
\begin{align*}
	(t 2^n)^{-1} & \int (t 2^n)^{-1}\drho^2\Phi \tvphi'((t 2^n)^{-1} \drho\Phi) \varphi(2^{-k}|r,z|) \varphi(2^{-k-p}r) \varphi(2^{-k-l}(z-\k^{-1})) H(x_1 r) r drdz \\
	&\lesssim (t 2^n)^{-1}  \int (t 2^n)^{-1}\drho^2\Phi \tvphi'((t 2^n)^{-1} \drho\Phi) \varphi(2^{-k-p}\rho\sin\phi) \varphi(2^{-k-l}\rho\cos\phi) H(x_1 \rho\sin\phi) \rho^2\sin\phi d\rho d\phi \\
	&\lesssim  (t 2^n)^{-1} \int (t 2^n)^{-1} \tvphi'((t 2^n)^{-1} y) \varphi(2^{l-p}\sin\phi) d\rho d\phi \cdot 2^{2k+p+l} \langle t 2^{k+2p-l} \rangle^{-\frac{1}{2}} dyd\phi \\
	&\lesssim (t 2^n)^{-1}  2^{p-l} 2^{2k+p+l} (t 2^{k+2p-l})^{-\frac{1}{2}} = 2^{\frac{3}{2}k + p + \frac{1}{2}l} 2^{-n} t^{-\frac{3}{2}}.
\end{align*}
Hence, summing up in $n$ for $2^n \gtrsim 1$, we get acceptable bound of $\lesssim 2^{\frac{3}{2}k + p + \frac{1}{2}l} t^{-\frac{3}{2}}$, and can assume $\drho\Phi \ll t$ from hereon. We repeat a similar process with
\begin{align*}
	I_{\Sigma,l}^m := & \int e^{i\Phi} \varphi((t 2^{k+p+\frac{1}{2}l +m})^{-1} \dphi\Phi) \varphi(2^{-k}|r,z|) \varphi(2^{-k-p}r) \varphi(2^{-k-l}(z-\k^{-1})) H(x_1 r) r drdz \\
	= & \hspace{1mm} i(t 2^{k+p+\frac{1}{2}l +m})^{-1} \\
	& \quad \cdot \int e^{i\Phi} \dphi[ \tvphi((t 2^{k+p+\frac{1}{2}l +m})^{-1} \dphi\Phi) \varphi(2^{-k}|r,z|) \varphi(2^{-k-p}r) \varphi(2^{-k-l}(z-\k^{-1})) H(x_1 r) r ] drdz,
\end{align*}
to obtain $|I_{\Sigma,l}^m| \lesssim 2^{\frac{3}{2}k - p + 2l - m} t^{-\frac{3}{2}}$. Thus we are left with the case where $|\drho\Phi| \ll t$ and $|\dphi\Phi| \ll t 2^{k+p+\frac{1}{2}l}$, and from \eqref{derivatives-La2-medfreq}
\begin{equation}
	\drho^2\Phi \sim t 2^{-k+2p-2l}, \quad \dphi^2\Phi \sim t 2^{k+l}, \quad \rho\drho\dphi\Phi - \dphi\Phi \lesssim t 2^{k+p+l}.
\end{equation}
The formula
\begin{equation}\label{Hessian-CoV}
	\rho^2|\nabla_{\rho,\phi}^2 \Phi| = \rho^4 |\nabla_{r,z}^2 \Phi|^2 - \rho^3\drho^2\Phi \cdot \drho\Phi - 2(\rho\drho\dphi\Phi - \dphi\Phi) \cdot \dphi\Phi - (\dphi\Phi)^2
\end{equation}
ensures us that we are now in the case \eqref{diag-dom} since the third derivative conditions are met by
\begin{align*}
	\dphi\drho^2\Phi & = t \frac{4\k^{-2}r}{d_1^2 d_2^5}[d_2^2(z-\k^{-1}) + 3\k^{-1}r^2] \lesssim |\drho^2\Phi| \cdot \frac{2^{3k+l}}{2^{3k+p}} \sim \frac{1}{\phi} |\drho^2\Phi|, \\
	d_1\drho\dphi^2\Phi - \dphi^2\Phi & = t \left( \frac{\k^{-1}(z-\k^{-1})}{d_2^2} + \frac{6\k^{-2}r^2}{d_2^4} \right) \frac{r^2+z^2-\k^{-2}}{d_2} - 2t\frac{\k^{-2}r^2}{d_2^3} \lesssim t 2^{k+2l} \lesssim \dphi^2\Phi,
\end{align*}
so that $|\drho\dphi^2\Phi| \lesssim \frac{1}{d_1} |\dphi^2\Phi|$. One can also check that
\begin{equation*}
L(\rho,\phi) = \varphi(2^{-k-p}\rho\sin\phi) \varphi(2^{-k-l}\rho\cos\phi) \rho^2\sin\phi H(x_1 \rho\sin\phi)
\end{equation*}
 satisfies the needed property with $M = 2^{2k+p+l} (t 2^{k+2p-l})^{-\frac{1}{2}}$. Therefore, using Theorem \ref{stationary-phase-thm}, we obtain
\begin{equation}
	|I_{\Sigma,l}| \lesssim 2^{\frac{3}{2}k - p + 2l} t^{-\frac{3}{2}}.
\end{equation}
Using, $p \leq l \leq 0$, one can see that the previous bounds are at least as strong as this, hence this is the single estimate that works for all $I_{\Sigma,l}$.\\

\noindent
\textbf{3-2.} $\boldsymbol{l \leq p.}$  We treat all $l$'s with a single localization $\bar{\varphi}(2^{-k-p}(z-\k^{-1})) := \underset{l \leq p}{\sum} \varphi(2^{-k-l}(z-\k^{-1}))$.
When $x_1 \nsim t$, from $\frac{|\dr^2\Phi|}{|\dr\Phi|} \lesssim 2^{-k-p}, \frac{|\dr^3\Phi|}{|\dr\Phi|} \lesssim 2^{-2k-2p}$, we have $|I_{\Sigma, \leq p}| \lesssim 2^{3k+3p} \cdot \min\{ 1, (t 2^{k+p})^{-2} \} \leq 2^{\frac{3}{2}k + \frac{3}{2}p} t^{-\frac{3}{2}}$. Hence we now assume $x_1 \sim t$ and start using the first derivatives first. Repeating what we did in the previous case with
\begin{align*}
	I_{\Sigma, \leq p,n} :=& \int e^{i\Phi} \varphi((t 2^n)^{-1} \drho\Phi) \varphi(2^{-k-p}r) \bvphi(2^{-k-p}(z-\k^{-1})) H(x_1 r) r drdz \\
	=& i(t 2^n)^{-1} \int e^{i\Phi} \drho[ \tvphi((t 2^n)^{-1} \drho\Phi) \varphi(2^{-k-p}r) \bvphi(2^{-k-p}(z-\k^{-1})) H(x_1 r) r ] drdz,
\end{align*}
gives $|I_{\Sigma, \leq p,n}| \lesssim 2^{\frac{3}{2}k + \frac{3}{2}p} 2^{-n} t^{-\frac{3}{2}}$. Also,
\begin{align*}
	I_{\Sigma, \leq p}^m :=& \int e^{i\Phi} \varphi((t 2^{k+\frac{3}{2}p +m})^{-1} \dphi\Phi) \varphi(2^{-k-p}r) \bvphi(2^{-k-p}(z-\k^{-1})) H(x_1 r) r drdz \\
	= & \hspace{1mm} i(t 2^{k+\frac{3}{2}p +m})^{-1} \int e^{i\Phi} \dphi[ \tvphi((t 2^{k+\frac{3}{2}p +m})^{-1} \dphi\Phi) \varphi(2^{-k-p}r) \bvphi(2^{-k-p}(z-\k^{-1})) H(x_1 r) r ] drdz \\
	\lesssim & \hspace{1mm} (t 2^{k+\frac{3}{2}p + m})^{-1} \cdot 2^{3k+3p} (t 2^{k+p})^{-\frac{1}{2}} = 2^{\frac{3}{2}k + p} 2^{-m} t^{-\frac{3}{2}}.
\end{align*}
Summing up in $n, m \geq 0$ gives acceptable results, and hence, now we are in a region where $\drho\Phi \ll t$ and $\dphi\Phi \ll t 2^{k+\frac{3}{2}p}$. Thus, from \eqref{Hessian-CoV} and
\begin{equation}
	\drho^2\Phi \sim t 2^{-k}, \quad \dphi^2\Phi \sim t 2^{k+p}, \quad \rho\drho\dphi\Phi - \dphi\Phi \lesssim t 2^{k+2p},
\end{equation}
which follows from \eqref{derivatives-La2-medfreq}, one can check that the first condition of \eqref{diag-dom} holds. The function
\begin{equation*}
L(\rho,\phi) = \varphi(2^{-k-p}r) \bvphi(2^{-k-p}(z-\k^{-1})) H(x_1\rho\sin\phi) \rho^2\sin\phi
\end{equation*}
 satisfies the conditions with size $2^{\frac{3}{2}k + \frac{3}{2}p} t^{-\frac{1}{2}}$ this time, and
\begin{align}
	\dphi\drho^2\Phi & = \drho^2\Phi \cdot \frac{2}{r d_2^2}[d_2^2 (z-\k^{-1}) + 3\k^{-1}r^2] \lesssim |\drho^2\Phi|, \\
	d_1\drho\dphi^2\Phi - \dphi^2\Phi & = t \left( \frac{\k^{-1}(z-\k^{-1})}{d_2^2} + \frac{6\k^{-2}r^2}{d_2^4} \right) \frac{r^2+z^2-\k^{-2}}{d_2} - 2t\frac{\k^{-2}r^2}{d_2^3} \lesssim t 2^{k+2p} \lesssim |\dphi^2\Phi|,
\end{align}
hence $|\drho\dphi^2\Phi| \lesssim \frac{1}{\rho} |\dphi^2\Phi|$. Therefore, all the conditions of Theorem \ref{stationary-phase-thm} are satisfied and this gives
\begin{equation}
	|I_{\Sigma, \leq p}| \lesssim (t 2^{\frac{1}{2}p})^{-1} 2^{\frac{3}{2}k + \frac{3}{2}p} t^{-\frac{1}{2}} = 2^{\frac{3}{2}k + p} t^{-\frac{3}{2}}.
\end{equation}
Again, this is the weakest among the estimates, and hence the final bound for $I_{\Sigma, \leq p}$.

\vspace{5mm}

\subsection{Decay of $\Omega$}\label{SecOmega}

We will now localize in $z$ in addition by default, and investigate the decay of
\begin{equation}
	I_\Omega := \int e^{i\Phi} \vphi(2^{-k}|r,z|) \vphi(2^{-k-p}r) \vphi(2^{-k-q}z) H(x_1 r) r \hspace{1mm} drdz.
\end{equation}
$\Phi = t\Omega - (r,z)\cdot(x_1,x_3)$, and when needed, we will assume $z \geq 0$ using the symmetry without loss of generality. We again give some derivatives of $\Phi$ that will be used several times:
\begin{equation}\label{size-of-Phi-Omega}
\begin{gathered}
	\dr\Phi = -t \frac{r\Omega}{d_1 d_2} - x_1, \\
	\dz\Phi = \frac{t}{2} \left( \frac{z+\k^{-1}}{d_2} - \frac{z-\k^{-1}}{d_1} \right) - x_3, \\
	\dr^2\Phi = \frac{t}{2} \left( \frac{(z+\k^{-1})^2}{d_2^3} - \frac{(z-\k^{-1})^2}{d_1^3} \right), \\
	\dz^2\Phi = \frac{t}{2} \left( \frac{r^2}{d_2^3} - \frac{r^2}{d_1^3} \right), \\
	\dr^3\Phi = -\frac{3tr}{2} \left( \frac{(z+\k^{-1})^2}{d_2^5} - \frac{(z-\k^{-1})^2}{d_1^5} \right), \\
	\dr\dz^2\Phi = \frac{tr}{2}\left( \frac{2(z+\k^{-1})^2 - r^2}{d_2^5} - \frac{2(z-\k^{-1})^2 - r^2}{d_1^5} \right), \\
	\dz\dr^2\Phi = t\frac{z+\k^{-1}}{2}\frac{2r^2 - (z+\k^{-1})^2}{d_2^5} - t\frac{z-\k^{-1}}{2}\frac{2r^2 - (z-\k^{-1})^2}{d_1^5}.
\end{gathered}
\end{equation}

\subsubsection{Low frequency, $\k 2^k \ll 1$.} When $x_1 \nsim t\k^2 2^{2k+p+q}$, we have $|\dr\Phi| \gtrsim t\k^2 2^{2k+p+q}$ from $\dr\Omega = -\frac{r\Omega}{d_1 d_2} \sim \k^2 2^{2k+p+q}$. Integration by parts along $\dr$ gains $(t\k^2 2^{2k+p+q})^{-1} 2^{-k-p}$ from \eqref{basic-dr-loss} and
\begin{gather*}
	\left| \frac{\dr^2\Phi}{\dr\Phi} \right| \lesssim \frac{t\k^2 2^{k+q}}{t\k^2 2^{2k+p+q}} = 2^{-k-p}, \\
	\left| \frac{\dr^3\Phi}{\dr\Phi} \right| \lesssim \frac{t\k^4 2^{2k+p+q}}{t\k^2 2^{2k+p+q}} = \k^2 \ll 2^{-2k} \leq 2^{-2k-2p},
\end{gather*}
since by \eqref{size-of-Phi-Omega}
\begin{equation*}
\begin{split}
	\dr^2\Phi &= t\frac{z^2 + \k^{-2}}{2} \left(\frac{1}{d_2^3} - \frac{1}{d_1^3}\right) + 2t\k^{-1}z \left(\frac{1}{d_2^3} + \frac{1}{d_1^3}\right) \lesssim t(\k^{-2} \cdot \k^4 2^{k+q} + \k^{-1} 2^{k+q} \k^3) \sim t\k^2 2^{k+q}, \\
	|\dr^3\Phi| &= \frac{3tr}{2}\left| \frac{(z-\k^{-1})^2}{d_1^5} - \frac{(z+\k^{-1})^2}{d_2^5} \right| \leq \frac{3tr}{2}(z^2+\k^{-2})\left|\frac{1}{d_1^5} - \frac{1}{d_2^5}\right| + tr\k^{-1}|z|\left(\frac{1}{d_1^5} + \frac{1}{d_2^5}\right) \\&\lesssim t\k^4 2^{2k+p+q}.
\end{split}
\end{equation*}
Thus, integration by parts along $\dr$ gives
\begin{equation}
	|I_\Omega| \lesssim 2^{3k+2p+q} \cdot \min\{ 1, (t\k^2 2^{3k+2p+q})^{-2}\} \leq \k^{-3} 2^{-\frac{3}{2}k - p -\frac{q}{2}} t^{-\frac{3}{2}}.
\end{equation}
Now we can assume $x_1 \sim t\k^2 2^{2k+p+q}$. From the size of $\dz^2\Phi$
\begin{equation*}
	\dz^2\Phi = t\frac{r^2}{2} \frac{(d_1 - d_2)(d_1^2 + d_1 d_2 + d_2^2)}{d_1^3 d_2^3} = -t r^2 \frac{\Omega(d_1^2 + d_1 d_2 + d_2^2)}{d_1^3 d_2^3} \sim t\k^4 2^{3k+2p+q},
\end{equation*}
we can also remove the regions with large $\dz\Phi$ by setting
\begin{gather*}
	L(r,z) \coloneqq \vphi(2^{-k}|r,z|) \vphi(2^{-k-p}r) \vphi(2^{-k-q}z) H(x_1 r) r, \\
	I_{\Omega,n} \coloneqq \int e^{i\Phi} \varphi\left( \frac{\partial_z \Phi}{t \kappa^2 2^{2k+2p+n}} \right) 
	L(r,z) \, drdz,
\end{gather*}
since $I_{\Omega,n}$ give the desired estimate after summation in $n$ from
\begin{align*}
	|I_{\Omega,n}| & = \left| \frac{i}{t \kappa^2 2^{2k+2p+n}} \int e^{i\Phi} \partial_z \left[ \tilde{\varphi}\left( \frac{\partial_z \Phi}{t \kappa^2 2^{2k+2p+n}} \right) L(r,z) \right] \, drdz \right| \\
	& \leq \left|
	\frac{i}{t \kappa^2 2^{2k+2p+n}} \int e^{i\Phi} \tilde{\varphi}^\prime\left( \frac{\partial_z \Phi}{t \kappa^2 2^{2k+2p+n}} \right) \cdot \frac{\partial_z^2 \Phi}{t \kappa^2 2^{2k+2p+n}} \cdot L(r,z) \, drdz
	\right| \\
	& \quad + \left|
	\frac{i}{t \kappa^2 2^{2k+2p+n}} \int e^{i\Phi} \tilde{\varphi}\left( \frac{\partial_z \Phi}{t \kappa^2 2^{2k+2p+n}} \right) \partial_z L(r,z) \, drdz
	\right| \\
	& \lesssim t^{-\frac{3}{2}} \kappa^{-3} 2^{-\frac{3}{2}k-p-\frac{1}{2}q} 2^{-n}.
\end{align*}
Now we only need to estimate $\tilde{I}_\Omega = \int e^{i\Phi} 
\bar{\varphi}\left( (t \kappa^2 2^{2k+2p})^{-1}\partial_z \Phi \right) 
L(r,z) \, drdz$.
From the fact that $r(\dr^2\Omega + \dz^2\Omega) = \dr\Omega \sim \k^2 2^{2k+p+q}$, one can conclude
\begin{align*}
	\dr^2 \Omega & \sim \k^2 2^{k+q},\\
	\dz^2 \Omega & \sim \k^4 2^{3k+2p+q}, \\
	\dr\dz\Omega & \sim \k^2 2^{k+p},
\end{align*}
which means now we are in a position to use Theorem \ref{stationary-phase-thm} in the case with \eqref{nondiag-dom}. Since this is the first usage of off-diagonal-dominant case, we again explain in full detail here.

The 3rd order derivative conditions are satisfied since
\begin{gather*}
	2\dz\dr^2\Omega = \frac{2(z+\k^{-1})}{d_2^3} - \frac{2(z-\k^{-1})}{d_1^3} - 3\frac{(z+\k^{-1})^3}{d_2^5} + \frac{3(z-\k^{-1})^3}{d_1^5} \lesssim \k^2 \sim 2^{-k-p} |\dr\dz\Omega|, \\
	2\dr\dz^2\Omega = \frac{4}{r}\dz^2\Omega - \frac{3r^3}{d_2^5} + \frac{3r^3}{d_1^5} \lesssim \k^4 2^{2k+p} \ll 2^{-k} |\dr\dz\Omega|.
\end{gather*}
Having said that, we will use $(2^p\partial_r, 2^q \partial_z)$ instead of $(\partial_r, \partial_z)$. The sizes change into
\begin{gather*}
	|2^p \partial_r \Phi| \lesssim  t \kappa^2 2^{2k+2p+q},
	\quad
	|2^q \partial_z \Phi| \lesssim  t \kappa^2 2^{2k+2p+q},
	\\
	|2^{2p}\partial_r^2 \Phi| \sim t \kappa^2 2^{k+2p+q},
	\quad
	|2^{2q}\partial_z^2\Phi| \sim t \kappa^4 2^{3k+2p+3q},
	\quad
	|2^p2^q\partial_r \partial_z \Phi|
	\sim t \kappa^2 2^{k+2p+q},
	\\
	|\det\nabla_{2^p\dr, 2^q\dz}^2\Phi|
	\sim t^2 \kappa^4 2^{2k+4p+2q},
	\\
	|2^{2p}2^q\partial_z \partial_r^2\Phi| \lesssim 2^{-k} |2^p2^q\partial_r \partial_z \Phi|,
	\quad
	|2^p2^{2q}\partial_r \partial_z^2 \Phi|
	\lesssim 
	2^{-k+q} |2^p2^q\partial_r \partial_z \Phi|
	\leqslant
	2^{-k}|2^p2^q\partial_r \partial_z \Phi|,
\end{gather*}
which still satisfies \eqref{nondiag-dom}, but now with $f_r = f_z = 2^k$. We set
\begin{equation*}
	a \coloneqq t^{\frac{1}{2}} \kappa 2^{\frac{1}{2}k+p+ \frac{q}{2}} \sim 
	|2^p 2^q\partial_r \partial_z \Phi|^{\frac{1}{2}},
	\quad
	b \coloneqq t \kappa^2 2^{2k+2p+q}.
\end{equation*}
\textcircled{1} $a \geq b$ : This means we are already localized enough. Change of variable brings
\begin{align*}
	|\tilde{I}_\Omega| & = \left\vert \int e^{i\Phi} \bvphi(b^{-1}2^p \dr\Phi) \bvphi(b^{-1}2^q \dz\Phi) L(r,z) drdz \right\vert \\
	& \lesssim b2^{-p} b2^{-q} \cdot (a^4 2^{-2p}2^{-2q})^{-1} \cdot M \\
	& \leq a^2 2^{-p-q} \cdot a^{-4}2^{2p+2q} \cdot M = t^{-\frac{3}{2}} \kappa^{-3}2^{-\frac{3}{2}k-p-\frac{q}{2}},
\end{align*}
since we can choose $|L| \leq \langle t\k^2 2^{2k+p+q} \cdot 2^{k+p} \rangle^{-1/2} \cdot 2^{k+p} = t^{-\frac{1}{2}} \kappa^{-1} 2^{-\frac{1}{2}k-\frac{q}{2}} =: M$, and we are done.\\

\noindent
\textcircled{2} $a \leq b$ : Solving this in terms of $t$, we have $t^{\frac{1}{2}} \geqslant \kappa^{-1} 2^{-\frac{3}{2}k-p-\frac{q}{2}}$, which then gives $f_r^{-1} = f_z^{-1} = 2^{-k} \leq a$. First condition of \eqref{small-area-cond} is straightforward from the homogeneity of $L$, and hence, applying Theorem \ref{stationary-phase-thm} with $L, M$ the same as above, we get
\begin{equation*}
	|\tilde{I}_\Omega| \lesssim a^{-2}2^{p+q} \cdot M = \k^{-3} 2^{-\frac{3}{2}k - p -\frac{q}{2}} t^{-\frac{3}{2}}.
\end{equation*}

\vspace{5mm}

\subsubsection{High frequency, $\k 2^k \gg 1$.}

When $x_1 \nsim t\k^{-1}2^{-k+p+q}$, from
\begin{equation}
	\dr\Omega = \frac{1}{2}\left(\frac{r}{d_2} - \frac{r}{d_1}\right)
	= \frac{r}{2} \frac{d_1 - d_2}{d_1 d_2} = -\k^{-1} \frac{r z}{d_1 d_2 \Sigma} \sim \k^{-1} \frac{2^{2k+p+q}}{2^{3k}} = \ka^{-1} 2^{-k+p+q},
\end{equation}
we have $|\dr\Phi| \gtrsim t\k^{-1}2^{-k+p+q}$ in this case. Also, along with \eqref{basic-dr-loss},
\begin{equation*}
\begin{split}
	\dr^2\Omega& = \frac{t}{2}\left(\frac{(z+\k^{-1})^2}{d_2^3} - \frac{(z-\k^{-1})^2}{d_1^3}\right) = \frac{t}{2}\left[(z^2+\k^{-2})\left(\frac{1}{d_2^3} - \frac{1}{d_1^3}\right) + 2\ka^{-1}z \left(\frac{1}{d_2^3} + \frac{1}{d_1^3}\right)\right]\\
	& \lesssim t\k^{-1}2^{-2k+q}
	\end{split}
\end{equation*}
since
\begin{equation*}
	\frac{1}{d_2^3} - \frac{1}{d_1^3} = \frac{(d_1 - d_2)(d_1^2 + d_1 d_2 + d_2^2)}{d_1^3 d_2^3} = -\frac{2\Omega (d_1^2 + d_1 d_2 + d_2^2)}{d_1^3 d_2^3} \sim \k^{-1} 2^{-4k+q},
\end{equation*}
so that $\left|\frac{\dr^2\Phi}{\dr\Phi}\right| \lesssim \frac{t\k^{-1}2^{-2k+q}}{t\ka^{-1}2^{-k+p+q}} = 2^{-k-p}$, and
\begin{align*}
	|\dr^3 \Phi| & = \frac{t}{2}\left|\frac{r(z-\k^{-1})^2}{d_1^5} - \frac{r(z+\k^{-1})^2}{d_2^5}\right| \\
	& = tr\left|(z^2+\k^{-2}) \frac{\k^{-1}z (d_1^4 + d_1^3 d_2 + d_1^2 d_2^2 + d_1 d_2^3 + d_2^4)}{d_1^5 d_2^5 \hspace{0.5mm} \Sigma} -\k^{-1}z \left(\frac{1}{d_1^5} + \frac{1}{d_2^5}\right) \right| \\
	& \lesssim t2^{k+p}[ (2^{2k+2q} + \k^{-2})\frac{\k^{-1}2^{k+q} 2^{4k}}{2^{11k}} \k^{-1} 2^{k+q} 2^{-5k} ] \lesssim t\k^{-1} 2^{-3k+p+q},
\end{align*}
so that one can check that from each integration by parts we gain $(t\k^{-1}2^{-k+p+q})^{-1} 2^{-k-p} = t^{-1}\k 2^{-2p-q}$ up to twice. Therefore,
\begin{equation}
	|I_\Omega| \lesssim 2^{3k+2p+q}\cdot \min\{ 1, (t^{-1}\k 2^{-2p-q})^2 \} \leq \k^{\frac{3}{2}} 2^{3k-p-\frac{q}{2}} t^{-\frac{3}{2}}.
\end{equation}
Hence, now we assume $x_1 \sim t\k^{-1}2^{-k+p+q}$. When $2^q \ll 1$, we have the sizes
\begin{gather*}
	\dr^2\Phi \sim t\k^{-1} 2^{-2k+q}, \\
	\dz^2\Phi = t\frac{r^2}{2}\left(\frac{1}{d_2^3} - \frac{1}{d_1^3}\right) \sim t\k^{-1} 2^{-2k+q}, \\
	\dr\dz\Phi = -\frac{3tr}{2}\left(\frac{z+\ka^{-1}}{d_2^3} - \frac{z-\k^{-1}}{d_1^3}\right)
	= -\frac{3tr}{2}\left(z\left[\frac{1}{d_2^3} - \frac{1}{d_1^3}\right] + \k^{-1}\left[\frac{1}{d_2^3} + \frac{1}{d_1^3}\right]\right) \sim t\k^{-1} 2^{-2k},
\end{gather*}
which satisfy \eqref{nondiag-dom} as the 3rd order derivatives also have bounds
\begin{equation}
	\dr\dz^2\Phi = \frac{2}{r}\dz^2\Phi + \frac{3tr^3}{2}\left(\frac{1}{d_1^5} - \frac{1}{d_2^5}\right) \lesssim 2^{-k} |\dz^2\Phi| \ll 2^{-k-q} |\dr\dz\Phi|,
\end{equation}
\begin{align*}
	\dz\dr^2\Omega & = \frac{1}{2}\left( \frac{2(z+\k^{-1})}{d_2^3} - \frac{2(z-\k^{-1})}{d_1^3} - \frac{3(z+\k^{-1})^3}{d_2^5} + \frac{3(z-\k^{-1})^3}{d_1^5} \right) \\
	& = z\left( \frac{1}{d_2^3} - \frac{1}{d_1^3} \right)
	+ \k^{-1}\left( \frac{1}{d_2^3} + \frac{1}{d_1^3} \right) 
	- 3z(z^2+3\k^{-2})\left( \frac{1}{d_2^5} - \frac{1}{d_1^5} \right) 
	- 3\k^{-1}(3z^2+\k^{-2})\left( \frac{1}{d_2^5} + \frac{1}{d_1^5} \right) \\
	& \lesssim \k^{-1} 2^{-3k} \lesssim 2^{-k} |\dr\dz\Omega|.
\end{align*}
As in the low frequency case, the estimates above enables us to use $(\dr, 2^q\dz)$. When $(t\k^{-1}2^{-2k+q})^{1/2} \leq 2^{-k}$, simple set size estimate is already enough to give $|I_\Omega| \lesssim 2^{3k+q} \leq \k^{\frac{3}{2}} 2^{3k-p-\frac{q}{2}} t^{-\frac{3}{2}}$. Hence, assuming $(t\k^{-1}2^{-2k+q})^{1/2} \geq 2^{-k} = f_r^{-1} = f_z^{-1}$, condition \eqref{small-area-cond} has been met, and we can use Theorem \ref{stationary-phase-thm} with $L(r,z) = \vphi(2^{-k}|r,z|) \vphi(2^{-k}r) \vphi(2^{-k-q}z) H(x_1 r) r$ and $M = \k^{\frac{1}{2}} 2^{k-\frac{q}{2}} t^{-\frac{1}{2}}$, which gives
\begin{equation}
	|I_\Omega| \lesssim (t\k^{-1} 2^{-2k+p})^{-1} \k^{\frac{1}{2}} 2^{k-\frac{q}{2}} t^{-\frac{1}{2}} = \k^{\frac{3}{2}} 2^{3k-p-\frac{q}{2}} t^{-\frac{3}{2}}.
\end{equation}
When $2^q \sim 1$, we use the spherical coordinates. We first consider
\begin{equation}
	I_{\Omega,n} := \int e^{i\Phi} \vphi(t^{-1}\k 2^{k-2p-q-n} \drho\Phi)
	\vphi(2^{-k}\rho) \vphi(2^{-p} \sin\phi) \vphi(2^{-q}\cos\phi) H(x_1 \rho\sin\phi) \rho^2 \sin\phi \hspace{1mm} d\rho d\phi.
\end{equation}
By integration by parts along $\drho$,
\begin{align*}
	|I_{\Omega,n}| = & \hspace{1mm} t^{-1}\k 2^{k-2p-q-n} \\
	& \quad \cdot \left\vert \int e^{i\Phi} \drho\left[ \tvphi(t^{-1}\k 2^{k-2p-q-n} \drho\Phi) \vphi(2^{-k}\rho) \vphi(2^{-p}\sin\phi) \vphi(2^{-q}\cos\phi) H(x_1 \rho\sin\phi) \rho^2 \sin\phi \right] d\rho d\phi \right\vert \\
	\leq & \hspace{1mm} (t^{-1}\k 2^{k-2p-q-n})^2 \int |\drho^2\Phi| \tvphi'(t^{-1}\k 2^{k-2p-q-n} \drho\Phi) \vphi(2^{-p}\sin\phi) \vphi(2^{-q}\cos\phi) \langle x_1 2^{k+p} \rangle^{-\frac{1}{2}} 2^{2k+p} d\rho d\phi \\
	& + t^{-1}\k 2^{k-2p-q-n} \int |\drho\left[ \vphi(2^{-k}\rho) \vphi(2^{-p}\sin\phi) \vphi(2^{-q}\cos\phi) H(x_1 \rho\sin\phi) \rho^2 \sin\phi \right]| d\rho d\phi \\
	\lesssim & \hspace{1mm} \k^{\frac{3}{2}} 2^{3k-p-\frac{q}{2}} 2^{-n} t^{-\frac{3}{2}},
\end{align*}
giving acceptable contribution by summing over $n$ for $2^n \gtrsim 1$. Similarly,
\begin{align*}
	|I_\Omega^m| := & \hspace{1mm} \left\vert \int e^{i\Phi} \vphi(t^{-1}\k 2^{-p-m}\dphi\Phi) \vphi(2^{-k}\rho) \vphi(2^{-p}\sin\phi) \vphi(2^{-q}\cos\phi) H(x_1\rho\sin\phi) \rho^2 \sin\phi d\rho d\phi \right\vert \\
	\leq & \hspace{1mm} t^{-1}\k 2^{-p-m} \\
	& \quad \cdot \int \left\vert \dphi \left[ \vphi(t^{-1}\k 2^{-p-m}\dphi\Phi) \vphi(2^{-k}\rho) \vphi(2^{-p}\sin\phi) \vphi(2^{-q}\cos\phi) H(x_1\rho\sin\phi) \rho^2 \sin\phi \right] \right\vert d\rho d\phi \\
	\lesssim & \hspace{1mm} t^{-1}\k 2^{-p-m} \langle t\k^{-1} 2^{2p+q} \rangle^{-\frac{1}{2}} 2^{3k+p} \leq \k^{\frac{3}{2}} 2^{3k - p - \frac{q}{2}} t^{-\frac{3}{2}} 2^{-m},
\end{align*}
and again we have acceptable contribution for $2^m \gtrsim 1$. Therefore, now we have the sizes
\begin{gather*}
	\drho^2\Phi = t\frac{r^2}{2\k^2|\xi|^2}\left(\frac{1}{d_2^3} - \frac{1}{d_1^3}\right) \sim t\k^{-3} 2^{-4k+2p+q}, \\
	\dphi^2\Phi = -\frac{t}{2}\left[ \frac{\rho^2\Omega}{d_1 d_2} + \k^{-2}r^2 \left(\frac{1}{d_2^3} - \frac{1}{d_1^3}\right) \right] - \rho\drho\Phi \sim t\k^{-1} 2^q \\
	\drho\dphi\Phi = \frac{t}{\rho} \frac{\k^{-1}r}{2} \left( \frac{1}{d_2} + \frac{1}{d_1} - \frac{\k^{-1}(z+\k^{-1})}{d_2^3} + \frac{\k^{-1}(z-\k^{-1})}{d_1^3} \right) + \frac{1}{\rho}\dphi\Phi \sim t\k^{-1} 2^{-k+p},	
\end{gather*}
where the second line holds due to $2^q \sim 1$ (so actually the size is simply $t\k^{-1}$). Using $(2^k\drho, 2^p\dphi)$ instead of $(\drho,\dphi)$ similarly as previous cases shows that this satisfy \eqref{nondiag-dom} along with
\begin{equation*}
	2^k \dphi\drho^2\Phi = 2^k \frac{2z}{r}\drho^2\Phi + 2^k \frac{tr^2}{\k^2 \rho^2} \left( \frac{3\k^{-1}r}{d_2^5} + \frac{3\k^{-1}r}{d_1^5} \right) \lesssim 2^{k-p} |\drho^2\Phi| \lesssim |\drho\dphi\Phi|,
\end{equation*}
\begin{align*}
	2^p \drho\dphi^2\Phi = -\frac{t\k^{-1} 2^p}{2\rho} & \left( -\frac{\k^{-1}r^2}{d_2^3} + 3\frac{\k^{-1}r^2(z+\k^{-1})^2}{d_2^5} - z\k^{-1}\frac{2(z+\k^{-1})d_2^2 - 3(z+\k^{-1})^3}{d_2^5} \right. \\
	& \left. + \frac{\k^{-1}r^2}{d_1^3} - 3\frac{\k^{-1}r^2(z-\k^{-1})^2}{d_1^5}
	+ z\k^{-1} \frac{2(z-\k^{-1})d_1^2 - 3(z-\k^{-1})^3}{d_1^5} \right) \\
	\lesssim t\k^{-1} 2^{-k+p} & \lesssim 2^{-k+p} |\dphi^2\Phi| \lesssim |\drho\dphi\Phi|.
\end{align*}
Thus, proceeding similarly as in the low frequency case to use Theorem \ref{stationary-phase-thm} with $M = 2^{2k-\frac{q}{2}} \k^{\frac{1}{2}} t^{-\frac{1}{2}}$ and $L(\rho,\phi) = \vphi(2^{-k}\rho) \vphi(2^{-p} \sin\phi) \vphi(2^{-q}\cos\phi) H(x_1 \rho\sin\phi) \rho^2 \sin\phi$, we have
\begin{equation}
	|I_\Omega| \lesssim (t\k^{-1} 2^{-k+p})^{-1} \cdot 2^{2k-\frac{q}{2}} \k^{\frac{1}{2}} t^{-\frac{1}{2}} = 2^{3k-p-\frac{q}{2}} \k^{\frac{3}{2}} t^{-\frac{3}{2}}.
\end{equation}

\vspace{5mm}

\subsubsection{Medium frequency} The situations are subtly different according to the regions, and we work in three different regimes : $2^q \sim 1, p \leq l \leq 0$ ; $2^q \sim 1, l \leq p$ ; and $2^q \ll 1$, where $l$ is introduced again for $z-\k^{-1} \sim 2^{k+l}$ as in case 5.1.3. This implies we will only be working on frequencies on the upper half plane.\\

\noindent \textcircled{1} $2^q \sim 1, p \leq l \leq 0$ : Here, we work with the integral
\begin{equation}
	I_{\Omega,l} := \int e^{i\Phi} \vphi(2^{-k}|r,z|) \vphi(2^{-k-p}r) \vphi(2^{-k}z) \vphi(2^{-k-l}(z-\k^{-1})) H(x_1 r) r dr dz.
\end{equation}
When $x_1 \nsim t 2^{p-l}$, $|\dr\Phi| \gtrsim t 2^{p-l}$ since $\dr\Omega = \frac{r}{2}\left(\frac{1}{d_2} - \frac{1}{d_1}\right) = -\frac{r\Omega}{d_1 d_2} \sim 2^{p+q-l}$, and
\begin{align*}
	|\dr^2\Phi| & = \frac{t}{2}\left\vert \frac{(z+\k^{-1})^2}{d_2^3} - \frac{(z-\k^{-1})^2}{d_1^3} \right\vert \lesssim t 2^{-k-l}, \\
	|\dr^3\Phi| & = \frac{3t}{2}\left\vert \frac{r(z+\k^{-1})^2}{d_2^5} - \frac{r(z-\k^{-1})^2}{d_1^5} \right\vert \lesssim t 2^{-2k-2l},
\end{align*}
we have $\frac{|\dr^2\Phi|}{|\dr\Phi|} \lesssim 2^{-k-p}$ and $\frac{|\dr^3\Phi|}{|\dr\Phi|} \lesssim 2^{-2k-p-l} \leq 2^{-2k-2p}$. Hence,
\begin{equation}
	|I_{\Omega,l}| \lesssim 2^{3k+2p+l} \cdot \min\{1, (t 2^{p-l} \cdot 2^{k+p})^{-2}\} \leq 2^{\frac{3}{2}k - p + \frac{5}{2}l} t^{-\frac{3}{2}}.
\end{equation}
So we assume $x_1 \sim t 2^{p-l}$ and now change to spherical coordinates centered at $(0, \k^{-1})$, denoted by $(\rho,\phi)$ again as in the analysis of medium frequency of $\Sigma$. From
\begin{align*}
	|I_{\Omega,l,n}| & := \left\vert \int e^{i\Phi} \vphi(t^{-1}2^{-n} \drho\Phi) \vphi(2^{-k-p}\rho\sin\phi) \vphi(2^{-k-l}\rho\cos\phi) H(x_1 \rho\sin\phi) \rho^2\sin\phi d\rho d\phi \right\vert \\
	& \leq t^{-1}2^{-n} \int \left\vert \drho \left[ \vphi(t^{-1}2^{-n} \drho\Phi) \vphi(2^{-k-p}\rho\sin\phi) \vphi(2^{-k-l}\rho\cos\phi) H(x_1 \rho\sin\phi) \rho^2\sin\phi \right] d\rho d\phi \right\vert \\
	& \lesssim t^{-1} 2^{-n} ( 2^{-k-l} 2^{3k+2p+l} + 2^{2k+2l} ) \langle t 2^{k+2p-l} \rangle^{-\frac{1}{2}} \leq 2^{\frac{3}{2}k - p + \frac{5}{2}l} 2^{-n} t^{-\frac{3}{2}},
\end{align*}
we get acceptable contributions for summations over $n$. $\vphi(2^{-k}|r,z|)$ and $\vphi(2^{-k}z)$ were omitted for notational convenience and will be dropped from hereon, because
\begin{gather*}
	\bigg| \drho[\vphi(2^{-k}|r,z|), \hspace{1mm} \vphi(2^{-k}z)] \bigg| \lesssim \hspace{1mm} \frac{z}{\rho} 2^{-k}, \hspace{1mm} \frac{|z-\k^{-1}|}{\rho} 2^{-k} \leq 2^{-k}, \\
	\bigg| \dphi[\vphi(2^{-k}|r,z|), \hspace{1mm} \vphi(2^{-k}z)] \bigg| \lesssim \hspace{1mm} 0, \hspace{1mm} r 2^{-k} \lesssim 2^p,
\end{gather*}
and hence are terms of smaller sizes. Also from
\begin{align*}
	|I_{\Omega,l}^m| & := \left\vert \int e^{i\Phi} \vphi(t^{-1}2^{-k-p-\frac{1}{2}l-m} \dphi\Phi) \vphi(2^{-k-p}\rho\sin\phi) \vphi(2^{-k-l}\rho\cos\phi) H(x_1 \rho\sin\phi) \rho^2\sin\phi d\rho d\phi \right\vert \\
	& \leq t^{-1}2^{-k-p-\frac{1}{2}l-m} \\
	& \quad \cdot \int \left\vert \dphi \left[ \vphi(t^{-1}2^{-k-p-\frac{1}{2}l-m} \dphi\Phi) \vphi(2^{-k-p}\rho\sin\phi) \vphi(2^{-k-l}\rho\cos\phi) H(x_1 \rho\sin\phi) \rho^2\sin\phi \right] \right\vert d\rho d\phi \\
	& \lesssim t^{-1}2^{-k-p-\frac{1}{2}l-m} ( 2^{l-p} 2^{3k+2p+l} + 2^{3k+p+2l} ) \langle t 2^{k+2p-l} \rangle^{-\frac{1}{2}} \leq 2^{\frac{3}{2}k - p + 2l} 2^{-m} t^{-\frac{3}{2}},
\end{align*}
we again get acceptable contributions for summations over $m$. Therefore, now we are in a region where $\drho\Phi \ll t$ and $\dphi\Phi \ll t 2^{k+p+\frac{1}{2}l}$, which yields
\begin{equation}
	\drho^2\Phi \sim t 2^{-k+2p-2l}, \quad
	\dphi^2\Phi \sim t 2^{k+l}, \quad
\end{equation}
from
\begin{equation}\label{derivatives-La0-medfreq}
\begin{split}
	\drho^2\Phi & = \frac{t}{2}\left( \frac{1}{d_2} - \frac{(r^2+z^2-\k^{-2})^2}{d_1^2 d_2^3} \right) = 2t \frac{\k^{-2}r^2}{d_1^2 d_2^3}, \\
	\rho\drho\dphi\Phi - \dphi\Phi & = t\k^{-1}r \frac{r^2 + (z+\k^{-1})(z-\k^{-1})}{d_2^3}, \\
	\dphi^2\Phi + \rho\drho\Phi & = \frac{t}{2}\left( -d_1 + \frac{r^2+z^2+\k^{-2}}{d_2} - \frac{2\k^{-1}(z-\k^{-1})}{d_2} - \frac{4(\k^{-1}r)^2}{d_2^3} \right) \\
	& = \frac{t}{2} \left( -d_1 + \frac{d_1^2}{d_2} - 4\frac{(\k^{-1}r)^2}{d_2^3} \right) \\
	& = -t \left( \frac{d_1}{d_2}\Omega + 2\frac{(\k^{-1}r)^2}{d_2^3} \right),
\end{split}
\end{equation}
since this shows that $|\rho\drho\dphi\Phi - \dphi\Phi| \lesssim t 2^{k+p+l}$ and $|\dphi^2\Phi + \rho\drho\Phi| \sim t 2^{k+l}$. Last conditions of \eqref{diag-dom} can be checked by
\begin{align*}
	\dphi\drho^2\Phi & = \drho^2\Phi \cdot \frac{2}{r d_2^2}[d_2^2(z-\k^{-1}) + 3\k^{-1}r^2] \lesssim 2^{l-p} |\drho^2\Phi|, \\
	\frac{1}{t}[\rho\drho\dphi^2\Phi - \dphi^2\Phi] & = \left( \frac{\k^{-1}(z-\k^{-1})}{d_2^2} + 6\frac{\k^{-2}r^2}{d_2^4} \right) \frac{r^2+z^2-\k^{-2}}{d_2} - \frac{2\k^{-2}r^2}{d_2^3} \lesssim 2^{k+2l} \lesssim \frac{1}{t} |\dphi^2\Phi|
\end{align*}
Therefore, with $L(\rho,\phi) = \vphi(2^{-k-p}\rho\sin\phi) \vphi(2^{-k-l}\rho\cos\phi) H(x_1 \rho\sin\phi) \rho^2\sin\phi$ and $M = 2^{\frac{3}{2}k + \frac{3}{2}l} t^{-\frac{1}{2}}$ gives us
\begin{equation}
	|I_{\Omega,l}| \lesssim (t 2^{p-\frac{1}{2}l})^{-1} \cdot 2^{\frac{3}{2}k + \frac{3}{2}l} t^{-\frac{1}{2}} = 2^{\frac{3}{2}k - p + 2l} t^{-\frac{3}{2}}.
\end{equation}
This is the weakest estimate among the bounds obtained above, and therefore the final estimate.\\

\noindent
\textcircled{2} $2^q \sim 1, l \leq p$ : As in the analysis of $\Sigma$, we don't localize in $l$ individually anymore, but consider
\begin{equation}
	I_{\Omega, \leq p} := \int e^{i\Phi} \vphi(2^{-k}|r,z|) \vphi(2^{-k-p}r) \vphi(2^{-k}z) \bvphi(2^{-k-p}(z-\k^{-1})) H(x_1 r) r dr dz.
\end{equation}
Important change from case \textcircled{1} is that now $d_1 \sim 2^{k+p}$. If $x_1 \nsim t$,
\begin{gather*}
	|\dr^2\Phi| = \frac{t}{2}\left\vert \frac{(z+\k^{-1})^2}{d_2^3} - \frac{(z-\k^{-1})^2}{d_1^3} \right\vert \lesssim t 2^{-k-p}, \\
	|\dr^3\Phi| = \frac{3t}{2}\left\vert \frac{r(z+\k^{-1})^2}{d_2^5} - \frac{r(z-\k^{-1})^2}{d_1^5} \right\vert \lesssim t 2^{-2k-2p},
\end{gather*}
and hence $\frac{|\dr^2\Phi|}{|\dr\Phi|} \lesssim 2^{-k-p}$, $\frac{|\dr^3\Phi|}{|\dr\Phi|} \lesssim 2^{-2k-2p}$. With the other usual loss from $\dr$ derivative, this gives
\begin{equation}
	|I_{\Omega, \leq p}| \lesssim 2^{3k+3p} \cdot \min\{ 1, (t2^{k+p})^{-2}\} \leq 2^{\frac{3}{2}k + \frac{3}{2}p} t^{-\frac{3}{2}}.
\end{equation}
Now we assume $x_1 \sim t$ and start using polar coordinates around $(0,\k^{-1})$ again. From
\begin{align*}
	|I_{\Omega, \leq p, n}| & = \left\vert \int e^{i\Phi} \vphi(t^{-1}2^{-n} \drho\Phi) \vphi(2^{-k-p}\rho\sin\phi) \bvphi(2^{-k-p}\rho\cos\phi) H(x_1 \rho\sin\phi) \rho^2\sin\phi d\rho d\phi \right\vert \\
	& = t^{-1}2^{-n} \int \left\vert \drho \left[ \vphi(t^{-1}2^{-n} \drho\Phi) \vphi(2^{-k-p}\rho\sin\phi) \bvphi(2^{-k-p}\rho\cos\phi) H(x_1 \rho\sin\phi) \rho^2\sin\phi \right] d\rho d\phi \right\vert \\
	& \lesssim t^{-1} 2^{-n} ( 2^{-k-p} 2^{3k+3p} + 2^{2k+2p} ) \langle t 2^{k+p} \rangle^{-\frac{1}{2}} \leq 2^{\frac{3}{2}k + \frac{3}{2}p} 2^{-n} t^{-\frac{3}{2}},
\end{align*}
summation over $2^n \gtrsim 1$ gives acceptable bound, and hence we assume $\drho\Phi \ll t$ now. Also,
\begin{align*}
	|I_{\Omega,\leq p}^m| & = \left\vert \int e^{i\Phi} \vphi(t^{-1}2^{-k-\frac{3}{2}p-m} \dphi\Phi) \vphi(2^{-k-p}\rho\sin\phi) \vphi(2^{-k-p}\rho\cos\phi) H(x_1 \rho\sin\phi) \rho^2\sin\phi d\rho d\phi \right\vert \\
	& \leq t^{-1}2^{-k-\frac{3}{2}p-m} \\
	& \quad \cdot \int \left\vert \dphi \left[ \vphi(t^{-1}2^{-k-\frac{3}{2}p-m} \dphi\Phi) \vphi(2^{-k-p}\rho\sin\phi) \vphi(2^{-k-p}\rho\cos\phi) H(x_1 \rho\sin\phi) \rho^2\sin\phi \right] \right\vert d\rho d\phi \\
	& \lesssim t^{-1}2^{-k-\frac{3}{2}p-m} \cdot 2^{3k+3p} \langle t 2^{k+p} \rangle^{-\frac{1}{2}} \leq 2^{\frac{3}{2}k + p} 2^{-m} t^{-\frac{3}{2}},
\end{align*}
so by the same reasoning, we now assume $\dphi\Phi \ll t 2^{k+\frac{3}{2}p}$. These, together with \eqref{derivatives-La0-medfreq} gives
\begin{equation}
	\drho^2\Phi \sim t 2^{-k}, \quad
	\dphi^2\Phi \sim t 2^{k+p}.
\end{equation}
Also, from
\begin{gather*}
	\dphi\drho^2\Phi = \drho^2\Phi \cdot \frac{2}{r d_2^2}[d_2^2(z-\k^{-1}) + 3\k^{-1}r^2] \lesssim |\drho^2\Phi|, \\
	\frac{1}{t}[\rho\drho\dphi^2\Phi - \dphi^2\Phi] = \left( \frac{\k^{-1}(z-\k^{-1})}{d_2^2} + 6\frac{\k^{-2}r^2}{d_2^4} \right) \frac{r^2+z^2-\k^{-2}}{d_2} - \frac{2\k^{-2}r^2}{d_2^3} \lesssim 2^{k+2p} \lesssim \frac{1}{t} |\dphi^2\Phi|,
\end{gather*}
the phase and $L(\rho,\phi) = \vphi(2^{-k-p}\rho\sin\phi) \vphi(2^{-k-p}\rho\cos\phi) H(x_1 \rho\sin\phi) \rho^2\sin\phi$ with $M = 2^{\frac{3}{2}k + \frac{3}{2}p} t^{-\frac{1}{2}}$ satisfy \eqref{diag-dom}, and thus
\begin{equation}
	|I_{\Omega, \leq p}| \lesssim t^{-1} 2^{-\frac{1}{2}p} \cdot 2^{\frac{3}{2}k + \frac{3}{2}p} t^{-\frac{1}{2}} = 2^{\frac{3}{2}k + p} t^{-\frac{3}{2}}.
\end{equation}
Again, other bounds are as strong as this last one, and hence this is the final estimate for case \textcircled{2}.\\

\noindent
\textcircled{3} $2^q \ll 1$ : When $x_1 \nsim t 2^q$, we have $|\dr\Phi| \gtrsim t 2^q$,
\begin{align*}
	|\dr^2\Phi| & = \frac{t}{2}\left| (z^2+\k^{-2})\left(\frac{1}{d_2^3} - \frac{1}{d_1^3}\right) + 2\k^{-1}z \left(\frac{1}{d_2^3} + \frac{1}{d_1^3}\right) \right| \lesssim t 2^{-k+q-3l}, \\
	|\dr^3\Phi| & = \frac{3rt}{2}\left| (z^2+\k^{-2})\left( \frac{1}{d_2^5} - \frac{1}{d_1^5} \right) + 2\k^{-1}z\left(\frac{1}{d_2^5} + \frac{1}{d_1^5}\right) \right| \lesssim 2^{-2k+p+q-5l}.
\end{align*}
Noting that $p=l=0$ when $2^q \ll 1$, this means $\dr$ derivative causes $2^{-k}$ loss, and therefore, $|I_\Omega| \lesssim 2^{\frac{3}{2}k - \frac{1}{2}q} t^{-\frac{3}{2}}$ via the usual reasoning. Now we assume $x_1 \sim t 2^q$. Second order derivatives satisfy
\begin{equation}
\begin{split}
	\dr^2\Phi & \lesssim t 2^{-k+q}, \\
	\dz^2\Phi & = t\frac{r^2}{2}\left(\frac{1}{d_2^3} - \frac{1}{d_1^3}\right) = t\frac{r^2\Omega(d_1^2 + d_1 d_2 + d_2^2)}{d_1^3 d_2^3} \sim t 2^{-k+q}, \\
	\dr\dz\Phi & = -t\frac{r}{2}\left[ \frac{z+\k^{-1}}{d_2} - \frac{z-\k^{-1}}{d_1} \right] = -t\frac{r}{2}\left[z\left(\frac{1}{d_2^3} - \frac{1}{d_1^3}\right) +\k^{-1}\left(\frac{1}{d_2^3} + \frac{1}{d_1^3}\right)\right] \sim t 2^{-k},
\end{split}
\end{equation}
when $2^q \ll 1$. Along with
\begin{gather*}
	\dr\dz^2\Phi = tr\left(\frac{1}{d_2^3} - \frac{1}{d_1^3}\right) + \frac{3tr^3}{2}\left( \frac{1}{d_1^5} - \frac{1}{d_2^5} \right) \lesssim t 2^{-2k} \lesssim 2^{-k} |\dr\dz\Phi|, \\
	\dz\dr^2\Phi = t\left(\frac{z+\k^{-1}}{d_2^3} - \frac{z-\k^{-1}}{d_1^3}\right) + \frac{3t}{2}\left( \frac{(z-\k^{-1})^3}{d_1^5} - \frac{(z+\k^{-1})^3}{d_2^5} \right) \lesssim t 2^{-2k} \lesssim 2^{-k} |\dr\dz\Phi|,
\end{gather*}
the phase satisfy \eqref{nondiag-dom} and can be made to satisfy \eqref{small-area-cond} by using $(\dr, 2^q\dz)$ and removing small time region similarly as in the low frequency case. Therefore, together with $L(r,z) = \vphi(2^{-k}|r,z|) \vphi(2^{-k-p}r) \vphi(2^{-k}z) H(x_1 r) r$ and $M = 2^{\frac{1}{2}k - \frac{1}{2}q} t^{-\frac{1}{2}}$, we obtain
\begin{equation}
	|I| \lesssim t^{-1} 2^k \cdot 2^{\frac{1}{2}k - \frac{1}{2}q} t^{-\frac{1}{2}} = 2^{\frac{3}{2}k - \frac{1}{2}q} t^{-\frac{3}{2}}
\end{equation}
from Theorem \ref{stationary-phase-thm}.

\section{High speeds and long time -- proof of Theorem \ref{main-theorem} and Theorem \ref{MainThmEC}}\label{SecMainThm}

We conclude by giving the proof of Theorem \ref{main-theorem} and Theorem \ref{MainThmEC} in this section.

\begin{proof}[Proof of Theorem \ref{main-theorem}]
We recall from \eqref{eq:projections} the projection operators $\PP_{\mu\La}$, $\mu \in \{+,-\}$, $\La \in \{\Sigma,\Omega\}$, onto the eigenspaces corresponding to $i\mu \La$. From \eqref{eq:CER2} we obtain that
\begin{equation}\label{eq:projected-system}
	\dt \PP_{\mu\La} \begin{pmatrix}
		\rho \\ u
	\end{pmatrix} = i\mu \La \PP_{\mu\La} \begin{pmatrix}
		\rho \\ u
	\end{pmatrix} - \PP_{\mu\La} \begin{pmatrix}
		u\cdot\nabla\rho+\alpha\rho\div( u) \\ u\cdot\nabla u + \alpha\rho\nabla \rho
	\end{pmatrix}.
\end{equation}
Employing the shorthand $W:=(\rho,u)^\intercal$, Duhamel's formula gives
\begin{equation}\label{Duhamel}
	\begin{aligned}
		\PP_{\mu\La}W(t) & = e^{ it\mu\Lambda} \PP_{\mu\La} W(0) -  \int_0^t e^{ i(t-\tau)\mu\Lambda} \PP_{\mu\La}\mathcal{N}(W, W)(\tau)d\tau,\\
		\mathcal{N}(W, W)&=\begin{pmatrix}
			u\cdot\nabla\rho+\alpha\rho\div( u) \\ u\cdot\nabla u + \alpha\rho\nabla \rho
		\end{pmatrix}.
	\end{aligned}
\end{equation}

The remainder of this proof gives a lower bound on the time of existence by closing a bootstrap for the energy estimates \eqref{energy-estimate}. More precisely, by standard local well-posedness and \eqref{energy-estimate}, it suffices to establish an a priori bound on $\norm{\nabla W}_{L_t^1(0,T; L_x^\infty)}$, although we will bound $\norm{(W, \nabla W)}_{L_t^1(0,T; L_x^\infty)}$ for technical reasons. For notational simplicity, in the remaining part of the proof, we omit the time-region of integration and simply write, e.g. $\Vert \nabla W\Vert_{L^1_tL^\infty_x}$\\
\textcircled{1} High regularity case $m > \frac{7}{2}$: In order to put us in a position where we can invoke the Strichartz estimates from Theorem \ref{Strichartz}, we observe that by H\"older in time and Sobolev embedding in space, for any $(s,R)\in(0,\infty)\times(2,\infty)$ satisfying $\frac{s}{3} - \frac{1}{R} > 0$ and $q>2$ with $\frac{1}{q} + \frac{1}{R} = \frac{1}{2}$ there holds that
\begin{equation}
	\norm{\nabla W}_{L_t^1L_x^\infty} \lesssim \norm{\nabla W}_{L_t^1W_x^{s,R}} \leq T^{1-\frac{1}{q}} \norm{\nabla W}_{L_t^qW_x^{s,R}}.
\end{equation}
Moreover, by Proposition \ref{prop:Lp-equivalence} and \eqref{Duhamel} we have that
\begin{equation}\label{eq:Duhamel-used}
\begin{aligned}
	\norm{\nabla W}_{L_t^qW_x^{s,R}}&\lesssim \sum_{\substack{\mu\in\{-,+\},\\\La\in\{\Sigma,\Omega\}}}\norm{\PP_{\mu\La}\nabla W}_{L_t^qW_x^{s,R}}\\
	&\lesssim \sum_{\substack{\mu\in\{-,+\},\\\La\in\{\Sigma,\Omega\}}}\norm{e^{ it\mu\Lambda} \PP_{\mu\La} \nabla W(0)}_{L_t^qW_x^{s,R}}+\norm{\int_0^t e^{ i(t-\tau)\mu\Lambda} \nabla\PP_{\mu\La}\mathcal{N}(W, W)(\tau)d\tau}_{L_t^qW_x^{s,R}}.
\end{aligned}    
\end{equation}

To bound the linear terms in \eqref{eq:Duhamel-used}, we note that by the Strichartz estimates \eqref{eq:L2Strichartz} there holds that for $k\in\Z$
\begin{equation}
	\norm{P_k e^{ it\mu\Lambda} \PP_{\mu\La} \nabla W(0)}_{L_t^q W_x^{s,R}} \lesssim 2^{(1+s)k^+} (\langle 2^k \ka\rangle^3 \varep\k^{-3})^{\frac{1}{q}} \norm{P_k \PP_{\mu\La} W(0)}_{L^2}, 
\end{equation}
where $k^+ = \max\{0,k\}$. Thus,
\begin{equation}\label{eq:linear_bd_with_gradient}
\begin{aligned}    
	\norm{e^{ it\mu\Lambda} \PP_{\mu\La} \nabla W(0)}_{L_t^q W_x^{s,R}}
	& \lesssim \left\{\sum_{k \in \mathbb{Z}} \left[ 2^{(1+s)k^+} \varep^{\frac{1}{q}} \max\{2^{3k}, \ka^{-3}\}^{\frac{1}{q}} \norm{P_k \PP_{\mu\La} W(0)}_{L^2}\right]^2 \right\}^{1/2} \\
	& \lesssim \varepsilon^{\frac{1}{q}} \norm{\PP_{\mu\La} W(0)}_{H^{1+s+\frac{3}{q}}} + \left( \frac{\varepsilon}{\ka^3} \right)^{\frac{1}{q}} \norm{\PP_{\mu\La} W(0)}_{H^{1+s}}.
\end{aligned}    
\end{equation}
Likewise,
\begin{equation}\label{eq:linear_bd_wo_gradient_localized}
	\norm{P_k e^{ it\mu\Lambda} \PP_{\mu\La} W(0)}_{L_t^q W_x^{s,R}} \lesssim 2^{s k^+} (\langle 2^k \ka\rangle^3 \varep\k^{-3})^{\frac{1}{q}} \norm{P_k \PP_{\mu\La} W(0)}_{L^2}, 
\end{equation}
which results in
\begin{equation}\label{eq:linear_bd_wo_gradient}
\begin{aligned}    
	\norm{e^{ it\mu\Lambda} \PP_{\mu\La} W(0)}_{L_t^q W_x^{s,R}}
	& \lesssim \left\{\sum_{k \in \mathbb{Z}} \left[ 2^{s k^+} \varep^{\frac{1}{q}} \max\{2^{3k}, \ka^{-3}\}^{\frac{1}{q}} \norm{P_k \PP_{\mu\La} W(0)}_{L^2}\right]^2 \right\}^{1/2} \\
	& \lesssim \varepsilon^{\frac{1}{q}} \norm{\PP_{\mu\La} W(0)}_{H^{s+\frac{3}{q}}} + \left( \frac{\varepsilon}{\ka^3} \right)^{\frac{1}{q}} \norm{\PP_{\mu\La} W(0)}_{H^s}.
\end{aligned}    
\end{equation}

As for the nonlinear terms in \eqref{eq:Duhamel-used}, by the Strichartz estimates \eqref{eq:dualStrichartz} there holds that
\begin{equation}
\begin{aligned}
	&\left\Vert \int_0^t e^{ i(t-\tau)\mu\Lambda} \nabla \right.  \PP_{\mu\La} \mathcal{N}(W, W)(\tau)d\tau \bigg\Vert_{L_t^qW_x^{s,R}}\\
	&\quad \lesssim
	\left\{ \sum_{k \in \mathbb{Z}} \left[ 2^{(1+s)k^+} \varepsilon^{\frac{1}{q}} \max\{2^{3k}, \ka^{-3}\}^{\frac{1}{q}} \Vert P_k \PP_{\mu\La} \mathcal{N}(W, W) \Vert_{L_t^1 L_x^2} \right]^2 \right\}^{1/2} \\
	&\quad  \lesssim \varepsilon^{\frac{1}{q}} \Vert \mathcal{N}(W, W) \Vert_{L_t^1 H_x^{1+s+\frac{3}{q}}}
	+ \left( \frac{\varepsilon}{\ka^3} \right)^{\frac{1}{q}} \Vert \mathcal{N}(W, W) \Vert_{L_t^1 H_x^{1+s}} \\
	&\quad\lesssim \varepsilon^{\frac{1}{q}}[ \Vert W \Vert_{L_t^1 L_x^\infty} \Vert \nabla W \Vert_{L_t^{\infty} H_x^{1+s+\frac{3}{q}}} + \Vert W \Vert_{L_t^\infty H_x^{1+s+\frac{3}{q}}} \Vert \nabla W \Vert_{L_t^1 L_x^\infty} ] \\
	& \quad\qquad+ \left( \frac{\varepsilon}{\ka^3} \right)^{\frac{1}{q}} [\Vert W \Vert_{L_t^1 L_x^\infty} \Vert \nabla W \Vert_{L_t^{\infty} H_x^{1+s}} + \Vert W \Vert_{L_t^{\infty} H_x^{1+s}} \Vert \nabla W \Vert_{L_t^1 L_x^\infty} ],
\end{aligned}
\end{equation}
since the nonlinear terms are all of the form $f\cdot\nabla g$ with $f,g \in \{\rho, u\}$, where $\nabla$ and multiplication are understood in a suitable sense according to their dimensions. Similarly,
\begin{equation}\label{eq:nonlinear_bd_wo_grad}
\begin{aligned}
	\left\Vert \int_0^t e^{ i(t-\tau)\mu\Lambda} \right. \PP_{\mu\La}  \mathcal{N}(W, W)(\tau)d\tau \bigg\Vert_{L_t^qW_x^{s,R}}
	&\lesssim\left\{ \sum_{k \in \mathbb{Z}} \left[ 2^{s k^+} \varepsilon^{\frac{1}{q}} \max\{2^{3k}, \ka^{-3}\}^{\frac{1}{q}} \Vert P_k \PP_{\mu\La} \mathcal{N}(W, W) \Vert_{L_t^1 L_x^2} \right]^2 \right\}^{1/2} \\
	&\lesssim \varepsilon^{\frac{1}{q}} \Vert \mathcal{N}(W, W) \Vert_{L_t^1 H_x^{s+\frac{3}{q}}}
	+ \left( \frac{\varepsilon}{\ka^3} \right)^{\frac{1}{q}} \Vert \mathcal{N}(W, W) \Vert_{L_t^1 H_x^{s}} \\
	 &\lesssim\varepsilon^{\frac{1}{q}}[ \Vert W \Vert_{L_t^1 L_x^\infty} \Vert \nabla W \Vert_{L_t^{\infty} H_x^{s+\frac{3}{q}}} + \Vert W \Vert_{L_t^\infty H_x^{s+\frac{3}{q}}} \Vert \nabla W \Vert_{L_t^1 L_x^\infty} ] \\
	&\quad + \left( \frac{\varepsilon}{\ka^3} \right)^{\frac{1}{q}} [\Vert W \Vert_{L_t^1 L_x^\infty} \Vert \nabla W \Vert_{L_t^{\infty} H_x^{s}} + \Vert W \Vert_{L_t^{\infty} H_x^{s}} \Vert \nabla W \Vert_{L_t^1 L_x^\infty} ].
\end{aligned}
\end{equation}
As a result,
\begin{equation}\label{eq:nonlinear_bd}
\begin{split}
	& \left\Vert \int_0^t e^{ i(t-\tau)\mu\Lambda} \PP_{\mu\La} \mathcal{N}(W, W)(\tau)d\tau \right\Vert_{L_t^qW_x^{s,R}} + \left\Vert \int_0^t e^{ i(t-\tau)\mu\Lambda} \nabla \PP_{\mu\La} \mathcal{N}(W, W)(\tau)d\tau \right\Vert_{L_t^qW_x^{s,R}} \\
	&\qquad \lesssim \varepsilon^{\frac{1}{q}} \norm{W}_{L_t^1W_x^{1,\infty}} \norm{W}_{L_t^\infty H_x^{2+s+\frac{3}{q}}} + (\varep\ka^{-3})^{\frac{1}{q}} \norm{W}_{L_t^1W_x^{1,\infty}} \norm{W}_{L_t^\infty H_x^{2+s}}.
\end{split}
\end{equation}

In conclusion, by combining the bounds \eqref{eq:linear_bd_with_gradient}, \eqref{eq:linear_bd_wo_gradient}, and \eqref{eq:nonlinear_bd} with the energy estimates \eqref{energy-estimate}, we thus have for any $m>0$ that
\begin{equation}\label{eq:summary}
\begin{aligned}
	\norm{W}_{L^\infty_tH^m_x} \leq & \norm{W(0)}_{H^m_x} \exp\left(K\norm{\nabla W}_{{L^1_t L^{\infty}_x}}\right),\\
	\norm{W}_{L_t^1W_x^{1,\infty}} \leq & C\Big[(T^{q-1} \varepsilon)^{\frac{1}{q}} \norm{W(0)}_{H_x^{1+s+\frac{3}{q}}}
	+ (T^{q-1} \varepsilon\ka^{-3})^{\frac{1}{q}} \norm{W(0)}_{H_x^{1+s}} \\
	& \qquad + (T^{q-1} \varepsilon)^{\frac{1}{q}} \norm{W}_{L_t^1 W_x^{1,\infty}} \norm{W}_{L_t^\infty H_x^{2+s+\frac{3}{q}}} \\
	& \qquad + (T^{q-1} \varep\ka^{-3})^{\frac{1}{q}} \norm{W}_{L_t^1 W_x^{1,\infty}} \norm{W}_{L_t^\infty H_x^{2+s}}\Big],
\end{aligned}
\end{equation}
where $K = K(m)$ and $C = C(q,s)$.

For $m>\frac{7}{2}$ given, we set $s = m-2-\frac{3}{q}$. Note that since $s > \frac{3}{R} = \frac{3}{2} - \frac{3}{q}$, $m = 2 + s + \frac{3}{q}$ can match any number greater than $\frac{7}{2}$, by choosing appropriate $s$. Then, until the time given by the assumption \eqref{eq:main-assumption-high-reg} of Theorem \ref{main-theorem}
\begin{equation}
	T = M  \eps^{-\frac{1}{q-1}}\min\{1,(c\eps)^{^{\frac{3}{q-1}}}\}\norm{W(0)}_{H^m}^{-\frac{q}{q-1}}.
\end{equation}
we obtain from \eqref{eq:summary} that
\begin{equation}
	\norm{W}_{L^\infty_tH^m_x} \leq 2 e^K \norm{W(0)}_{H^m},\qquad
	\norm{W}_{L_t^1W_x^{1,\infty}} \leq 1,
\end{equation}
provided that $M=M(m,q)>0$ is small enough. This closes the bootstrap and finishes the proof of Theorem \ref{main-theorem} for $m > \frac{7}{2}$.\\

\noindent
\textcircled{2} Low regularity case $ \frac{5}{2} < m < \frac{7}{2}$: Here, we control the high frequency of $W$ differently to prevent the loss of derivatives. We first note that we may assume $T \gtrsim \norm{W(0)}_{H^m}^{-1}$ from the local well-posedness. Hence, if $\norm{W(0)}_{H^m} \gtrsim \eps^{-1}$, we are done when $\kappa \gtrsim 1$ since
\begin{equation*}
	\eps^{-\frac{ m-\frac{5}{2} }{ q - (m-\frac{5}{2}) }} \norm{W(0)}_{H^m}^{-\frac{q}{ q - (m-\frac{5}{2}) }} \lesssim \norm{W(0)}_{H^m}^{-1},
\end{equation*}
so that the theorem is true by simple local theory. Likewise, if $\norm{W(0)}_{H^m} \gtrsim \eps^2 c^3$ when $\ka \lesssim 1$,
\begin{equation*}
	\eps^{\frac{ 2(m-\frac{5}{2}) }{ q - (m-\frac{5}{2}) }} c^{\frac{ 3(m-\frac{5}{2}) }{ q - (m-\frac{5}{2}) }} \norm{W(0)}_{H^m}^{-\frac{q}{ q - (m-\frac{5}{2}) }} \lesssim \norm{W(0)}_{H^m}^{-1},
\end{equation*}
so that we are again already done. Thus, we assume $\norm{W(0)}_{H^m} \lesssim \eps^{-1}(1+\k^{-3})^{-1}$, and in particular, assume that $T \geq \eps(1+\k^{-3})$. For $k_0 \geq 0$ to be determined, we control the blow-up norm by
\begin{equation*}
\begin{split}
	\Vert W \Vert_{L^1_t W^{1,\infty}_x} & \lesssim \Vert P_{\ge k_0} W \Vert_{L^1_t \dot{W}^{1,\infty}} + \Vert P_{\le k_0} W \Vert_{L^1_t W^{1,\infty}_x} \\
	& \lesssim \sum_{k\ge k_0} 2^{\frac{5}{2}k} \Vert P_k W \Vert_{L^1_t L^2_x} + T^{1-\frac{1}{q}} \sum_{k \leq k_0} \Vert P_{k} W \Vert_{L^q_t W^{1,\infty}_x} \\
	& \lesssim T 2^{(\frac{5}{2}-m)k_0} \Vert W \Vert_{L^\infty_t H^{m}_x} + T^{1-\frac{1}{q}} \sum_{k \leq k_0} 2^{\frac{3}{R}k} \Vert P_k W \Vert_{L^q_t W^{1,R}_x}.
\end{split}
\end{equation*}
The low frequency part can be analyzed in a similar way that was used in the high regularity case. Using \eqref{eq:linear_bd_wo_gradient_localized} with $s=1$,
\begin{align*}
	\sum_{k \leq k_0} 2^{\frac{3}{R}k} \Vert P_k e^{ it\mu\Lambda} \PP_{\mu\La} & W(0) \Vert_{L_t^q W_x^{1,R}} \lesssim \sum_{k \leq k_0} 2^{\frac{3}{R}k} 2^{k^+} (\langle 2^k \ka\rangle^3 \varep\k^{-3})^{\frac{1}{q}} \norm{P_k \PP_{\mu\La} W(0)}_{L^2} \\
	& \lesssim \eps^{\frac{1}{q}} \sum_{k \leq k_0} (2^{\frac{3}{2}k + k^+} + 2^{\frac{3}{R}k + k^+}\k^{-\frac{3}{q}}) \norm{P_k W(0)}_{L^2} \\
	& \lesssim \eps^{\frac{1}{q}} \sum_{k \leq k_0} (2^{\frac{3}{2}k + k^+ + (1-m)k^+} \norm{W(0)}_{H^{m-1}} + 2^{\frac{3}{R}k + k^+ + (1-m)k^+} \k^{-\frac{3}{q}} \norm{W(0)}_{H^{m-1}}) \\
	& \lesssim \eps^{\frac{1}{q}} (1+ \k^{-\frac{3}{q}}) 2^{(\frac{7}{2}-m)k_0} \norm{W(0)}_{H^{m-1}},
\end{align*}
and similarly
\begin{align*}
	\sum_{k \leq k_0} 2^{\frac{3}{R}k} \left\Vert P_k \int_0^t e^{ i(t-\tau)\mu\Lambda} \right.  \PP_{\mu\La} \mathcal{N}(W, W)(\tau)d\tau \bigg\Vert_{L_t^qW_x^{1,R}}
	& \lesssim \sum_{k \leq k_0} 2^{\frac{3}{R}k + k^+} (\langle 2^k \ka\rangle^3 \varep\k^{-3})^{\frac{1}{q}} \norm{\mathcal{N}(W, W)}_{L_t^1 L_x^2} \\
	& \lesssim \eps^{\frac{1}{q}} (1+ \k^{-\frac{3}{q}}) 2^{(\frac{7}{2}-m)k_0} \norm{\mathcal{N}(W, W)}_{L_t^1 H_x^{m-1}} \\
	& \lesssim \eps^{\frac{1}{q}} (1+ \k^{-\frac{3}{q}}) 2^{(\frac{7}{2}-m)k_0} \norm{W}_{L_t^1 W_x^{1,\infty}} \norm{W}_{L_t^\infty H_x^m}.
\end{align*}
Putting together, we obtain
\begin{equation*}
	\begin{split}
		\norm{W}_{L_t^1 W_x^{1,\infty}} \lesssim & \hspace{1mm} T (2^{k_0})^{(\frac{5}{2}-m)} \norm{W}_{L_t^\infty H^m_x} \\
		& + T^{1-\frac{1}{q}} \eps^{\frac{1}{q}} (1+\ka^{-\frac{3}{q}}) (2^{k_0})^{(\frac{7}{2}-m)} (\norm{W(0)}_{H^{m-1}} + \norm{W}_{L_t^1W_x^{1,\infty}} \norm{W}_{L_t^\infty H_x^m}).
	\end{split}
\end{equation*}
Therefore, choosing $2^{k_0} = (T\eps^{-1})^{\frac{1}{q}} (1+\k^{-\frac{3}{q}})^{-1}$, we get the bound
\begin{equation}\label{eq:low-reg-final-bound}
	\norm{W}_{L_t^1 W_x^{1,\infty}} \lesssim (T^{q-(m-\frac{5}{2})} \eps^{m-\frac{5}{2}} (1 + \ka^{-3})^{(m-\frac{5}{2})} )^{\frac{1}{q}} (\norm{W}_{L_t^\infty H^m_x} + \norm{W}_{L_t^1W_x^{1,\infty}} \norm{W}_{L_t^\infty H_x^m}).
\end{equation}
Hence, until the time given by the assumption \eqref{eq:main-assumption-low-reg}
\begin{equation}
	T = M \eps^{-\frac{ m-\frac{5}{2} }{ q - (m-\frac{5}{2}) }} \min\{1,\ka^{\frac{3(m-\frac{5}{2})}{ q - (m-\frac{5}{2}) } }\} \norm{W(0)}_{H^m}^{-\frac{q}{ q - (m-\frac{5}{2}) }},
\end{equation}
we obtain from \eqref{eq:low-reg-final-bound} and the energy estimates that
\begin{equation}
	\norm{W}_{L^\infty_t H^m_x} \leq 2 e^K \norm{W(0)}_{H^m},\qquad
	\norm{W}_{L_t^1W_x^{1,\infty}} \leq 1,
\end{equation}
so that the bootstrap will work until such time. This finishes the proof of Theorem \ref{main-theorem}.

\end{proof}

\begin{proof}[Proof of Theorem \ref{MainThmEC}]

Using the Strichartz estimates in \eqref{EC-Strichartz}, as well as the inhomogeneous counterpart, we can exploit the fact that nonlinearities of \eqref{EC} and \eqref{eq:CER2} have similar structure where they are both quadratic with one derivative, to derive an ``incompressible limit'' of Theorem \ref{main-theorem}. We see  \cite{KLT14b} that Duhamel's formula for \eqref{EC} is given by
\begin{equation}
	\mathbb{P}_\pm u(t) = e^{\pm it\omega_\eps} \mathbb{P}_\pm u_0 - \int_0^t e^{\pm i\omega_\eps (t-\tau)} \mathbb{P}_\pm (u\cdot\nabla)u(\tau) d\tau, \quad \mathcal{F}(\mathbb{P}_\pm u) = \frac{1}{2}\left( (I-\frac{\xi \xi^\top}{|\xi|^2})\hat{u} \pm i \frac{\xi}{|\xi|} \times \hat{u} \right).
\end{equation}
Note that $\mathbb{P}_\pm$ is a standard Calder\'on-Zygmund operator and is bounded on $L^p$, $1<p<\infty$. We take the energy estimates in \cite{KLT14b} for granted,
\begin{equation}
	\norm{u(t)}_{H^m} \leq \norm{u_0}_{H^m} \exp(K_m \norm{\nabla u(\tau)}_{L_\tau^1(0,t; L_x^\infty)}), \quad m \geq 0,
\end{equation}
and use the same strategy of bootstrapping $\norm{u}_{L_t^\infty H^m_x}$ and $\norm{u}_{L_t^1 W_x^{1,\infty}}$ together. When $m > \frac{7}{2}$, the blowup norm of the linear term is treated by
\begin{align*}
	\norm{e^{\pm it\omega_\eps} \mathbb{P}_\pm u_0}_{L_t^1 W_x^{1,\infty}}
	& \lesssim T^{1-\frac{1}{q}} \norm{e^{\pm it\omega_\eps} u_0}_{L_t^q W_x^{1+s,R}} \\
	& \lesssim T^{1-\frac{1}{q}} \left\{\sum_{k \in \mathbb{Z}} \left[ 2^{(1+s)k^+} \varep^{\frac{1}{q}} 2^{\frac{3}{q}k} \norm{P_k u_0}_{L^2}\right]^2 \right\}^{1/2} \\
	& \lesssim T^{1-\frac{1}{q}} \eps^{\frac{1}{q}} \norm{u_0}_{H^{1+s+\frac{3}{q}}}.
\end{align*}
The nonlinear term follows similarly:
\begin{align*}
	& \left\Vert \int_0^t e^{\pm i\omega_\eps (t-\tau)} \mathbb{P}_\pm (u\cdot\nabla)u(\tau) d\tau \right\Vert_{L_t^1 W_x^{1,\infty}} \\
	& \lesssim T^{1-\frac{1}{q}} \left\{ \sum_{k \in \mathbb{Z}} \left[ 2^{(1+s)k^+} \varepsilon^{\frac{1}{q}} 2^{\frac{3}{q}k} \Vert P_k (u\cdot\nabla)u \Vert_{L_t^1 L_x^2} \right]^2 \right\}^{1/2} \\
	& \lesssim T^{1-\frac{1}{q}} \eps^{\frac{1}{q}} \norm{(u\cdot\nabla)u}_{L_t^1 H_x^{1+s+\frac{3}{q}}} \\
	& \lesssim T^{1-\frac{1}{q}} \eps^{\frac{1}{q}} [ \norm{u}_{L_t^\infty H^{1+s+\frac{3}{q}}} \norm{\nabla u}_{L_t^1 L_x^\infty} + \norm{\nabla u}_{L_t^\infty H^{1+s+\frac{3}{q}}} \norm{u}_{L_t^1 L_x^\infty}].
\end{align*}
We set $s = m-2-\frac{3}{q}$, which gives us,
\begin{gather*}
	\norm{u}_{L_t^\infty H^m} \leq \norm{u_0}_{H^m} \exp(K_m \norm{u}_{L_t^1W_x^{1,\infty}}), \\
	\norm{u}_{L_t^1W_x^{1,\infty}} \leq C_{q,m} T^{1-\frac{1}{q}} \eps^{\frac{1}{q}} (\norm{u_0}_{H^m} + \norm{u}_{L_t^\infty  H^m} \norm{u}_{L_t^1W_x^{1,\infty}}).
\end{gather*}
By the same argument as before, this shows that the solution to \eqref{EC} exists at least up to
\begin{equation}
	T \gtrsim \eps^{-\frac{1}{q-1}} \norm{u_0}_{H^m}^{-\frac{q}{q-1}}.
\end{equation}
The case $m \in (\frac{5}{2}, \frac{7}{2})$ can be shown similarly.
\end{proof}

\subsection*{Acknowledgments}
H.\ Ko and B.\ Pausader are supported in part by NSF grant DMS-2154162 and DMS-2452275. R.\ Takada was supported in part by JSPS KAKENHI Grant Numbers
JP22K03388, JP23K20805 and JP23K20222. K.\ Widmayer gratefully acknowledges support of the SNSF through grant PCEFP2\_203059. The authors are very thankful to the anonymous referees for constructive comments which helped improve the clarity of the paper.

\appendix

\section{$L^p$ equivalence of change of variables}\label{appdx-Lp-equiv}

Here, we finish the proof of Proposition \ref{prop:Lp-equivalence} by showing that the components of the matrix in \eqref{UtoV-transformation-matrix} are all $L^p$ Fourier multipliers. Together with Lemma \ref{transformation-in-angles}, this will show that the norms are equivalent. The following elementary facts will be often used to ensure that the equivalence does not depend on $\kappa$.
\begin{itemize}
\item If $m(\xi)$ is a H\"ormander-Mikhlin multiplier, then so is $m(a\xi)$ for any $a\in\mathbb{R}$, and the boundedness is independent of $a$.
\item If $f(\xi)$ is a smooth function in a compact set $K$, then it is a H\"ormander-Mikhlin multiplier. Furthermore, the boundedness is indepndent of $a > 0$ for $f(a\xi)$ in $a^{-1} K$.
\end{itemize}
\vspace{4mm}

\noindent
A.1. $b_{\rho,1}$ and $b_{\rho,2}$.
For the matrix in \eqref{From-U12-to-rho-beta}, it's easier to investigate $b_{\rho,1}$ and $b_{\rho,2}$. We decompose $b_{\rho,1}$ into
\begin{align*}
b_{\rho,1} & = b_{\rho,1} \psi(2\k d_1) \phi(\frac{\xi_3-\k^{-1}}{r}) + b_{\rho,1} \psi(2\k d_1) (1-\phi(\frac{\xi_3-\k^{-1}}{r})) \\
& + b_{\rho,1} \psi(2\k d_2) \phi(\frac{\xi_3+\k^{-1}}{r}) + b_{\rho,1} \psi(2\k d_2) (1-\phi(\frac{\xi_3+\k^{-1}}{r})) + b_{\rho,1} (1-\psi(2\k d_1)) (1-\psi(2\k d_2)),
\end{align*}
where $\psi$ is the usual smooth bump function with $\text{supp } \psi \subset [-1,1]$ and $\phi \in C^{\infty}(\mathbb{R})$ is a smoothed version of Heaviside function such that $\phi(x) = 0$ for $x \leq -\frac{1}{10}$, $\phi(x) = 1$ for $x \geq \frac{1}{10}$ and monotone. We also decompose $b_{\rho,2}$ using the multipliers localized in the same 5 different regions. \\

\noindent
\textcircled{1} $d_1 \leq \frac{\k^{-1}}{2}, \xi_3 - \k^{-1} \geq -\frac{1}{10}r$ : It will be convenient to write $\mathbf{x} = (x,y,z) = (\xi_1,\xi_2,\xi_3 - \k^{-1})$ in this region, which makes $d_1 = |\mathbf{x}|$. Using \eqref{evalue-squares}, $b_{\rho,1}$ can be written as
\begin{equation*}
b_{\rho,1} = \frac{r}{\{ [z + \k \frac{r^2+z^2}{2} + \sqrt{(z + \k \frac{r^2+z^2}{2})^2 + r^2}]^2 + r^2 \}^{1/2}}.
\end{equation*}
By expanding the squares in the denominator and collecting terms of same order,
\begin{align*}
b_{\rho,1}^2 & = \frac{r^2}{2 d_1} \frac{1}{d_1+z + d_1(\k z + \frac{\k^2}{4}d_1^2) + z(\sqrt{1+\k z + \frac{\k^2}{4}d_1^2} - 1) + \frac{\k}{2}d_1^2 \sqrt{1+\k z + \frac{\k^2}{4}d_1^2}} \\
& = \frac{r^2}{2 d_1(d_1+z)} \frac{1}{1 + p_1(\k d_1,\k z)}.
\end{align*}
$p_1$ is continuous at the origin as a sum of product of $0$th order rational function and a smooth function, and is greater than $-1$ under the assumed localization.(crude estimate gives a lower bound $-\frac{11}{100}$) Therefore, the problem reduces to proving that $\frac{r}{\sqrt{d_1(d_1+z)}} = \frac{r}{\sqrt{|\mathbf{x}|(|\mathbf{x}|+z)}}$ is a $L^p$ Fourier multiplier. By the Marcinkiewicz multiplier theorem,
\begin{align*}
\nabla_{x,y} \frac{r}{\sqrt{|\mathbf{x}|(|\mathbf{x}|+z)}} & = \frac{1}{r\sqrt{|\mathbf{x}|(|\mathbf{x}|+z)}} \begin{pmatrix} x \\ y \end{pmatrix}
- \frac{r}{|\mathbf{x}|\sqrt{|\mathbf{x}|(|\mathbf{x}|+z)}} \frac{1}{|\mathbf{x}|} \begin{pmatrix} x \\ y \end{pmatrix}
- \frac{r}{(|\mathbf{x}|+z)\sqrt{|\mathbf{x}|(|\mathbf{x}|+z)}} \frac{1}{|\mathbf{x}|} \begin{pmatrix} x \\ y \end{pmatrix} \\
& \lesssim \frac{1}{|\mathbf{x}|}, \\
\partial_z \frac{r}{\sqrt{|\mathbf{x}|(|\mathbf{x}|+z)}} & = -\frac{r}{|\mathbf{x}|\sqrt{|\mathbf{x}|(|\mathbf{x}|+z)}} \frac{z}{|\mathbf{x}|}
- \frac{r}{(|\mathbf{x}|+z)\sqrt{|\mathbf{x}|(|\mathbf{x}|+z)}} \frac{|\mathbf{x}|+z}{|\mathbf{x}|} \lesssim \frac{1}{|z|}, \\
\partial_{xy} \frac{r}{\sqrt{|\mathbf{x}|(|\mathbf{x}|+z)}} & =
-\frac{xy}{r^3} \frac{1}{\sqrt{|\mathbf{x}|(|\mathbf{x}|+z)}}
- \frac{2xy}{r|\mathbf{x}|}\left( \frac{1}{|\mathbf{x}|\sqrt{|\mathbf{x}|(|\mathbf{x}|+z)}} + \frac{1}{(|\mathbf{x}|+z)\sqrt{|\mathbf{x}|(|\mathbf{x}|+z)}} \right) \\
& + \frac{xy}{|\mathbf{x}|^3} \left( \frac{r}{|\mathbf{x}|\sqrt{|\mathbf{x}|(|\mathbf{x}|+z)}} + \frac{r}{(|\mathbf{x}|+z)\sqrt{|\mathbf{x}|(|\mathbf{x}|+z)}} \right) \\
& - \frac{x}{|\mathbf{x}|}\left( \frac{1}{|\mathbf{x}|} + \frac{1}{|\mathbf{x}|+z} \partial_y\left( \frac{r}{\sqrt{|\mathbf{x}|(|\mathbf{x}|+z)}} \right) \right) \\
& + \frac{x}{|\mathbf{x}|}\left( \frac{y}{|\mathbf{x}|^3} + \frac{y}{|\mathbf{x}|(|\mathbf{x}|+z)^2} \right) \frac{r}{\sqrt{|\mathbf{x}|(|\mathbf{x}|+z)}} \quad \lesssim \frac{1}{|xy|} \\
\partial_z \nabla_{x,y} \frac{r}{\sqrt{|\mathbf{x}|(|\mathbf{x}|+z)}} & = \partial_z \left( \frac{r}{\sqrt{|\mathbf{x}|(|\mathbf{x}|+z)}} \right) \frac{1}{r^2} \begin{pmatrix} x \\ y \end{pmatrix} \\
& + \left( \frac{z}{|\mathbf{x}|^3} + \frac{1}{|\mathbf{x}|(|\mathbf{x}|+z)} \right) \frac{r}{\sqrt{|\mathbf{x}|(|\mathbf{x}|+z)}} \frac{1}{|\mathbf{x}|} \begin{pmatrix} x \\ y \end{pmatrix} \\
& - \left( \frac{1}{|\mathbf{x}|} + \frac{1}{|\mathbf{x}|+z} \right) \partial_z \left( \frac{r}{\sqrt{|\mathbf{x}|(|\mathbf{x}|+z)}} \right) \frac{1}{|\mathbf{x}|} \begin{pmatrix} x \\ y \end{pmatrix} \quad \lesssim \frac{1}{r|z|} \\
\end{align*}
\begin{align*}
\partial_{xyz} \frac{r}{\sqrt{|\mathbf{x}|(|\mathbf{x}|+z)}} = & -\frac{2xy}{r^4} \partial_z \left( \frac{r}{\sqrt{|\mathbf{x}|(|\mathbf{x}|+z)}} \right)
+ \frac{x}{r^2} \partial_{yz} \left( \frac{r}{\sqrt{|\mathbf{x}|(|\mathbf{x}|+z)}} \right) \\
& - \frac{xy}{|\mathbf{x}|^3} \left( \frac{z}{|\mathbf{x}|^3} + \frac{1}{|\mathbf{x}|(|\mathbf{x}|+z)} \right) \frac{r}{\sqrt{|\mathbf{x}|(|\mathbf{x}|+z)}} \\
& - \frac{x}{|\mathbf{x}|} \left( \frac{3yz}{|\mathbf{x}|^5} + \frac{y}{|\mathbf{x}|^3 (|\mathbf{x}|+z)} + \frac{y}{|\mathbf{x}|^2 (|\mathbf{x}|+z)^2} \right) \frac{r}{\sqrt{|\mathbf{x}|(|\mathbf{x}|+z)}} \\
& + \frac{x}{|\mathbf{x}|} \left( \frac{z}{|\mathbf{x}|^3} + \frac{1}{|\mathbf{x}|(|\mathbf{x}|+z)} \right) \partial_y \left(\frac{r}{\sqrt{|\mathbf{x}|(|\mathbf{x}|+z)}}\right) \\
& + \frac{xy}{|\mathbf{x}|^3} \left( \frac{1}{|\mathbf{x}|} + \frac{1}{|\mathbf{x}|+z} \right) \partial_z \left( \frac{r}{\sqrt{|\mathbf{x}|(|\mathbf{x}|+z)}} \right) \\
& + \frac{x}{|\mathbf{x}|}\left( \frac{y}{|\mathbf{x}|^3} + \frac{y}{|\mathbf{x}|(|\mathbf{x}+z)^2} \right) \partial_z \left(\frac{r}{\sqrt{|\mathbf{x}|(|\mathbf{x}|+z)}}\right) \\
& - \frac{x}{|\mathbf{x}|}\left( \frac{1}{|\mathbf{x}|} + \frac{1}{|\mathbf{x}|+z} \right) \partial_{yz} \left( \frac{r}{\sqrt{|\mathbf{x}|(|\mathbf{x}|+z)}} \right) \\
\lesssim & \hspace{1mm} \frac{1}{|xyz|},
\end{align*}
gives us the desired conclusion.\\

To show that $b_{\rho,2}$ is also a $L^p$ bounded operator, we use the identity
\begin{equation*}
\left[ \k^2 \left( \Sigma^2 - \k^{-2} \right)^2 + r^2 \right] \cdot \left[ \k^2 \left( \Omega^2 - \k^{-2} \right)^2 + r^2 \right] = (d_1 d_2 \k r)^2.
\end{equation*}
Then,
\begin{align*}
b_{\rho,2} = \frac{1}{\k d_1 d_2} \left[ \k^2 \left( \Sigma^2 - \k^{-2} \right)^2 + r^2 \right]^{1/2}
= \frac{\sqrt{2}}{\k d_2} \left[1 + \k z + \frac{\k^2}{4}d_1^2 + (\frac{z}{d_1} + \frac{\k}{2}d_1) \sqrt{1 + \k z + \frac{\k^2}{4}d_1^2} \right]^{1/2}.
\end{align*}
Both $\k d_2$ and the term inside the square root are smooth and nonzero in the area under consideration.(crude estimate gives a lower bound $\frac{\sqrt{7}}{4} - \frac{7}{20}$) Hence, it satisfies the H\"ormander-Mihklin condition, and therefore a $L^p$ Fourier multiplier. Note that the right-hand side can be written as a function in $\k\mathbf{x}$ so that the boundedness is uniform. \\

\noindent
\textcircled{2} $d_1 \leq \frac{\k^{-1}}{2}, \xi_3 - \k^{-1} \leq \frac{1}{10}r$ : The formula for $b_{\rho,2}$ gives
\begin{align*}
b_{\rho,2}^2 & = \frac{r^2}{2 d_1} \frac{1}{d_1 - z + d_1(\k z + \frac{\k^2}{4}d_1) - z(\sqrt{1+\k z + \frac{\k^2}{4}d_1^2} - 1) - \frac{\k^2}{2}d_1^2 \sqrt{1+\k z + \frac{\k^2}{4}d_1^2}} \\
& = \frac{r^2}{2 d_1(d_1-z)} \frac{1}{1 + p_2(\k d_1,\k z)}.
\end{align*}
Again, $p_2$ is always greater than $-1$, and $(1+p_2(s,t))^{-1}$ has a removable singularity at the origin, making it a bounded Fourier multiplier. The different region $z \leq \frac{1}{10}r$ instead of $z \geq -\frac{1}{10}r$ makes $b_{\rho,2}$ the favorable one to apply this formula, and it's straightforward to adjust the computations in \textcircled{1} to see that $\frac{r}{\sqrt{d_1(d_1-z)}}$ is a $L^p$ Fourier multiplier in this region.

Boundedness of $b_{\rho,1}$ follows the similar reasoning how we deduced $b_{\rho,2}$ is bounded in region \textcircled{1}. By the same identity $\left[ \k^2 \left( \Sigma^2 - \k^{-2} \right)^2 + r^2 \right] \cdot \left[ \k^2 \left( \Omega^2 - \k^{-2} \right)^2 + r^2 \right] = (d_1 d_2 \k r)^2$,
\begin{equation*}
b_{\rho,1} = \frac{1}{\k d_1 d_2} \left[ \k^2 \left( \Omega^2 - \k^{-2} \right)^2 + r^2 \right]^{1/2}
= \frac{\sqrt{2}}{\k d_2} \left[ 1 + \k z + \frac{\k^2}{4}d_1^2 - (\frac{z}{d_1} + \frac{\k}{2}d_1) \sqrt{1 + \k z + \frac{\k^2}{4}d_1^2} \right]^{1/2}.
\end{equation*}
Again, the term inside the square root is nonzero under the current localization, and hence the H\"ormander-Mihklin theorem applies directly. \\

\noindent
\textcircled{3}, \textcircled{4} : These are the cases when $d_2 \leq \frac{\k^{-1}}{2}$ and $\xi_3 + \k^{-1} \geq -\frac{1}{10}r$ or $\xi_3 + \k^{-1} \leq \frac{1}{10}r$, which are `upside-down' cases of \textcircled{1} and \textcircled{2}, respectively. Exchanging the role of $d_1$ and $d_2$, one gets the same results. \\

\noindent
\textcircled{5} $d_1, d_2 \geq \frac{\k^{-1}}{4}$ : Firstly, in the intersection with $|\xi| \leq 2\k^{-1}$, the multipliers have no singularity, so the multipliers and their derivatives are bounded by some constants simply by compactness. Hence, we assume $|\xi| \geq 2\k^{-1}$ instead. By writing $b_{\rho,1}, b_{\rho,2}$ as functions in $\frac{\xi}{|\xi|^2}$, we can invert the region into a compact area and proceed similarly. To illustrate, if we translate the origin to $(0,0,\k^{-1})$ again,
\begin{equation}\label{multiplier-near-infinity}
b_{\rho,1} = \frac{d_1}{d_2} \left[\frac{2 \sqrt{\frac{1}{4} + \frac{z}{\k d_1^2} + \frac{1}{\k^2 d_1^2}} }{\sqrt{\frac{1}{4} + \frac{z}{\k d_1^2} + \frac{1}{\k^2 d_1^2}} + \frac{1}{2} + \frac{z}{\k d_1^2} }\right]^{1/2} \frac{r}{\k d_1^2}
\end{equation}
can be written as $\frac{|\mathbf{x}|}{d_2} F(\frac{\mathbf{x}}{\k|\mathbf{x}|^2})$ for some smooth function $F$, if we put $\mathbf{x} = \xi - \k^{-1}e_3$ again for the translated variable. Since $\frac{\mathbf{x}}{\k|\mathbf{x}|^2} = f(\k\mathbf{x})$ and $f(\mathbf{y}) = \frac{\mathbf{y}}{|\mathbf{y}|^2}$ satisfies the H\"ormander-Mihklin condition in the exterior domain $B_2^c$, so does $F(\frac{\mathbf{x}}{\k|\mathbf{x}|^2})$ in $(B_{2\k^{-1}})^c$. Hence, we only need to check $\frac{|\mathbf{x}|}{d_2}$.
\begin{align*}
\nabla \frac{|\mathbf{x}|}{d_2} & = \frac{(z+2\k^{-1})^2}{d_2^3 |\mathbf{x}|}\mathbf{x} - 2\k^{-1}\frac{|\mathbf{x}|}{d_2^3}e_3 \quad \lesssim \frac{1}{|\mathbf{x}|}, \\
\nabla^2 \frac{|\mathbf{x}|}{d_2} & = 2\frac{z+2\k^{-1}}{d_2^3 |\mathbf{x}|} \mathbf{x}
-3\frac{(z+2\k^{-1})^2}{d_2^5 |\mathbf{x}|} \mathbf{x} (\mathbf{x} + 2k^{-1})^\intercal
- \frac{(z+2\k^{-1})^2}{d_2^3 |\mathbf{x}|^3} \mathbf{x} \mathbf{x}^\intercal \\
& + \frac{(z+2\k^{-1})^2}{d_2^3 |\mathbf{x}|} \text{Id}
+ 6\k^{-1} \frac{|\mathbf{x}|}{d_2^5} e_3 (x,y,z+2\k^{-1})^\intercal 
- \frac{2\k^{-1}}{d_2^3 |\mathbf{x}|} e_3 \mathbf{x}^{\intercal} \quad \lesssim \frac{1}{|\mathbf{x}|^2}
\end{align*}
hold on $d_1, d_2 \geq \frac{\k^{-1}}{4}$ since $d_1 \sim d_2$. Therefore, by the H\"ormander-Mihklin theorem, $\frac{d_1}{d_2}$ is a $L^p$ Fourier multiplier, and this finishes the proof.\\

\noindent
A.2. $b_{\gamma,1}$ and $b_{\gamma,2}$
For the matrix in \eqref{From-U12-to-alpha-gamma}, it's easier to use formulas for $b_{\gamma,1}$ and $b_{\gamma,2}$ because $b_{\gamma,1} = -\k\Sigma b_{\rho,1}$ and $b_{\gamma,2} = -i\k\Omega b_{\rho,2}$. This time, we decompose them into
\begin{gather*}
b_{\gamma,1} = b_{\gamma,1} \psi(2\k d_1) + b_{\gamma,1} \psi(2\k d_2) + b_{\gamma,1} (1-\psi(2\k d_1)) (1-\psi(2\k d_2)), \\
b_{\gamma,2} = b_{\gamma,2} \phi(2\k\xi_3) + b_{\gamma,2} (1 - \phi(2\k\xi_3)).
\end{gather*}
However, $b_{\gamma,1}$ and $b_{\gamma,2}$ are not as symmetric as $b_{\rho,1}$ and $b_{\rho,2}$. Hence, we investigate these two separately.\\

\noindent
\textcircled{1} $L^p$ boundedness of $b_{\gamma,2}$ : It's enough to prove that $\k\Omega$ is a $L^p$ Fourier multiplier. From the definition,
\begin{equation*}
\k\Omega = \frac{\k}{2}(d_2 - d_1) = \frac{2\xi_3}{d_1 + d_2}.
\end{equation*}
For $b_{\gamma,2} \phi(2\k\xi_3)$, we use $\mathbf{x} = \xi - \k^{-1}e_3$ again to center at the singularity $(0,0,\k^{-1})$. Then,
\begin{align*}
\nabla_{x,y} \k\Omega & = - \frac{2(z+\k^{-1})}{d_1 d_2(d_1+d_2)} \begin{pmatrix} x \\ y \end{pmatrix} \quad \lesssim \frac{1}{d_1} = \frac{1}{|\mathbf{x}|}, \\
\partial_z \k\Omega & = -\frac{2(z+\k^{-1})^2}{d_1 d_2 (d_1+d_2)} + \frac{d_1+d_2}{2d_1 d_2} \quad \lesssim \frac{1}{|z|}, \\
\partial_{xy} \k\Omega & = \frac{d_1^2 + d_1 d_2 + d_2^2}{d_1^3 d_2^3} \frac{2xy(z+\k^{-1})}{d_1 + d_2} \quad \lesssim \frac{1}{d_1 d_2} \leq \frac{1}{|xy|}, \\
\partial_z \nabla_{x,y} \k\Omega & = \frac{d_1^2 + d_1 d_2 + d_2^2}{d_1^3 d_2^3} \frac{2(z+\k^{-1})^2}{d_1 + d_2} \begin{pmatrix} x \\ y \end{pmatrix} - \frac{d_1^3 + d_2^3}{d_1^3 d_2^3} \begin{pmatrix} x \\ y \end{pmatrix} \quad \lesssim \frac{1}{d_1^2}, \\
\partial_{xyz} \k\Omega & = \frac{d_1^4 + d_1^3 d_2 + d_1^2 d_2^2 + d_1 d_2^3 + d_2^4}{d_1^5 d_2^5} \frac{6xy(z+\k^{-1})^2}{d_1 + d_2} + \frac{3xy(d_1^5 + d_2^5)}{2d_1^5 d_2^5} \lesssim \frac{1}{d_1^3},
\end{align*}
and hence the Marcinkiewicz multiplier theorem can be applied. For $ b_{\gamma,2} (1 - \phi(2\k\xi_3))$, by setting $z = \xi_3 + \k^{-1}$, the derivatives are of the same form in $z-\k^{-1}$ instead of $z+\k^{-1}$, and one can apply the Marcinkiewicz theorem again. \\

\noindent
\textcircled{2} $b_{\gamma,1}$ when $d_1, d_2 \geq \frac{\k^{-1}}{4}$ : Since $\k\Sigma$ itself is unbounded, $b_{\gamma,1}$ cannot be treated in the same way as $b_{\gamma,2}$. We first look at the boundedness of the piece away from the singularities. Using the formula \eqref{multiplier-near-infinity} and the variable $\mathbf{x} = \xi - \k^{-1}e_3$ again,
\begin{align*}
-b_{\gamma, 1} & = b_{\rho,1} \k\Sigma
= \frac{d_1}{d_2} \left[\frac{2 \sqrt{\frac{1}{4} + \frac{z}{\k d_1^2} + \frac{1}{\k^2 d_1^2}} }{\sqrt{\frac{1}{4} + \frac{z}{\k d_1^2} + \frac{1}{\k^2 d_1^2}} + \frac{1}{2} + \frac{z}{\k d_1^2} }\right]^{1/2} \frac{r}{\k d_1^2} \cdot \k\frac{d_1+d_2}{2} \\
& = \left[\frac{2\sqrt{\frac{1}{4} + \frac{z}{\k d_1^2} + \frac{1}{\k^2 d_1^2}} }{\sqrt{\frac{1}{4} + \frac{z}{\k d_1^2} + \frac{1}{\k^2 d_1^2}} + \frac{1}{2} + \frac{z}{\k d_1^2} }\right]^{1/2} \frac{r}{2d_1} \left(\frac{d_1}{d_2} +1\right).
\end{align*}
We have already seen that the term with square root and $\frac{d_1}{d_2}$ are $L^p$ Fourier multipliers. Hence, it's enough to prove the same for $\frac{r}{d_1} = \frac{r}{|\mathbf{x}|}$. We can use the Marcinkiewicz multiplier theorem again:
\begin{align*}
\nabla_{x,y} \frac{r}{|\mathbf{x}|} & = \frac{1}{r|\mathbf{x}|} \begin{pmatrix} x \\ y \end{pmatrix} - \frac{r}{|\mathbf{x}|^3} \begin{pmatrix} x \\ y \end{pmatrix} \quad \lesssim \frac{1}{|\mathbf{x}|}, \\
\partial_z \frac{r}{|\mathbf{x}|} & = -\frac{rz}{|\mathbf{x}|^3} \quad \lesssim \frac{1}{|\mathbf{x}|}, \\
\partial_{xy} \frac{r}{|\mathbf{x}|} & = -\frac{xy}{r^3|\mathbf{x}|} - \frac{2xy}{r|\mathbf{x}|^3} + \frac{3rxy}{|\mathbf{x}|^5} \quad \lesssim \frac{1}{|xy|}, \\
\partial_z \nabla_{x,y} \frac{r}{|\mathbf{x}|} & = \frac{3rz}{|\mathbf{x}|^5} \begin{pmatrix} x \\ y \end{pmatrix} - \frac{z}{r|\mathbf{x}|^3} \begin{pmatrix} x \\ y \end{pmatrix} \quad \lesssim \frac{1}{r|\mathbf{x}|}, \\
\partial_{xyz} \frac{r}{|\mathbf{x}|} & = -\frac{15rxyz}{|\mathbf{x}|^7} + \frac{6xyz}{r|\mathbf{x}|^5} + \frac{xyz}{r^3|\mathbf{x}|^3} \quad \lesssim \frac{1}{|xyz|}.
\end{align*}
This finishes the proof. \\

\noindent
\textcircled{3} $b_{\gamma,1}$ when $d_1 \leq \frac{\k^{-1}}{2}$ or $d_2 \leq \frac{\k^{-1}}{2}$ : Except around infinity, $\k\Sigma$ is also bounded, so now we can treat it separately. $d_2$ is a smooth function in $d_1 \leq \frac{\k^{-1}}{2}$, and it is well known that $d_1 = |\mathbf{x}|$ is also a $L^p$ Fourier multiplier in a compact region. These make $\k\Sigma = \frac{\k}{2}(d_1 + d_2)$ a bounded Fourier multiplier in $d_1 \leq \frac{\k^{-1}}{2}$, independent of $\k$. By swapping the roles of $d_1$ and $d_2$ in $d_2 \leq \frac{\k^{-1}}{2}$, this proves that $b_{\gamma,1}$ is a $L^p$ Fourier multiplier.

\section{Optimality of the decay rate in Corollary \ref{cor:lin_decay}}\label{appdx-optimal-decay-rate}

\begin{proposition}\label{OptimalDecay}
For any function $f$ whose Fourier transform $g \in C_c^{\infty}(\mathbb{R}^3)$ is a nonnegative spherically symmetric function centered at $(0, \kappa^{-1})$, such that $g|_{B_{(4\ka)^{-1}}(0, \kappa^{-1})} \equiv 1$, and {\normalfont supp} $g \subset B_{(2\ka)^{-1}}(0, \kappa^{-1})$,
\begin{equation}
	\vert (e^{it\Lambda_{\kappa}} f)(0) \vert \gtrsim t^{-1} (c\ka^2)^{-1} - t^{-2} (c^2\ka)^{-1}
\end{equation}
holds for sufficiently large $t$ for $\Lambda \in \{\Sigma,\Omega\}$. The same is true when we assume $g$ is symmetric around $(0, -\kappa^{-1})$.
In particular,
$t^{-1}$ is the optimal decay rate for $||e^{it\Lambda}f||_\infty$ for general functions.
\end{proposition}

\begin{proof}
Without loss of generality, and for notational convenience, we may assume $g$ is spherically symmetric around the origin by considering $f(x)e^{i\frac{x_3}{\ka}}$, and abuse the notation to write $g(\rho)$ for $\hat{f}(\xi)$. Translation affects the dispersion relations to take the form
\begin{equation*}
	2\Sigma = c(\sqrt{\rho^2 + (2\ka^{-1})^2 + 4\ka^{-1}\rho\cos(\phi)} + \rho), \quad 2\Omega = c(\sqrt{\rho^2 + (2\ka^{-1})^2 + 4\ka^{-1}\rho\cos(\phi)} - \rho),
\end{equation*}
where $(\rho, \theta, \phi)$ are spherical coordinates of $\xi$.
\begin{align*}
	(e^{it\Sigma}f)(0) & = \int_0^{\infty} \int_0^\pi \int_0^{2\pi} e^{\frac{1}{2}itc\rho + \frac{1}{2}itc\sqrt{\rho^2 + (2\ka^{-1})^2 + 4\ka^{-1}\rho\cos(\phi)}} g(\rho) \rho^2 \sin(\phi) d\theta d\phi d\rho \\
	& = \int_0^\infty 2\pi g(\rho) \rho^2 e^{\frac{1}{2}itc\rho} \left( \int_0^\pi e^{\frac{1}{2}itc\sqrt{\rho^2 + (2\ka^{-1})^2 + 4\ka^{-1}\rho\cos(\phi)}} \sin(\phi) d\phi \right) d\rho
\end{align*}
The nested integral has a closed form. Putting $s = \cos(\phi)$ and $\Phi(s) = \sqrt{\rho^2 + (2\ka^{-1})^2 + 4\ka^{-1}\rho s}$,
\begin{align*}
	\int_{-1}^1 e^{\frac{1}{2}itc\Phi(s)} ds & = \int_{-1}^1 \frac{2}{itc\Phi'(s)} \frac{d}{ds} e^{\frac{1}{2}itc\Phi(s)} ds \\
	& = \left[ \frac{2}{itc\Phi'(s)} e^{\frac{1}{2}itc\Phi(s)} \right]^1_{-1} + \frac{2}{itc} \int_{-1}^1 e^{\frac{1}{2}itc\Phi(s)} \frac{\Phi''(s)}{[\Phi'(s)]^2} ds.
\end{align*}
From $\Phi'(s) = \frac{2\ka^{-1}\rho}{\Phi(s)}, \Phi''(s) = - \frac{4\ka^{-2}\rho^2}{\Phi(s)^3}$, $\frac{\Phi''(s)}{[\Phi'(s)]^2} = -\frac{1}{\Phi(s)} = -\frac{\Phi'(s)}{2\ka^{-1}\rho}$, hence, $\Phi(1) = \rho + 2\ka^{-1}, \Phi(-1) = 2\ka^{-1} - \rho$, and
\begin{align*}
	\int_0^\pi e^{\frac{1}{2}itc\sqrt{\rho^2 + (2\ka^{-1})^2 + 4\ka^{-1}\rho\cos(\phi)}} \sin(\phi) d\phi
	& = \frac{2}{itc}\left( \frac{2\ka^{-1} + \rho}{2\ka^{-1}\rho} e^{\frac{1}{2}itc(2\ka^{-1} + \rho)} - \frac{2\ka^{-1} - \rho}{2\ka^{-1}\rho} e^{\frac{1}{2}itc(2\ka^{-1} - \rho)} \right) \\
	& + \frac{2}{t^2 c^2 \ka^{-1} \rho} (e^{\frac{1}{2}itc(2\ka^{-1} + \rho)} - e^{\frac{1}{2}itc(2\ka^{-1} - \rho)}) \\
	& = \frac{e^{itc\ka^{-1}}}{itc\rho} \left[ 4i\sin(\frac{tc}{2}\rho) + 2\rho\ka\cos(\frac{tc}{2}\rho) \right] + \frac{4\ka e^{itc\ka^{-1}}}{t^2 c^2 \rho} i\sin(\frac{tc}{2}\rho).
\end{align*}
This gives us
\begin{equation*}
	(e^{it\Sigma}f)(0) = 2\pi e^{itc\ka^{-1}}\int_0^\infty g(\rho) \rho^2 e^{\frac{1}{2}itc\rho} \left( \frac{2}{tc\rho} \left[ 2\sin(\frac{tc}{2}\rho) - i\ka\rho\cos(\frac{tc}{2}\rho) \right] + \frac{4\ka}{t^2 c^2 \rho} i\sin(\frac{tc}{2}\rho) \right) d\rho.
\end{equation*}
We will only bound imaginary part of the integral from below which will be enough to prove the proposition. Contribution from the first term with $\sin(\frac{tc}{2}\rho)$ is
\begin{align*}
	\int_0^\infty g(\rho) \rho^2 \frac{4}{tc\rho} \sin^2(\frac{tc}{2}\rho) d\rho
	\geq \frac{4}{tc} \int_0^{(4\ka)^{-1}} g(\rho) \rho \sin^2(\frac{tc}{2}\rho) d\rho = \frac{4}{tc} \int_0^{(4\ka)^{-1}} \rho \sin^2(\frac{tc}{2}\rho) d\rho,
\end{align*}
using the assumption on $g$. Elementary calculations show us that
\begin{equation*}
	\int_0^a x \sin^2(x) dx = \frac{1}{4}a^2 - \frac{1}{4}a\sin(2a) + \frac{1}{4} \sin^2(a).
\end{equation*}
Hence, we get a lower bound by
\begin{equation*}
	\frac{4}{tc} \int_0^{(4\ka)^{-1}} \rho \sin^2(\frac{tc}{2}\rho) d\rho
	= \frac{1}{tc}\frac{1}{16\ka^2} - \frac{1}{(tc)^2} \frac{1}{2\ka}\sin(\frac{tc}{4\ka}) + \frac{1}{(tc)^3} 4\sin^2(\frac{tc}{8\ka}).
\end{equation*}
Similarly, the second term with $\cos(\frac{tc}{2}\rho)$ gives
\begin{equation*}
	-\int_0^\infty g(\rho) \rho^2 \frac{2\ka}{tc}\cos^2(\frac{tc}{2}\rho) d\rho \geq - \frac{\ka}{tc} \int_0^{(2\ka)^{-1}} g(\rho) \rho^2 \cos^2(\frac{tc}{2}\rho) d\rho
	\geq - \frac{2\ka}{tc} \int_0^{(2\ka)^{-1}} \rho^2 \cos^2(\frac{tc}{2}\rho) d\rho.
\end{equation*}
By another elementary integration
\begin{equation*}
	\int_0^a x^2 \cos^2(x) dx = \frac{1}{6}a^3 + \frac{1}{4}a^2\sin(2a) + \frac{1}{4} a\cos(2a) - \frac{1}{8} \sin(2a),
\end{equation*}
we get a lower bound
\begin{equation*}
	- \frac{2\ka}{tc} \int_0^{(2\ka)^{-1}} \rho^2 \cos^2(\frac{tc}{2}\rho) d\rho
	= -\frac{1}{tc} \frac{1}{24\ka^2} - \frac{1}{(tc)^2} \frac{1}{4\ka}\sin(\frac{tc}{2\ka}) - \frac{1}{(tc)^3} \cos(\frac{tc}{2\ka}) + \frac{2\ka}{(tc)^4} \sin(\frac{tc}{2\ka}).
\end{equation*}
The last term does not have a definite sign, but can be estimated in a crude way due to the presence of $t^2$ in the denominator:
\begin{equation*}
	\int_0^\infty g(\rho) \rho^2 \frac{4\ka}{(tc)^2 \rho} \sin(2tc\rho) d\rho
	\geq - \frac{4\ka}{t^2 c^2}\int_0^{(2\ka)^{-1}} \rho d\rho = - \frac{1}{(tc)^2} \frac{1}{2\ka}.
\end{equation*}
Therefore, we have the following.
\begin{equation*}
	\frac{1}{2\pi}\left\vert (e^{2it\Sigma} f)(0) \right\vert \geq \frac{1}{tc}\frac{1}{48\ka^2} - \frac{1}{t^2} \frac{1}{2c^2\ka}\left[ \sin(\frac{tc}{4\ka}) + \frac{1}{2}\sin(\frac{tc}{2\ka}) + 1 \right] + O(\frac{1}{(tc)^3}),
\end{equation*}
which is the desired estimate.

Same decay estimate holds for $\Omega$ with having $e^{-\frac{1}{2}itc\rho}$ instead of $e^{\frac{1}{2}itc\rho}$, and one can immediately observe this shouldn't affect the proof at all. That the functions symmetric around $(0, -\ka^{-1})$ should also behave similarly can be easily seen by swapping the roles of $d_1$ and $d_2$.
\end{proof}

\begin{remark}
Note that the leading decay rate $t^{-1} (c\ka^2)^{-1} = \varepsilon^{-2} c^{-3} t^{-1}$ is exactly the same decay rate in Corollary \ref{cor:lin_decay}. $\langle \ka 2^k \rangle^3$ is dropped since $g$ is localized where $\ka 2^k \leq 2$.
\end{remark}

\bibliographystyle{plain}
\bibliography{CER_Strichartz}

\end{document}